# FINITE SIZE SCALING FOR THE CORE OF LARGE RANDOM HYPERGRAPHS

BY AMIR DEMBO[1] AND ANDREA MONTANARI[2]

*Stanford University*

The (two) core of a hypergraph is the maximal collection of hyperedges within which no vertex appears only once. It is of importance in tasks such as efficiently solving a large linear system over GF[2], or iterative decoding of low-density parity-check codes used over the binary erasure channel. Similar structures emerge in a variety of NP-hard combinatorial optimization and decision problems, from vertex cover to satisfiability.

For a uniformly chosen random hypergraph of $m = n\rho$ vertices and $n$ hyperedges, each consisting of the same fixed number $l \geq 3$ of vertices, the size of the core exhibits for large $n$ a first-order phase transition, changing from $o(n)$ for $\rho > \rho_c$ to a positive fraction of $n$ for $\rho < \rho_c$, with a transition window size $\Theta(n^{-1/2})$ around $\rho_c > 0$. Analyzing the corresponding "leaf removal" algorithm, we determine the associated finite-size scaling behavior. In particular, if $\rho$ is inside the scaling window (more precisely, $\rho = \rho_c + rn^{-1/2}$), the probability of having a core of size $\Theta(n)$ has a limit strictly between 0 and 1, and a leading correction of order $\Theta(n^{-1/6})$. The correction admits a sharp characterization in terms of the distribution of a Brownian motion with quadratic shift, from which it inherits the scaling with $n$. This behavior is expected to be universal for a wide collection of combinatorial problems.

**1. Introduction.** The $k$-core of a nondirected graph $G$ is the unique subgraph obtained by recursively removing all vertices of degree less than $k$. In particular, the 2-core, hereafter called the *core* of $G$, is the maximal collection of edges having no vertex appearing in only one of them. With an abuse

Received January 2007; revised December 2007.
[1]Research supported in part by NSF Grants DMS-04-06042 and DMS-FRG-0244323.
[2]Research supported in part by the European Union under the IP EVERGROW.

of language we shall use the same term for the induced subgraph. The core of a hypergraph is analogously defined and plays an important role in the analysis of many combinatorial problems.

In the first of such applications, Karp and Sipser [24] (hereafter KS) considered the problem of finding the largest possible matching (i.e., vertex disjoint set of edges) in a graph $G$. They proposed an algorithm that recursively selects an edge $(i,j) \in G$ for which the vertex $i$ has degree 1. If no such edge exists, the algorithm declares a failure. Otherwise it includes it in the matching and removes it from the graph together with all the edges incident on $j$ (that cannot belong to the matching). Whenever the algorithm successfully matches all vertices, the resulting matching can be proved to have maximal size. KS analyze the performance of such an algorithm on uniformly random graphs with $N$ vertices and $M = \lfloor Nc/2 \rfloor$ edges as $N \to \infty$, using the ODE asymptotic approximation for random processes, based on [25] (cf. [2, 16] for recent contributions).

It is easy to realize that the algorithm is successful if and only if a properly constructed hypergraph $\widetilde{G}$ does not contain a core. The hypergraph $\widetilde{G}$ includes a node $\widetilde{e}$ for each edge $e$ in $G$, and a hyperedge $\widetilde{i}$ for each vertex $i$ of degree 2 or more in $G$. The hyperedge $\widetilde{i}$ is incident on $\widetilde{e}$ in $\widetilde{G}$ if and only if $e$ is incident on $i$ in $G$.

A more recent application is related to the XOR-SAT problem, a simplified version of satisfiability introduced in [11]. One is given a linear system over $m$ binary variables, composed of $n$ equations modulo 2, each involving exactly $l \geq 3$ variables. The authors of [12, 29] propose a simple "leaf removal" algorithm to solve such a linear system. The algorithm recursively selects a variable that appears in a single equation, and eliminates the corresponding equation from the system. In fact, such an equation can be eventually satisfied by properly setting the selected variable. If all the equations are removed by this procedure, a solution can be constructed by running the process backward and fixing along the way the selected variables.

A hypergraph $G$ is associated to the linear system by including a vertex for each variable, and a hyperedge for each equation. Hyperedge $e$ is incident on vertex $i$ if and only if the corresponding equation involves the $i$th variable with nonzero coefficient. It is easy to realize that the leaf removal algorithm is successful if and only if the corresponding hypergraph $G$ does not contain a core.

Uniformly random linear systems with $n$ equations and $m = \rho n$ variables are considered in [12, 29]. It is proved there that the algorithm is successful with high probability if $\rho$ is larger than a critical value $\rho_c$, and fails with high probability if $\rho < \rho_c$. See Figure 1 for an illustration of this phenomenon. Further, it is shown there that the structure of the set of solutions of the linear system changes dramatically at $\rho_c$, exhibiting a "clustering effect" when $\rho < \rho_c$.

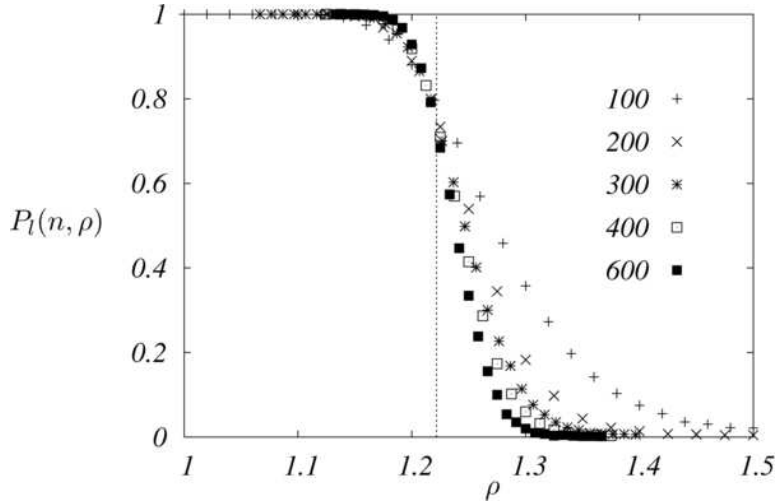

FIG. 1. *Probability that a random $l = 3$-hypergraph with $m$ vertices and $n = m/\rho$ hyperedges has a nonempty 2-core estimated numerically for $m = 100, \ldots, 600$. The vertical line corresponds to the asymptotic threshold $\rho_c \approx 1.2218$.*

The same "solution clustering" phenomenon has been conjectured for a variety of random combinatorial decision problems, on the basis of nonrigorous statistical mechanics calculations. The most studied among these problems is the random $K$-satisfiability, for which some indication of clustering is rigorously proved in [5, 15]. Several authors suggest that the solution clustering phenomenon is related to the poor performance of search algorithms on properly chosen ensembles of random instances. Still within random $K$-satisfiability, the performance of certain standard solution heuristics (such as the "pure-literal" rule) is also related to the appearance of properly defined cores (see [28]).

We conclude with the application to the analysis of low-density parity-check code ensembles, used for communication over the binary erasure channel, which is the most relevant motivation for our work. The decoding of a noisy message amounts in this case to finding the *unique* solution of a linear system over GF [2] (the solution exists by construction, but is not necessarily unique, in which case decoding fails). If the linear system includes an equation with only one variable, we thus determine the value of this variable, and substitute it throughout the system. Repeated recursively, this procedure either determines all the variables, thus yielding the unique solution of the system, or halts on a linear subsystem each of whose equations involves at least two variables. While such an algorithm is not optimal (when it halts, the resulting linear subsystem might still have a unique solution), it is the simplest instance of the widely used belief propagation decoding strategy, that has proved extremely successful. For example, on properly optimized code ensembles, this algorithm has been shown to achieve the theoretical limits for reliable communication, that is, Shannon's channel capacity (see [26]).

Once again, one can construct a hypergraph $G$ by associating a hyperedge to each variable, and a vertex to each equation (notice that this representation is "dual" with respect to the one used for XOR-SAT). Decoding is successful (it finds the unique solution) if and only if this hypergraph does not contain a core. For a "reasonable" code ensemble the probability of this event approaches 1 (resp. 0) when the noise level is smaller (larger) than a certain critical value. See [26] for an explicit characterization of the critical noise value via an application of the ODE method (based again on [25, 31]). Though this result has been successfully used for code design, it is often a poor approximation for the moderate code block-length (say, $n = 10^2$ to $10^5$) that is relevant in practice.

To overcome this problem, a finite-size scaling law is derived in [3], providing the probability of successfully decoding in the double limit of large size $n$, and noise level approaching the critical value. In [3] the authors also conjecture a "refined" law that describes how the finite-size scaling limit is approached, and demonstrate empirically that this refined scaling formula is very accurate already for short message lengths $n \approx 100$, opening the way to an efficient code design procedure (cf. [4]).

In this paper we resolve the conjecture by rigorously proving the refined scaling law. To simplify the exposition we focus on a specific choice of the random ensemble of equations, but our proof generalizes without much difficulty to a large variety of other cases, and in particular to all those mentioned in the conjecture of [3]. (In the coding language, the example we consider corresponds to LDPC ensembles with regular left and Poisson right degree; it also coincides with the random XOR-SAT ensemble introduced in [11] and treated in [12, 29].) In graph-theoretical terms, we determine the probability that a uniformly random[3] $l$-hypergraph (i.e., a hypergraph with hyperedges of size $l$) with $n$ hyperedges and $m = n\rho$ vertices has a nonempty core as $n$ grows and $\rho = \rho(n)$ approaches $\rho_c$. In the process of establishing the refined scaling law we gain much insight about the core of such random hypergraphs. For example, we determine the fluctuations in the size of the core at criticality (see Remark 2.6), and show that if the hypergraph is built one hyperedge at a time, then its core size jumps from zero to a positive fraction of $m$ at a random time $n_c$, the distribution of which we explicitly determine (cf. Remark 2.5).

---

[3]Indeed, we work with a properly defined "configurational" model (somewhat similar to the one introduced in [8]) to be defined in the next section.

Our proof strategy should apply without conceptual changes to other phase transitions within the same class, such as $k$-core percolation on random graphs (with $k \geq 3$), or the pure-literal rule threshold in random $k$-SAT (with $k \geq 3$; cf. [17]). Even beyond this family of closely related phenomena, the form of the refined scaling law (in particular, the scaling with $n$ of the scaling window and of the first correction) are likely to be quite universal. For instance, in [3] it has been empirically found to hold for iterative decoding of LDPC codes over general channels. Within statistical physics, several core phase transitions have been studied as special examples of "mean field dynamical glass transition" [33]. It is possible that the refined finite-size scaling law generalizes to this (quite large) class as well.

Finite-size scaling has been the object of several investigations in statistical physics and in combinatorics. Most of these studies estimate the size of the corresponding scaling window. That is, fixing a small value of $\varepsilon > 0$, they find the amount of change in some control parameter which moves the probability of a relevant event from $\varepsilon$ to $1 - \varepsilon$. A remarkably general result in this direction is the rigorous formulation of a "Harris criterion" in [10, 35]. Under mild assumptions, this implies that the scaling window has to be at least $\Omega(n^{-1/2})$ for a properly defined control parameter (e.g., the ratio $\rho$ of the number of nodes to hyper-edges in our problem). A more precise result has recently been obtained for the satisfiable-unsatisfiable phase transition for the random 2-SAT problem, yielding a window of size $\Theta(n^{-1/3})$ [9]. Note, however, that statistical physics arguments suggest that the phase transition we consider here is not from the same universality class as the satisfiable-unsatisfiable transition for the random 2-SAT problem.

In contrast with the preceding, we provide a much sharper characterization, yielding beyond the scaling window and the limiting scaling function, also the asymptotic form of corrections to this limit. In this respect, our work is closer in its level of precision to that for the scaling behavior in the emergence of the giant component in Erdős–Rényi random graphs (for more on the latter, see [22] and the references therein).

At the level of degenerate (or zero–one) fluid-limits, the asymptotic size of $k$-core of random graphs is determined by [31] via the ODE method. See also [28] for a general approach for deriving such results without recourse to ODE approximations (using instead a method analogous to the "density evolution" technique from coding theory).

Darling and Norris determine in [14] the asymptotic size of the 2-core of a random hypergraph which is the "dual" of the model we consider here. Indeed, the hyperedges in their model are of random, Poisson distributed, sizes, which allows for a particularly simple Markovian description of the recursive algorithm that constructs the core. Dealing as we do with random hypergraphs at the critical point, where the asymptotic core size exhibits a discontinuity, they describe the fluctuations around the deterministic limit via a certain linear SDE. In doing so, they heavily rely on the powerful theory of weak convergence, in particular in the context of convergence of Markov processes. For further results that are derived along the same line of reasoning; see [13, 18, 19].

In contrast, as we outline in the next section, the focus of this paper is on correction terms and rates of convergence. These are beyond the scope of weak convergence theory. In the context of our main result, Theorem 2.3, these only provide the limit, as $n \to \infty$ and $\rho_n$ near its critical value, of the probability that a uniformly chosen random hypergraph with $n$ hyperedges over $n\rho_n$ vertices has a nonempty 2-core.

The need to estimate correction terms is why many steps in our proof involve delicate coupling arguments, expanding and keeping track of the rate of decay of approximation errors (in terms of $n$). Our technique can be extended to provide rates of convergence (in the sup-norm) as $n$ grows, for distributions of inhomogeneous Markov chains on $\mathbb{R}^d$ whose transition kernels $W_{t,n}(x_{t+1} - x_t = y | x_t = x)$ are approximately (in $n$) linear in $x$, and "strongly elliptic" of uniformly bounded support with respect to $y$.

**2. Main result and outline of proof.** We consider hypergraphs with $n$ hyperedges over $m = \lfloor n\rho \rfloor$ vertices, $\rho > 0$. Each hyperedge is an ordered list of $l \geq 3$, *not necessarily distinct* vertices chosen independently and uniformly at random with replacement. We are interested in the probability $P_l(n, \rho)$ that a random hypergraph from this ensemble has a nonempty 2-core (i.e., the existence of a nonempty list of hyperedges such that, if a vertex appears in this list, then it does so at least twice).

In the next section we construct an inhomogeneous Markov chain $\{\vec{z}(\tau) = (z_1(\tau), z_2(\tau)), n \geq \tau \geq 0\}$, where $z_1(\tau)$ and $z_2(\tau)$ keep track, respectively, of the number of vertices of degree 1 and of degree at least 2 after $\tau$ steps of the decimation algorithm. As we show in Section 5, in the large $n$ limit, this chain is well approximated by a simpler chain with transition probabilities,

$$\widehat{\mathbb{P}}_{n,\rho}\{\vec{z}(\tau+1) = \vec{z} + (q_1 - q_0, -q_1) | \vec{z}(\tau) = \vec{z}\}$$
(2.1)
$$= \binom{l-1}{q_0 - 1, q_1, q_2} \mathfrak{p}_0^{q_0 - 1} \mathfrak{p}_1^{q_1} \mathfrak{p}_2^{q_2}.$$

For $\vec{x} = \vec{z}/n$, $\theta = \tau/n$,

$$(2.2) \quad \mathfrak{p}_0 = \frac{\max(x_1, 0)}{l(1-\theta)}, \qquad \mathfrak{p}_1 = \frac{x_2 \lambda^2}{l(1-\theta)(e^\lambda - 1 - \lambda)}, \qquad \mathfrak{p}_2 = \frac{x_2 \lambda}{l(1-\theta)},$$

where for $x_2 > 0$, we set $\lambda$ as the unique positive solution of

$$(2.3) \qquad f_1(\lambda) \equiv \frac{\lambda(e^\lambda - 1)}{e^\lambda - 1 - \lambda} = \frac{l(1-\theta) - \max(x_1, 0)}{x_2}$$

enforcing $\mathfrak{p}_0 + \mathfrak{p}_1 + \mathfrak{p}_2 = 1$, while for $x_2 = 0$, we instead set by continuity $\mathfrak{p}_1 = 0$ and $\mathfrak{p}_2 = 1 - \mathfrak{p}_0$.

Further, we show in Lemma 4.4 that $n^{-1}\vec{z}(0)$ converges to the nonrandom vector

$$(2.4) \qquad \vec{y}(0) = (le^{-l/\rho}, \rho(1 - e^{-l/\rho}) - le^{-l/\rho}).$$

Since the chain (2.1) has bounded increments, and the corresponding probabilities depend on the state only through the macroscopic variables $\vec{x}$ and $\theta$, it is not hard to verify that the scaled process $n^{-1}\vec{z}(\theta n)$ concentrates around the solution of the ODE

$$(2.5) \qquad \frac{d\vec{y}}{d\theta}(\theta) = \vec{F}(\vec{y}(\theta), \theta),$$

where $\vec{F}(\vec{x}, \theta) = (-1 + (l-1)(\mathfrak{p}_1 - \mathfrak{p}_0), -(l-1)\mathfrak{p}_1)$ is the mean of $\vec{z}(\tau+1) - \vec{z}(\tau)$ under the transitions of (2.1); see, for instance, [12, 26, 29]. The solution of this ODE will be denoted by $\vec{y}(\theta, \rho)$, often using the shorthand $\vec{y}(\theta)$ (where the fixed value of $\rho$ is clear from the context). From the solution, one finds that $y_1(\theta)$ remains strictly positive for all $\theta \in [0,1)$ if and only if $\rho > \rho_c$ [see (4.3)]. As shown in [26] this indicates that, with high probability, the algorithm successfully decimates the whole hypergraph without ever running out of degree 1 vertices if $\rho > \rho_c$. Vice versa, for $\rho < \rho_c$, the solution $\vec{y}(\theta)$ crosses the $y_1 = 0$ plane; this is shown to imply that the algorithm stops and returns a large core with high probability. In the critical case $\rho = \rho_c$, the solution $\vec{y}(\theta)$ touches the $y_1 = 0$ plane at the unique time $\theta = \theta_c \in (0,1)$ (see Proposition 4.2). The principal conclusion of Section 5 is that, near criticality, $P_l(n, \rho)$ can be estimated by the probability that $\vec{z}(\tau)$ is small in a neighborhood of $\tau = n\theta_c$. More precisely:

PROPOSITION 2.1. *Let $\beta \in (3/4, 1)$, $J_n = [n\theta_c - n^\beta, n\theta_c + n^\beta]$ and $|\rho - \rho_c| \leq n^{\beta'-1}$ with $\beta' < 2\beta - 1$. Then for $\varepsilon_n = A \log n$,*

$$(2.6) \quad \widehat{\mathbb{P}}_{n,\rho}\left\{\inf_{\tau \in J_n} z_1(\tau) \leq -\varepsilon_n\right\} - \delta_n \leq P_l(n, \rho) \leq \widehat{\mathbb{P}}_{n,\rho}\left\{\inf_{\tau \in J_n} z_1(\tau) \leq \varepsilon_n\right\} + \delta_n,$$

*where $\delta_n \equiv D n^{-1/2}(\log n)^2$.*

At the critical point (i.e., for $\rho = \rho_c$ and $\theta = \theta_c$) the solution of the ODE (2.5) is tangent to the $y_1 = 0$ plane and fluctuations in the $y_1$ direction determine whether a large core exists or not. Further, in a neighborhood of $\theta_c$, we have $y_1(\theta) \simeq \frac{1}{2}\widetilde{F}(\theta - \theta_c)^2$, for some $\widetilde{F} > 0$. In the same neighborhood, the contribution of fluctuations to the change of $z_1$ is approximately $\sqrt{\widetilde{G}n(\theta - \theta_c)}$, with $\widetilde{G} > 0$. Comparing these two contributions we see that the relevant scaling is $X_n(t) = n^{-1/3}[z_1(n\theta_c + n^{2/3}t) - z_1(n\theta_c)]$, which for large $n$ converges, by strong approximation, to $X(t) = \frac{1}{2}\widetilde{F}t^2 + \sqrt{\widetilde{G}}W(t)$, for a standard two-sided Brownian motion $W(t)$ (see Lemma 6.1 for a precise quantitative statement). Clearly,

$$(2.7) \qquad \widetilde{F} \equiv \frac{d^2 y_1}{d\theta^2}(\theta_c) = \frac{dF_1}{d\theta}(\vec{y}(\theta_c), \theta_c) = \frac{\partial F_1}{\partial \theta} + \frac{\partial F_1}{\partial y_2}F_2.$$

In the last expression we adopted the convention (to be followed hereafter) of omitting the arguments whenever they refer to the critical point $\theta = \theta_c$, $\vec{y} = \vec{y}(\theta_c)$ and the trajectory considered is the critical one, that is, $\rho = \rho_c$.

Fluctuations of $\vec{z}(n\theta_c)$ around $n\vec{y}(\theta_c)$ are accumulated in $n\theta_c$ stochastic steps, and are therefore of order $\sqrt{n}$. As shown in Section 6, the rescaled variable $(\vec{z}(n\theta) - n\vec{y}(\theta))/\sqrt{n}$ converges to a Gaussian random variable. Its covariance matrix $\mathbb{Q}(\theta, \rho) = \{Q_{ab}(\theta, \rho); 1 \leq a, b \leq 2\}$ is the symmetric positive definite solution of the ODE:

$$(2.8) \qquad \frac{d\mathbb{Q}(\theta)}{d\theta} = \mathbb{G}(\vec{y}(\theta), \theta) + \mathbb{A}(\vec{y}(\theta), \theta)\mathbb{Q}(\theta) + \mathbb{Q}(\theta)\mathbb{A}(\vec{y}(\theta), \theta)^\dagger,$$

where $\mathbb{A}(\vec{x}, \theta) \equiv \{A_{ab}(\vec{x}, \theta); 1 \leq a, b \leq 2\}$ for $A_{ab}(\vec{x}, \theta) = \partial_{x_b}F_a(\vec{x}, \theta)$, and $\mathbb{G}(\vec{x}, \theta)$ is the covariance of $\vec{z}(\tau+1) - \vec{z}(\tau)$ under the transitions (2.1), that is, the nonnegative definite symmetric matrix with entries

$$(2.9) \qquad \begin{cases} G_{11}(\vec{x}, \theta) = (l-1)[\mathfrak{p}_0 + \mathfrak{p}_1 - (\mathfrak{p}_0 - \mathfrak{p}_1)^2], \\ G_{12}(\vec{x}, \theta) = -(l-1)[\mathfrak{p}_0\mathfrak{p}_1 + \mathfrak{p}_1(1 - \mathfrak{p}_1)], \\ G_{22}(\vec{x}, \theta) = (l-1)\mathfrak{p}_1(1 - \mathfrak{p}_1). \end{cases}$$

Here again we use the convention $\mathbb{Q}(\theta) \equiv \mathbb{Q}(\theta, \rho)$ when the value of $\rho$ is clear from the context. The positive definite initial condition $\mathbb{Q}(0)$ for (2.8) is computed on the original graph ensemble, and given by

$$(2.10) \qquad \begin{cases} Q_{11}(0) = \frac{l}{\gamma}\gamma e^{-2\gamma}(e^\gamma - 1 + \gamma - \gamma^2), \\ Q_{12}(0) = -\frac{l}{\gamma}\gamma e^{-2\gamma}(e^\gamma - 1 - \gamma^2), \\ Q_{22}(0) = \frac{l}{\gamma}e^{-2\gamma}[(e^\gamma - 1) + \gamma(e^\gamma - 2) - \gamma^2(1 + \gamma)], \end{cases}$$

where $\gamma = l/\rho$ (see Section 4.2 for details).

The parameter describing the fluctuations of $z_1(n\theta) - z_1(n\theta_c)$ for $\theta$ near $\theta_c$ is simply $\widetilde{G} = G_{11}(\vec{y}(\theta_c), \theta_c)$. As we show in Section 6, this analysis allows us to approximate the probability that $\vec{z}(\tau)$ approaches the $z_1 = 0$ plane, by replacing $\{\vec{z}(\tau)\}$ by an appropriately constructed Gaussian process.

PROPOSITION 2.2. *Let $X(t) = \frac{1}{2}\widetilde{F}t^2 + \sqrt{\widetilde{G}}W(t)$ where $W(t)$ is a doubly infinite standard Brownian motion conditioned to $W(0) = 0$. Further, let* $\xi(r)$ be a normal random variable of mean $(\frac{\partial y_1}{\partial \rho})r$ and variance $Q_{11}$ (both evaluated at $\theta = \theta_c$ and $\rho = \rho_c$), which is independent of $W(t)$.

For some $\beta \in (3/4, 1)$, any $\eta < 5/26$, all $A > 0$, $r \in \mathbb{R}$ and $n$ large enough, if $\rho_n = \rho_c + r n^{-1/2}$ and $\varepsilon_n = A \log n$, then

$$(2.11) \quad \left| \widehat{\mathbb{P}}_{n, \rho_n} \left\{ \inf_{\tau \in J_n} z_1(\tau) \leq \pm \varepsilon_n \right\} - \mathbb{P}\left\{ n^{1/6} \xi + \inf_t X(t) \leq 0 \right\} \right| \leq n^{-\eta}.$$

We note in passing that within the scope of weak convergence, Aldous [1] pioneered the use of Brownian motion with quadratic drift [à la $X(t)$ of Proposition 2.2], to examine the near-critical behavior of the giant component in Erdös–Rényi random graphs, and his method was extended by Goldschmidt [19] to the giant set of identifiable vertices in Poisson random hypergraph models.

The key to the validity of Proposition 2.2 at the $o(n^{-1/6})$ level of accuracy relevant here, is the fact that within the critical time window $J_n$ the Markov chain of transition probabilities (2.1) is well approximated by the chain

$$(2.12) \quad \vec{z}'(\tau + 1) = \vec{z}'(\tau) + \widetilde{\mathbb{A}}_\tau (n^{-1} \vec{z}'(\tau) - \vec{y}(\tau/n)) + \Delta_\tau$$

with $\widetilde{\mathbb{A}}_\tau \equiv \mathbb{I}_{\tau < \tau_n} \mathbb{A}(\vec{y}(\tau/n, \rho), \tau/n)$ for $\tau_n \equiv \lfloor n\theta_c - n^\beta \rfloor$, and independent random variables $\{\Delta_\tau\}$ of mean $\vec{F}(\vec{y}(\tau/n), \tau/n)$ and covariance $\mathbb{G}(\vec{y}(\tau/n), \tau/n)$ (cf. Proposition 5.5). In particular, taking

$$(2.13) \quad \widetilde{\mathbb{B}}_\sigma^\tau \equiv \left(\mathbb{I} + \frac{1}{n}\widetilde{\mathbb{A}}_\tau\right) \cdot \left(\mathbb{I} + \frac{1}{n}\widetilde{\mathbb{A}}_{\tau-1}\right) \cdots \left(\mathbb{I} + \frac{1}{n}\widetilde{\mathbb{A}}_\sigma\right),$$

for integers $0 \leq \sigma \leq \tau$ (while $\widetilde{\mathbb{B}}_\sigma^\tau \equiv \mathbb{I}$ in case $\tau < \sigma$), we see that

$$(2.14) \quad \vec{z}'(\tau) = \widetilde{\mathbb{B}}_0^{\tau-1} \vec{z}'(0) + \sum_{\sigma=0}^{\tau-1} \widetilde{\mathbb{B}}_{\sigma+1}^{\tau-1} (\Delta_\tau - \widetilde{\mathbb{A}}_\sigma \vec{y}(\sigma/n))$$

is a sum of (bounded) independent random variables, hence of approximately normal distribution. Further, the mean and covariance of $\vec{z}'(\tau)$ are given by discretized versions of (2.5) and (2.8), hence are sufficiently close to the solutions $\vec{y}(\theta, \rho)$ and $\mathbb{Q}(\theta, \rho)$ of these ODEs (cf. Lemma 4.3).

Combining Propositions 2.1 and 2.2, we are now able to estimate the desired probability $P_l(n, \rho)$ in terms of the distribution of the global minimum of the process $\{X(t)\}$ (i.e., a Brownian motion plus a quadratic shift). The latter has been determined already in [20], yielding the following conclusion, which is our main result. Figure 2 illustrates the accuracy of the finite-size scaling expression proved below, by comparing it with numerical simulations.

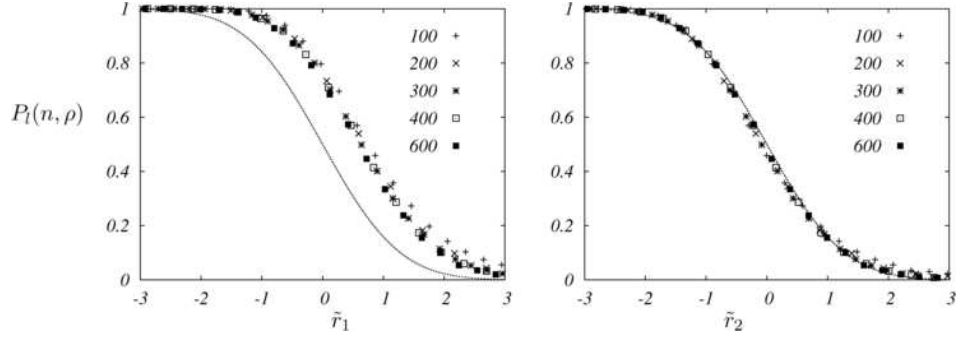

FIG. 2. *The numerical estimates for the core probabilities in Figure 1, plotted versus scaling variables $\tilde{r}_1$, $\tilde{r}_2$. On the left: $\tilde{r}_1 = \sqrt{n}(\rho - \rho_c)/\alpha_l$. On the right: $\tilde{r}_2 = \sqrt{n}(\rho - \rho_c - \delta_l n^{-2/3})/\alpha_l$ where $\delta_l = \alpha_l \beta_l \Omega$. According to Theorem 2.3, corrections to the asymptotic curve $\Phi(-\tilde{r})$ (dashed) are $\Theta(n^{-1/6})$ on the left, and $O(n^{-5/26+\varepsilon})$ on the right.*

THEOREM 2.3. *Let $l \geq 3$, and define $\alpha_l = \sqrt{Q_{11}}(\frac{\partial y_1}{\partial \rho})^{-1}$, $\beta_l = \frac{1}{\sqrt{Q_{11}}}\widetilde{G}^{2/3} \times \widetilde{F}^{-1/3}$, $\rho_n = \rho_c + rn^{-1/2}$. Then, for any $\eta < 5/26$*

$$(2.15) \quad P_l(n, \rho_n) = \Phi(-r/\alpha_l) + \beta_l \Omega \Phi'(-r/\alpha_l) n^{-1/6} + O(n^{-\eta}),$$

*where $\Phi(x)$ denotes the distribution function for a standard normal random variable, the finite constant $\Omega$ is given by the integral*

$$(2.16) \quad \Omega \equiv \int_0^\infty [1 - \mathcal{K}(z)^2]\, dz,$$

*where*

$$(2.17) \quad \mathcal{K}(z) \equiv \frac{1}{2} \int_{-\infty}^\infty \frac{\mathrm{Ai}(iy)\mathrm{Bi}(2^{1/3}z + iy) - \mathrm{Ai}(2^{1/3}z + iy)\mathrm{Bi}(iy)}{\mathrm{Ai}(iy)}\, dy,$$

*and $\mathrm{Ai}(\cdot)$, $\mathrm{Bi}(\cdot)$ are the Airy functions (as defined in [6], page 446).*

REMARK 2.4. The simulations in Figure 2 suggest that the approximation of $P_l(n, \rho_n)$ we provide in (2.15) is more accurate than the $O(n^{-5/26+\varepsilon})$ correction term suggests. Our proof shows that one cannot hope for a better, error estimate than $\Theta(n^{-1/3})$ as we neglect the second-order term in expanding $\Phi(-r/\alpha_l + Cn^{-1/6})$; see (2.18). We believe this is indeed the order of the next term in the expansion (2.15). Determining its form is an open problem.

REMARK 2.5. It is of interest to consider the (time) evolution of the core for the hypergraph process in which one hyperedge is added uniformly at random at each time step. In other words, $n$ increases with time, while the number of vertices $m$ is kept fixed. Let $S(n)$ be the corresponding (random) number of hyperedges in the core of the hypergraph at time $n$ and $n_c \equiv \min\{n : S(n) \geq 1\}$ the onset of a nonempty core. From Lemma 4.7 we have that for any $\rho > 0$ there exist $\kappa > 0$ and $C < \infty$ such that $\{S(n) : 0 \leq n \leq m/\rho\}$ intersects $[1, m\kappa]$ with probability at most $Cm^{1-l/2}$. Further, fixing $\rho < \rho_c$, the probability of an empty core, that is, $S(m/\rho) = 0$, decays (exponentially) in $m$. We thus deduce that for large $m$ most of the trajectories $\{S(n)\}$ jump from having no core to a linear (at least $m\kappa$) core size at the well-defined (random) critical edge number $n_c$. By the monotonicity of $S(n)$ we also know that $\mathbb{P}_m\{n_c \leq m/\rho\} = P_l(\rho, m/\rho)$. Therefore, Theorem 2.3 allows us to determine the asymptotic distribution of $n_c$. Indeed, expressing $n$ in terms of $m$ in (2.15) we get that for each fixed $x \in \mathbb{R}$,

$$\mathbb{P}\{n_c \leq m\rho_c^{-1} + m^{1/2}\rho_c^{-3/2}\alpha_l x\} = \Phi(x) + \beta_l \Omega \rho_c^{1/6} \Phi'(x) m^{-1/6} + O(m^{-\eta}),$$

whence we read off that $\hat{n}_c \equiv (n_c - m/\rho_c)/(\sqrt{m}\rho_c^{-3/2}\alpha_l) + \beta_l \Omega \rho_c^{1/6} m^{-1/6}$ converge in distribution to the standard normal law [and the corresponding distribution functions converge pointwise to $\Phi(x)$ at rate which is faster than $m^{-\eta}$ for any $\eta < 5/26$].

REMARK 2.6. Our techniques are applicable to many other properties of the core in the "scaling regime" $\rho_n = \rho_c + rn^{-1/2}$. For example, the distribution of the number of hyperedges $S$ in the core can be derived from the approximation of the trajectory of the decimation algorithm. Namely, as shown in Section 6, for such $\rho_n$, near the critical time $z_1(t) \simeq \sqrt{n}\xi(r) + X_n(t)$ for $\xi(r)$ and $X_n(t) \equiv n^{1/3}X(n^{-2/3}(t - n\theta_c))$ as in Proposition 2.2. With $\mathbb{E}X_n(t) = \frac{\widetilde{F}}{2n}(t - n\theta_c)^2$, upon noting that $n - S = \min\{t : z_1(t) = 0\}$, we obtain that, conditional to the existence of a nonempty core, $(S - n(1 - \theta_c))/n^{3/4}$ converges in distribution to $(4Q_{11}/\widetilde{F}^2)^{1/4}Z$ with $Z$ a nondegenerate random variable. Indeed, at the relevant time window $n\theta_c \pm O(n^{3/4})$ the contribution of $X_n(\cdot) - \mathbb{E}X_n(\cdot)$ to the fluctuations of $S$ is negligible in comparison with that of $\sqrt{n}\xi(r)$. So, more precisely, based on the explicit distribution of $\xi(r)$ we have that $Z \stackrel{d}{=} \sqrt{U - rb}$ for $b \equiv Q_{11}^{-1/2}\frac{\partial y_1}{\partial \rho}$ and $U$ a standard normal random variable conditioned to $U \geq rb$. In formulas, $Z$ is supported on $\mathbb{R}_+$ and admits there the probability density

$$p_Z(z) = \frac{2ze^{-(1/2)(rb+z^2)^2}}{\sqrt{2\pi}[1 - \Phi(rb)]}.$$

Naively one expects the core size to have $\Theta(n^{1/2})$ fluctuations. This is indeed the asymptotic behavior for a fixed $\rho < \rho_c$, but as usual in phase transitions, fluctuations are enhanced near the critical point.

The distribution of the fractions of vertices with a given degree within the core can be computed along the same lines.

REMARK 2.7. As already pointed out, our proof concerns a properly defined *configuration model* whereby each edge might include the same vertex more than once ("self-loop"), and two hyperedges might include the same vertices ("double edge"). We expect a result similar to Theorem 2.3 to hold for a uniformly random hypergraph, with forbidden self-loops and double edges.

The main difficulty in proving such a generalization would be the absence of an explicit representation for the kernel of the leaf removal process. In the present case, such an expression is known and provided by Lemma 3.1. Within the uniform model, one should resort to graph enumeration formulas as the ones in [27]. This would give rise to a new Markov chain that can nevertheless be coupled to the one defined in (2.1). The thesis would follow by bounding the expected maximum distance between the trajectories of the two chains.

PROOF OF THEOREM 2.3. Putting together Propositions 2.1 and 2.2, we get that

$$P_l(n, \rho_n) = \mathbb{P}\left\{n^{1/6}\xi + \inf_t X(t) \leq 0\right\} + O(n^{-\eta}).$$

By Brownian scaling, $X(t) = \widetilde{F}^{-1/3}\widetilde{G}^{2/3}\widetilde{X}(\widetilde{F}^{2/3}\widetilde{G}^{-1/3}t)$, where $\widetilde{X}(t) = \frac{1}{2}t^2 + \widetilde{W}(t)$ and $\widetilde{W}(t)$ is also a two-sided standard Brownian motion. With $Z = \inf_t \widetilde{X}(t)$, and $Y$ a standard normal random variable which is independent of $\widetilde{X}(t)$, we clearly have that

$$P_l(n, \rho_n) = \mathbb{P}\left\{n^{1/6}\left(\frac{\partial y_1}{\partial \rho}\right)r + n^{1/6}\sqrt{Q_{11}}Y + \widetilde{F}^{-1/3}\widetilde{G}^{2/3}Z \leq 0\right\}$$

(2.18) $\quad + O(n^{-\eta})$

$$= \mathbb{E}\left\{\Phi\left(-\frac{r}{\alpha_l} - \beta_l n^{-1/6}Z\right)\right\} + O(n^{-\eta}).$$

The proof of the theorem is thus completed by a first-order Taylor expansion of $\Phi(\cdot)$ around $-r/\alpha_l$, as soon as we show that $\mathbb{E}Z = -\Omega$, and $\mathbb{E}|Z|^2$ is finite. To this end, from [20], Theorem 3.1, we easily deduce that $Z$ has the continuous distribution function $F_Z(z) = 1 - \mathcal{K}(-z)^2$ for $z < 0$, while $F_Z(z) = 1$ for $z \geq 0$, resulting after integration by parts with the explicit formula (2.16) for $\Omega$. We note in passing that taking $c = 1/2$ and $s = 0$ in [20], (5.2), provides the explicit formula (2.17) for $\mathcal{K}(x)$, en-route to which [20] also proves the finiteness of the relevant integral. Further, [20], Corollary 3.4, shows that the probability that the minimum of $\widetilde{X}(t)$ is achieved as some $t \notin [-T, T]$ is at most $A_0^{-1}e^{-A_0 T^3}$ for a positive constant $A_0$. With $\widetilde{X}(t) \geq$

$\widetilde{W}(t)$ we therefore have that

$$F_Z(z) \equiv \mathbb{P}\{Z \leq z\} \leq \mathbb{P}\left\{\inf_{t \in [-T,T]} \widetilde{X}(t) \leq z\right\} + A_0^{-1} e^{-A_0 T^3} \leq e^{-z^2/2T} + A_0^{-1} e^{-A_0 T^3}.$$

Taking $T = \sqrt{z}$ we deduce that if $z < 0$, then $F_Z(z) < C^{-1} \exp(-C|z|^{3/2})$ for some $C > 0$, which yields the stated finiteness of each moment of $Z$ (and in particular, of $\mathbb{E}|Z|^2$ and $\Omega$). □

### 3. Ensembles and transition probabilities: exact expressions.

3.1. *Model for the (initial) graph.* Throughout the paper we follow the coding literature and identify the hypergraph with a bipartite graph with two types of nodes: v-nodes, corresponding to hyperedges, and c-nodes to vertices. A graph $G$ in the ensemble $\mathcal{G} = \mathcal{G}_l(n,m)$ consists of a set of v-nodes $V \equiv [n]$, a set of c-nodes $C \equiv [m]$ and an ordered list of edges, that is, couples $(i,a)$ with $i \in V$ and $a \in C$

$$E = [(1, a_1), (1, a_2), \ldots, (1, a_l); (2, a_{l+1}), \ldots; (n, a_{(n-1)l+1}), \ldots, (n, a_{nl})],$$

where a couple $(i,a)$ appears *before* $(j,b)$ whenever $i < j$ and each v-node $i$ appears *exactly* $l$ times in the list, with $l \geq 3$ a fixed integer parameter. The total number of graphs in this ensemble is thus

(3.1) $$|\mathcal{G}_l(n,m)| = m^{nl} = \mathsf{coeff}[(e^{\mathbf{x}})^m, \mathbf{x}^{nl}](nl)!.$$

The ensemble of graphs $\mathcal{G}$ is endowed with the uniform distribution. One way to sample from this distribution is by considering the v-nodes in order, $i = 1, \ldots, n$, where for each v-node and for $j = 1, \ldots, l$, we choose independently and uniformly at random a c-node $a = a_{(i-1)l+j} \in C$ and add the couple $(i,a)$ to the list $E$. An alternative way to sample from the same distribution is by first attributing $l$ sockets to each v-node, with sockets $(i-1)l+1, \ldots, il$ attributed to the $i$th v-node. Then, we attribute $k_a$ sockets to each c-node $a$, where $k_a$'s are mutually independent Poisson($\zeta$) random variables, conditioned upon their sum being $nl$ (these sockets are ordered using any pre-established convention). Finally, we connect the v-node sockets to the c-node sockets according to a permutation $\sigma$ of $\{1, \ldots, nl\}$ that is chosen uniformly at random and independently of the choice of $k_a$'s.

Throughout the *degree* of a v-node $i$ (or c-node $a$) will refer to the number of edges $(i,b)$ [resp. $(j,a)$] it belongs to. In the hypergraph description, this corresponds to counting hyperedges, and vertices *with their multiplicity*.

3.2. *Model for the graph produced by the algorithm.* The ensemble is characterized by the nonnegative integers $(z_1, z_2) \equiv \vec{z}$, $\tau$ and $l \geq 3$, $n, m$ and denoted [4] as $\mathcal{G}(\vec{z}, \tau)$. In order for $\mathcal{G}(\vec{z}, \tau)$ to be nonempty, we require either $z_2 \geq 1$ and $z_1 + 2z_2 \leq (n - \tau)l$ or $z_2 = 0$ and $z_1 = (n - \tau)l$. An element in the ensemble is a graph $G = (U, V; R, S, T; E)$ where $U, V$ are disjoint subsets of $[n]$ with $U \cup V = [n]$ and $R, S, T$ are disjoint subsets of $[m]$ with $R \cup S \cup T = [m]$, having the cardinalities $|U| = \tau$, $|V| = n - \tau$, $|R| = m - z_1 - z_2$, $|S| = z_1$, $|T| = z_2$. Finally, $E$ is an ordered list of $(n - \tau)l$ edges

$$E = [(i_1, a_1), \ldots, (i_1, a_l); (i_2, a_{l+1}), \ldots;$$
$$(i_{n-\tau}, a_{(n-\tau-1)l+1}), \ldots, (i_{n-\tau}, a_{(n-\tau)l})],$$

such that a couple $(i, a)$ appears before $(j, b)$ whenever $i < j$. Moreover, each $i \in V$ appears as the first coordinate of exactly $l$ edges in $E$, while each $j \in U$ does not appear in any of the couples in $E$. Similarly, each $a \in R$ does not appear in $E$, each $b \in S$ appears as the second coordinate of exactly one edge in $E$, and each $c \in T$ appears in at least two such edges. The total number of elements in $\mathcal{G}(\vec{z}, \tau)$ is thus

$$h(\vec{z}, \tau) \equiv |\mathcal{G}(\vec{z}, \tau)| = \binom{m}{z_1, z_2, \cdot} \binom{n}{\tau} \mathsf{coeff}[(e^{\mathbf{x}} - 1 - \mathbf{x})^{z_2}, \mathbf{x}^{(n-\tau)l - z_1}]((n - \tau)l)!.$$

The ensemble $\mathcal{G}(\vec{z}, \tau)$ is endowed with the uniform distribution. In order to sample from it, first partition $[n]$ into $U$ and $V$ uniformly at random under the constraints $|U| = \tau$ and $|V| = (n - \tau)$ [there are $\binom{n}{\tau}$ ways of doing this], and independently partition $[m]$ to $R \cup S \cup T$ uniformly at random under the constraints $|R| = m - z_1 - z_2$, $|S| = z_1$ and $|T| = z_2$ [of which there are $\binom{m}{z_1, z_2, \cdot}$ possibilities]. Then, attribute $l$ v-sockets to each $i \in V$ and number them from 1 to $(n - \tau)l$ according to some pre-established convention. Attribute one c-socket to each $a \in S$ and $k_a$ c-sockets to each $a \in T$, where $k_a$ are mutually independent Poisson($\zeta$) random variables conditioned upon $k_a \geq 2$, and further conditioned upon $\sum_{a \in T} k_a$ being $(n - \tau)l - z_1$. Finally, connect the v-sockets and c-sockets according to a uniformly random permutation on $(n - \tau)l$ objects, chosen independently of the $k_a$'s.

3.3. *Transition probabilities.* We consider the graph process $\{G(\tau), \tau \geq 0\}$, defined as follows. The initial graph $G(0)$ is a uniformly random element of $\mathcal{G}_l(n, m)$. At each time $\tau = 0, 1, \ldots$, if there is a nonempty set of c-nodes of degree 1, one of them (let us say $a$) is chosen uniformly at random. The corresponding edge $(i, a)$ is deleted, together with all the edges incident to the v-node $i$. The graph thus obtained is $G(\tau+1)$. In the opposite case, where there are no c-nodes of degree 1 in $G(\tau)$, we set $G(\tau+1) = G(\tau)$.

---

[4]Since $n$, $m$ and $l$ do not vary during the execution of the algorithm, we leave them implicit in the ensemble notation.

We define furthermore the process $\{\vec{z}(\tau) = (z_1(\tau), z_2(\tau)), \tau \geq 0\}$ on $\mathbb{Z}_+^2$. Here $z_1(\tau)$ and $z_2(\tau)$ are, respectively, the number of c-nodes in $G(\tau)$, having degree 1 or larger than 1, which necessarily satisfy that $(n-\hat{\tau})l \geq z_1(\tau) + 2z_2(\tau)$ for $\hat{\tau} \equiv \min(\tau, \inf\{\tau' \geq 0 : z_1(\tau') = 0\})$.

LEMMA 3.1. *The process $\{\vec{z}(\tau)\,\tau \geq 0\}$ is an inhomogeneous Markov process, whose transition probabilities, denoted by*

$$W_\tau^+(\Delta \vec{z}|\vec{z}) \equiv \mathbb{P}\{\vec{z}(\tau+1) = \vec{z} + \Delta \vec{z}|\vec{z}(\tau) = \vec{z}\}$$

*[here $\Delta \vec{z} \equiv (\Delta z_1, \Delta z_2)$], are such that $W_\tau^+(\Delta \vec{z}|\vec{z}) = \mathbb{I}(\Delta \vec{z} = 0)$ in case $z_1 = 0$, whereas for $z_1 > 0$,*

$$(3.2) \quad W_\tau^+(\Delta \vec{z}|\vec{z}) = \frac{h(\vec{z}', \tau+1)}{h(\vec{z}, \tau)}(\tau+1)l! \sum_{\mathcal{D}} \binom{m - z_1' - z_2'}{q_0, p_0, \cdot} \binom{z_1'}{q_1} \binom{z_2'}{q_2} \frac{q_0}{z_1}$$
$$\times \mathsf{coeff}[(e^{\mathbf{x}} - 1 - \mathbf{x})^{p_0}(e^{\mathbf{x}} - 1)^{q_1+q_2}, \mathbf{x}^{l-q_0}].$$

*Here $z_1' = z_1 + \Delta z_1$, $z_2' = z_2 + \Delta z_2$. Also, using the notation $z_0 = m - z_1 - z_2$ and $z_0' = m - z_1' - z_2'$, the collection $\mathcal{D}$ consists of all integers $p_0, q_0, q_1, q_2 \geq 0$, satisfying the equalities*

$$(3.3) \quad \begin{cases} z_0 = z_0' - q_0 - p_0, \\ z_1 = z_1' + q_0 - q_1, \\ z_2 = z_2' + p_0 + q_1, \end{cases}$$

*and the inequalities $(n-\tau)l - (z_1 + 2z_2) \geq l - (2p_0 + q_0 + q_1) \geq q_2$, $q_0 + p_0 \leq z_0'$, $q_1 \leq z_1'$ (equivalently, $q_0 \leq z_1$), $q_2 \leq z_2'$ (equivalently, $p_0 + q_1 + q_2 \leq z_2$).*

*Moreover, conditional on $\{\vec{z}(\tau'), 0 \leq \tau' \leq \tau\}$, the graph $G(\tau)$ is uniformly distributed over $\mathcal{G}(\vec{z}, \hat{\tau})$, that is,*

$$(3.4) \quad \mathbb{P}\{G(\tau) = G|\{\vec{z}(\tau'), 0 \leq \tau' \leq \tau\}\} = \frac{1}{h(\vec{z}, \hat{\tau})}\mathbb{I}(G \in \mathcal{G}(\vec{z}, \hat{\tau})).$$

PROOF. Fixing $\tau$, $\vec{z} = \vec{z}(\tau)$ such that $z_1 > 0$, $\vec{z}' = \vec{z}(\tau+1)$ and $G' \in \mathcal{G}(\vec{z}', \tau+1)$, let $N(G'|\vec{z}, \tau)$ count the pairs of graphs $G \in \mathcal{G}(\vec{z}, \tau)$ and choices of the deleted c-node from $S$ that result with $G'$ upon applying a single step of our algorithm. Obviously, $G$ and $G'$ must be such that $R \subset R'$, $S \subseteq R' \cup S'$ and $T' \subseteq T$. So, let $q_0 \geq 0$ denote the size of $R' \cap S$, $p_0 \geq 0$ the size of $R' \cap T$, and $q_1 \geq 0$ the size of $S' \cap T$. We have $q_0 + p_0 \leq m - z_1' - z_2'$, $q_1 \leq z_1'$, and the equalities of (3.3) follow as well. Let $T^*$ denote the set of c-nodes $a \in T'$ for which $k_a > k_a'$, and denote the size of $T^*$ by $q_2 \leq z_2'$. Observe that of the $l$ edges of the v-node $i$ deleted by the algorithm in the move from $G$ to $G'$, exactly one edge hits each of the nodes in $R' \cap S$, at least one edge hits each of the nodes in $S' \cap T$, and each of the nodes in $T^*$, while at least two edges hit each of the nodes in $R' \cap T$. Consequently, $2p_0 + q_0 + q_1 + q_2 \leq l$. Since $z_1 > 0$, we know that $\hat{\tau} = \tau$ and further, $(n - \tau - 1)l \geq z'_1 + 2z'_2$, which in view of (3.3) is equivalent to $(n - \tau)l - (z_1 + 2z_2) \geq l - (2p_0 + q_0 + q_1) \geq q_2$ as claimed.

To count $N(G'|\vec{z}, \tau)$ we first select the v-node $i$ to add to $G'$ from among the $\tau + 1$ elements of $U'$, and the order (permutation) of the $l$ sockets of $i$ that we use when connecting it to the c-nodes for creating $G \in \mathcal{G}(\vec{z}, \tau)$. Summing over the set $\mathcal{D}$ of allowed values of $p_0, q_0, q_1, q_2$, for each such value we have $\binom{m - z'_1 - z'_2}{q_0, p_0, \cdot}$ ways to select the nodes of $R'$ that are assigned to $S$, $T$ and $R$, then $\binom{z'_1}{q_1}$ ways to select those of $S'$ that are assigned to $T$ and $\binom{z'_2}{q_2}$ ways to select those of $T'$ that are assigned to $T^*$. We further have $\mathsf{coeff}[(e^{\mathbf{x}} - 1 - \mathbf{x})^{p_0}(e^{\mathbf{x}} - 1)^{q_1 + q_2}, \mathbf{x}^{l - q_0}]$ ways to select the precise number of edges ($\geq 2$) from $i$ that we are to connect to each of the $p_0$ nodes in $R' \cap T$, and the precise number of edges ($\geq 1$) from $i$ that we are to connect to each of the $q_1$ nodes in $S' \cap T$ and to each of the $q_2$ nodes in $T^*$, while allocating in this manner exactly $l - q_0$ edges out of $i$ (the remaining $q_0$ then used to connect to nodes in $R' \cap S$). Finally, noting that for each of the graphs $G$ thus created we have exactly $q_0$ ways to choose the deleted node from $S$ while still resulting with the graph $G'$, we conclude that

$$N(G'|\vec{z}, \tau) = (\tau + 1)l! \sum_{\mathcal{D}} \binom{m - z'_1 - z'_2}{q_0, p_0, \cdot} \binom{z'_1}{q_1} \binom{z'_2}{q_2}$$
$$\times q_0 \, \mathsf{coeff}[(e^{\mathbf{x}} - 1 - \mathbf{x})^{p_0}(e^{\mathbf{x}} - 1)^{q_1 + q_2}, \mathbf{x}^{l - q_0}].$$

We start at $\tau = 0$ with a uniform distribution of $G(0)$ within each possible ensemble $\mathcal{G}(\vec{z}(0), 0)$. Since $N(G'|\vec{\omega}, \tau)$ depends on $G'$ only via $\vec{\omega}'$, it follows by induction on $\tau = 1, 2, \ldots$ that this property, namely (3.4), is preserved as long as $\hat{\tau} = \tau$, since if $z_1(\tau) > 0$, then

$$\mathbb{P}\{G(\tau + 1) = G'|\{\vec{z}(\tau'), 0 \leq \tau' \leq \tau\}\} = \frac{1}{z_1} \frac{N(G'|\vec{z}(\tau), \tau)}{h(\vec{z}(\tau), \tau)}$$

is the same for all $G' \in \mathcal{G}(\vec{z} + \Delta\vec{z}, \tau + 1)$. Since there are exactly $h(\vec{z} + \Delta\vec{z}, \tau + 1)$ graphs in this ensemble, we thus recover also (3.2). Finally, noting that $G(\tau) = G(\hat{\tau})$ and $\vec{z}(\tau) = \vec{z}(\hat{\tau})$ we deduce that (3.4) holds also when $\hat{\tau} < \tau$. □

## 4. Asymptotic expressions.

4.1. *Properties of the ordinary differential equations.* We derive here the properties of solutions of the ODEs (2.5) and (2.8) that are needed for our analysis. This is based on the continuity of $(\vec{x},\theta) \mapsto \mathfrak{p}_a(\vec{x},\theta)$, $a = 0,1,2$ on the following compact subsets of $\mathbb{R}^2 \times \mathbb{R}_+$:

$$\widehat{q}(\varepsilon) \equiv \{(\vec{x},\theta): -l \leq x_1; 0 \leq x_2; \theta \in [0, 1-\varepsilon]; 0 \leq (1-\theta)l - \max(x_1,0) - 2x_2\},$$

and $\widehat{q}_+(\varepsilon) = \widehat{q}(\varepsilon) \cap \{x_1 \geq 0\}$, as stated in

LEMMA 4.1. *For any $\varepsilon > 0$, the functions $(\vec{x},\theta) \mapsto \mathfrak{p}_a(\vec{x},\theta)$, $a = 0, 1, 2$ are $[0,1]$-valued, Lipschitz continuous on $\widehat{q}(\varepsilon)$. Further, on $\widehat{q}_+(\varepsilon)$ the functions $(\vec{x},\theta) \mapsto \mathfrak{p}_a(\vec{x},\theta)$ have Lipschitz continuous partial derivatives.*

PROOF. Fixing $\varepsilon > 0$, the stated Lipschitz continuity holds for $\mathfrak{p}_0(\vec{x},\theta)$ since both $\max(x_1,0)$ and $1/(1-\theta)$ are Lipschitz continuous and bounded on $\widehat{q}(\varepsilon)$. Further, $\mathfrak{p}_0(\vec{x},\theta) \in [0,1]$ throughout $\widehat{q}(\varepsilon)$. Setting $f_1(0) = 2$, note that $f_1: \mathbb{R}_+ \to [2,\infty)$ of (2.3) is a monotone increasing, twice continuously differentiable function, with $f_1'(\lambda) = [(e^\lambda - 1)^2 - \lambda^2 e^\lambda]/(e^\lambda - 1 - \lambda)^2$ strictly positive and bounded away from zero throughout $\mathbb{R}_+$. Consequently, the inverse mapping $f_1^{-1}$ is well defined and twice continuously differentiable on $[2,\infty)$, from which we deduce that for each $\delta > 0$ the nonnegative function $\lambda(\vec{x},\theta)$ is well defined, bounded and continuously differentiable on the compact set $\widehat{q}(\varepsilon) \cap \{(\vec{x},\theta): x_2 \geq \delta\}$. Though $\lambda(\vec{x},\theta) \uparrow \infty$ as $x_2 \downarrow 0$, note that $\mathfrak{p}_2 = (1-\mathfrak{p}_0)(1-g(\lambda))$ for $g(\lambda) \equiv \lambda/(e^\lambda - 1)$. In particular, since $\mathfrak{p}_2 = 1 - \mathfrak{p}_0$ in case $x_2 = 0$, it follows that $\mathfrak{p}_2(\vec{x},\theta)$ is continuous throughout $\widehat{q}(\varepsilon)$. Since $\mathfrak{p}_0(\vec{x},\theta)$ is Lipschitz continuous on $\widehat{q}(\varepsilon)$, the Lipschitz continuity of $\mathfrak{p}_2$ follows by showing that, for $x_1 \neq 0$, $g(\lambda(\vec{x},\theta))$ has bounded derivatives as $x_2 \downarrow 0$. By letting $\vec{\xi} \equiv (\vec{x},\theta) \in \widehat{q}(\varepsilon)$, we have $\partial_{\xi_i} g(\lambda) = g'(\lambda)\partial_{\xi_i}\lambda$. Using the definition (2.3), and recalling that $f_1'(\lambda)$ is bounded away from zero, it follows that $|\partial_{\xi_i}\lambda| \leq C x_2^{-2}$ as $x_2 \downarrow 0$. On the other hand, $|g'(\lambda)| \leq Ce^{-\lambda} \leq Ce^{-C'/x_2}$ in the same limit thus implying that $\partial_{\xi_i} g(\lambda)$ is bounded.

Further, the identity (2.3) is equivalent to $\mathfrak{p}_0 + \mathfrak{p}_1 + \mathfrak{p}_2 = 1$, which thus implies that $\mathfrak{p}_1$ is also Lipschitz continuous on $\widehat{q}(\varepsilon)$. Finally, since both $\lambda(\vec{x},\theta)$ and $x_2$ are nonnegative throughout $\widehat{q}(\varepsilon)$, the same applies for $\mathfrak{p}_1$ and $\mathfrak{p}_2$, and consequently, $\mathfrak{p}_a \in [0,1]$ for $a = 0, 1, 2$.

Considering for the remainder of the proof $\vec{\xi} = (\vec{x},\theta) \in \widehat{q}_+(\varepsilon)$, we replace $\max(x_1,0)$ by $x_1$ in the definition of $(\mathfrak{p}_0, \mathfrak{p}_1, \mathfrak{p}_2)$. The stated regularity of $\mathfrak{p}_0$ is then obvious and as before the regularity of $\mathfrak{p}_1 = 1 - \mathfrak{p}_0 - \mathfrak{p}_2$ follows from that of $\mathfrak{p}_2 = (1-\mathfrak{p}_0)(1-g(\lambda))$. To this end, we see that it suffices to show that $\partial_{\xi_i} g(\lambda) = g'(\lambda)\partial_{\xi_i}\lambda$, are Lipschitz continuous in $\vec{\xi}$ on the compact set $\widehat{q}_+(\varepsilon)$. As seen already $\lambda \mapsto g'(\lambda)$ is bounded and Lipschitz continuous on $\mathbb{R}_+$, and $\partial_{\xi_i}\lambda(\vec{\xi})$ is bounded and has bounded derivatives on $\widehat{q}_+(\varepsilon) \cap \{\vec{\xi}: x_2 \geq \delta\}$. The proof is completed by showing that $\partial_{\xi_j}[g'(\lambda)\partial_{\xi_i}\lambda] = g''(\lambda)\partial_{\xi_i}\lambda\partial_{\xi_j}\lambda + g'(\lambda)\partial_{\xi_i}\partial_{\xi_j}\lambda$ converges to zero as $x_2 \to 0$. This is proved similarly to what was already done for the first-order derivatives. Indeed, the first and second derivatives of $\lambda \mapsto f_1(\lambda)$ as well as $\vec{\xi} \mapsto \partial_{\xi_i}[x_2 f_1(\lambda)]$ and its partial derivatives are all bounded, hence $|\partial_{\xi_i}\lambda| \leq Cx_2^{-2}$ and $|\partial_{\xi_i}\partial_{\xi_j}\lambda| \leq Cx_2^{-4}$ as $x_2 \to 0$, which since $\lambda \to \infty$ inversely proportional to $x_2 \to 0$, is more than compensated by the exponential decay in $\lambda$ of $g'$ and $g''$. □

Setting $h_\rho(u) \equiv u - 1 + \exp(-lu^{l-1}/\rho)$ and the finite and positive *critical density*

$$\rho_c \equiv \inf\{\rho > 0 : h_\rho(u) > 0 \ \forall u \in (0,1]\},$$

we have the following properties of the ODEs.

PROPOSITION 4.2. *For any $\rho > 0$, the ODE (2.5) admits a unique solution $\vec{y}$ subject to the initial conditions (2.4), and the ODE (2.8) admits a unique, positive definite, solution $\mathbb{Q}$ subject to the initial conditions (2.10), such that:*

(a) *For any $\varepsilon > 0$, $\theta < 1 - \varepsilon$, we have that $(\vec{y}(\theta,\rho),\theta)$ is in the interior of $\hat{q}(\varepsilon)$, with both functions $(\theta,\rho) \mapsto \vec{y}$ and $\theta \mapsto \mathbb{Q}$ Lipschitz continuous on $(\theta,\rho) \in [0,1-\varepsilon] \times [\varepsilon,1/\varepsilon]$.*

(b) *Let $u(\theta) \equiv (1-\theta)^{1/l}$ and $\theta_-(\rho) \equiv \inf\{\theta \geq 0 : h_\rho(u(\theta)) < 0\} \wedge 1$. Then, for $\theta \in [0,\theta_-(\rho)]$*

$$(4.1) \quad y_1(\theta,\rho) = lu(\theta)^{l-1}[u(\theta) - 1 + e^{-\gamma u(\theta)^{l-1}}],$$

$$(4.2) \quad y_2(\theta,\rho) = \frac{l}{\gamma}[1 - e^{-\gamma u(\theta)^{l-1}} - \gamma u(\theta)^{l-1} e^{-\gamma u(\theta)^{l-1}}]$$

*(where $\gamma = l/\rho$). In particular, $(\theta,\rho) \mapsto \vec{y}$ is infinitely continuously differentiable and $(\theta,\rho) \mapsto \mathbb{Q}$ is Lipschitz continuous on $\{(\theta,\rho) : \theta \leq \min(\theta_-(\rho), 1 - \varepsilon), \varepsilon \leq \rho \leq 1/\varepsilon\}$.*

(c) *Let $\theta_*(\rho) \equiv \inf\{\theta \geq 0 : h_\rho(u(\theta)) \leq 0\}$. Then, $\theta_*(\rho) = \sup\{\theta \leq 1 : y_1(\theta', \rho) > 0 \text{ for all } \theta' \in [0,\theta)\}$ and the critical density is such that*

$$(4.3) \quad \rho_c = \inf\{\rho > 0 : \theta_*(\rho) = 1\} = \inf\{\rho > 0 : y_1(\theta,\rho) > 0 \ \forall \theta \in [0,1)\}.$$

(d) *The critical time $\theta_c \equiv \theta_*(\rho_c)$ is in $(0,1)$, whereas $\theta_-(\rho_c) = 1$. For $\rho = \rho_c$ the infinitely continuously differentiable function $\theta \mapsto y_1(\theta)$ is positive for $\theta \neq \theta_c$, $\theta \neq 1$, with $y_1(\theta_c) = y_1'(\theta_c) = 0$, and $y_1''(\theta_c) > 0$, while $y_1(1-\delta) = l\delta + o(\delta)$ for any $\delta > 0$.*

PROOF. (a) For any $\rho > 0$ the initial condition $\vec{y}(0)$ of (2.4) is such that $(\vec{y}(0),0)$ is in the interior of $\hat{q}(\varepsilon)$. Further, fixing $\varepsilon > 0$, by Lemma 4.1 we have that $\vec{F}(\vec{x},\theta)$ is bounded and Lipschitz continuous on $\hat{q}(\varepsilon)$. Consequently, for $\theta \in [0,\theta_\varepsilon]$ there exists a unique solution $\vec{y}(\theta)$ of the ODE (2.5) [i.e., $\frac{d\vec{y}}{d\theta} = \vec{F}(\vec{y},\theta)$], starting at this initial condition, where $\theta_\varepsilon = \inf\{\theta > 0 : (\vec{y}(\theta),\theta) \notin \hat{q}(\varepsilon)\}$ is strictly positive, and $(\vec{y}(\theta_\varepsilon),\theta_\varepsilon)$ is necessarily on the boundary of $\hat{q}(\varepsilon)$. We proceed to verify that $\theta_\varepsilon = 1 - \varepsilon$ by showing that:

(i) $y_1(\theta_\varepsilon) > -l$. Indeed, since $y_1(0) > 0$ and $F_1(\vec{y}(\theta), \theta) \geq -l$, we have $y_1(\theta_\varepsilon) \geq -l\theta_\varepsilon > -l$.

(ii) $y_2(\theta_\varepsilon) > 0$. In fact $x_2 = 0$ implies $\mathfrak{p}_1(\vec{x}, \theta) = 0$, and therefore $F_2(\vec{x}, \theta) = 0$. By the Lipschitz continuity of $F_2$ on $\widehat{q}(\varepsilon)$ it follows that $F_2(\vec{x}, \theta) \geq -Cx_2$ for some finite $C$ and all $x_2$ in $\widehat{q}(\varepsilon)$. Therefore, $y_2(\theta_\varepsilon) \geq y_2(0)e^{-C\theta_\varepsilon} > 0$.

(iii) $v(\theta_\varepsilon) > 0$, where $v(\theta) = w(\vec{y}(\theta), \theta)$ for $w(\vec{x}, \theta) = l(1-\theta) - \max(x_1, 0) - 2x_2$. Indeed, note that $v(0) > 0$ and

$$\frac{dv}{d\theta} = [-l + 2(l-1)\mathfrak{p}_1(\vec{y}, \theta)]\mathbb{I}(y_1(\theta) \leq 0) - (l-1)\mathfrak{p}_2(\vec{y}, \theta)\mathbb{I}(y_1(\theta) > 0).$$

Further, recall that if $w(\vec{x}, \theta) = 0$, then $\mathfrak{p}_2(\vec{x}, \theta) = 0$, and if in addition $x_1 \leq 0$, then also $\mathfrak{p}_1(\vec{x}, \theta) = 1$. Hence, by the Lipschitz continuity of $\mathfrak{p}_1(\cdot, \cdot)$ and $\mathfrak{p}_2(\cdot, \cdot)$ on $\widehat{q}(\varepsilon)$ we have that $\mathfrak{p}_2(\vec{x}, \theta) \leq 2Cw(\vec{x}, \theta)$ and $\mathfrak{p}_1(\vec{x}, \theta) \geq (1 - Cw(\vec{x}, \theta))\mathbb{I}(x_1 \leq 0)$ for some finite $C > 0$, throughout $\widehat{q}(\varepsilon)$. Since $l \geq 2$, it follows that $\frac{dv}{d\theta} \geq -2(l-1)Cv(\theta)$ for all $\theta \in [0, \theta_\varepsilon]$, resulting with $v(\theta_\varepsilon) \geq v(0)e^{-2(l-1)C\theta_\varepsilon} > 0$.

Lemma 4.1 further implies that for any $a, b \in \{1, 2\}$ both $A_{ab}(\vec{x}, \theta) \equiv \partial_{x_b} F_a(\vec{x}, \theta)$ and $G_{ab}(\vec{x}, \theta)$ are uniformly bounded over $\widehat{q}(\varepsilon)$. The linear ODE (2.8) has these functions as its coefficients, for $\vec{x} = \vec{y}(\theta)$. We thus deduce that there exists a unique solution $\mathbb{Q}(\theta)$ of the initial value problem for this ODE at least for $\theta \in [0, \theta_\varepsilon]$. With $\theta_\varepsilon = 1 - \varepsilon$ and $\varepsilon > 0$ arbitrarily small we established the existence of a unique solution $(\vec{y}, \mathbb{Q})$ for $\theta \in [0, 1)$.

It also follows from the above discussion that $\theta_\varepsilon = 1 - \varepsilon$ and $\vec{y}(\theta, \rho)$ is Lipschitz continuous in $\theta$ on $[0, 1 - \varepsilon] \times [\varepsilon, 1/\varepsilon]$. Further, applying Gronwall's lemma, the Lipschitz continuity of $\vec{F}(\vec{x}, \theta)$ implies that the solution $\vec{y}(\theta)$ of the ODE is then also Lipschitz continuous with respect to the initial condition $\vec{y}(0)$, with a uniform in $\theta \leq 1 - \varepsilon$ bound on the corresponding Lipschitz norm. Clearly, $\vec{y}(0)$ of (2.4) is differentiable in $\rho$ with a uniformly bounded derivative when $\rho \in [\varepsilon, 1/\varepsilon]$. Consequently, we arrive at the stated Lipschitz continuity of $(\theta, \rho) \mapsto \vec{y}(\theta, \rho)$.

The same argument shows that the initial conditions (2.10) for the ODE (2.8) are bounded in $\rho \in [\varepsilon, 1/\varepsilon]$. Further, $\vec{y}(\theta, \rho)$ stays in $\widehat{q}(\varepsilon)$ and with the coefficients of the linear ODE (2.8) uniformly bounded on $[0, 1-\varepsilon] \times [\varepsilon, 1/\varepsilon]$, its solution $\mathbb{Q}$ is also Lipschitz continuous in $\theta$. Suppressing the dependence of the various matrices on $\rho$, set $\mathbb{B}_\zeta^\zeta = \mathbb{I}$ and, for $\theta \geq \zeta$

$$(4.4) \qquad \frac{d\mathbb{B}_\zeta^\theta}{d\theta} = \mathbb{A}(\vec{y}(\theta), \theta)\mathbb{B}_\zeta^\theta.$$

It is easy to check that the unique solution of (2.8) is given by

$$(4.5) \qquad \mathbb{Q}(\theta) = \mathbb{B}_0^\theta \mathbb{Q}(0)(\mathbb{B}_0^\theta)^\dagger + \int_0^\theta \mathbb{B}_\zeta^\theta \mathbb{G}(\vec{y}(\zeta), \zeta)(\mathbb{B}_\zeta^\theta)^\dagger \, d\zeta,$$

for the nonnegative definite matrix $\mathbb{G}(\vec{x}, \theta)$ of (2.9). In particular, starting from the symmetric, positive definite $\mathbb{Q}(0)$ of (2.10), this implies that $\mathbb{Q}(\theta)$ is nonnegative definite. Further, since $\det \mathbb{B}_0^0 = 1$ and

$$\frac{d(\det \mathbb{B}_0^\theta)}{d\theta} = (\det \mathbb{B}_0^\theta) \text{Trace}(\mathbb{A}(\vec{y}(\theta), \theta)),$$

with the entries of $\mathbb{A}(\vec{x}, \theta)$ uniformly bounded, it follows that $\det \mathbb{B}_0^\theta > 0$, hence the solution $\mathbb{Q}(\theta)$ of (2.8) is positive definite.

(b) Though this is a special case of a result of [26], we provide its short proof for the reader's convenience. We first check that $\vec{y}(\theta, \rho)$ of (4.1) and (4.2) is the unique solution of the ODE (2.5) for $\theta \in [0, \theta_-(\rho)]$. To this end, first note that for $\theta = 0$ the functions $\vec{y}(\theta, \rho)$ of (4.1) and (4.2) satisfy the initial condition (2.4). Further, the function $y_1(\theta, \rho)$ of (4.1) is nonnegative for $\theta \in [0, \theta_-(\rho)]$. Hence, upon substituting $y_1(\theta, \rho)$ for $\max(x_1, 0)$ and $y_2(\theta, \rho)$ for $x_2$ on the right-hand side of (2.3), and noticing that $(1 - \theta) = u(\theta)^l$, it is not hard to verify that this equation is satisfied by $\lambda(\vec{y}(\theta), \theta) = \gamma u(\theta)^{l-1}$. Using this value of $\lambda$ yields after some algebra that $F_1(\vec{y}(\theta), \theta) = -1 - \frac{(l-1)}{u}(u - 1 + e^{-\gamma u^{l-1}} - \gamma u^{l-1} e^{-\gamma u^{l-1}})$ and $F_2(\vec{y}(\theta), \theta) = -\gamma(l-1)u^{l-2}e^{-\gamma u^{l-1}}$. With $\frac{du}{d\theta} = -u^{1-l}/l$, it is then immediate to verify that the functions given by (4.1) and (4.2) indeed satisfy (2.5) as long as $\theta \leq \theta_-(\rho)$. Clearly, $\vec{y}(\theta, \rho)$ of (4.1) and (4.2) is infinitely continuously differentiable on $[0, 1-\varepsilon] \times [\varepsilon, 1/\varepsilon]$. With $\mathbb{Q}(\theta, \rho)$ Lipschitz continuous in $\theta$ [by (a)], it remains only to show that this function is Lipschitz continuous with respect to $\rho \in [\varepsilon, 1/\varepsilon]$. Since the ODE (2.8) is linear and of bounded coefficients, with initial condition $\mathbb{Q}(0)$ of (2.10) that is Lipschitz continuous in $\rho \in [\varepsilon, 1/\varepsilon]$ it suffices to show that the coefficients $A_{ab}(\vec{x}, \theta)$ and $G_{ab}(\vec{x}, \theta)$, are Lipschitz continuous in $\vec{x}$ on $\widehat{q}_+(\varepsilon)$. We deduce the latter property from Lemma 4.1 upon noting that these coefficients are smooth bounded functions of $\mathfrak{p}_a$ and $\partial_{x_b}\mathfrak{p}_a$.

(c) We turn to verify that $\rho_c$ satisfies (4.3). We have already seen that the solution of (2.5) starting at $\vec{y}(0)$ of (2.4) is given for $\theta \leq \theta_*(\rho) \leq \theta_-(\rho)$ by (4.1) and (4.2), and in particular is such that $y_1(\theta, \rho) > 0$ for all $\theta < \theta_*(\rho)$. Further, $\rho \mapsto \theta_*(\rho)$ is monotone nondecreasing, and since $u(1) = 0$, we see that $\theta_*(\rho) \leq 1$ for all $\rho > 0$. Thus, to complete the proof it suffices to assume that for some positive $\delta$ and $\rho_0$ the solution of the ODE (2.5) is such that $y_1(\theta, \rho_0) > 0$ for all $\theta \in [0, \theta_*(\rho_0) + \delta]$ and arrive at a contradiction. To this end, note that for $\rho = \rho_0$ and $\theta \leq \theta_*(\rho_0) + \delta$, the solution of (2.5) must also satisfy the modified ODE

$$\frac{d\vec{y}}{d\theta}(\theta) = \vec{F}^*(\vec{y}(\theta), \theta), \tag{4.6}$$

where $\vec{F}^*(\vec{x}, \theta) = (-1 + (l-1)(\mathfrak{p}_1^* - \mathfrak{p}_0^*), -(l-1)\mathfrak{p}_1^*)$ and $\mathfrak{p}_a^*$ are obtained by replacing $\max(x_1, 0)$ in (2.2) and (2.3) with $x_1$. Modifying the set $\widehat{q}(\varepsilon)$ in the same manner, it is easy to verify that the statement and proof of Lemma 4.1 remain valid for $\mathfrak{p}_a^*(\vec{x}, \theta)$ (apart from the fact that the latter are not $[0,1]$ valued). We also find that $\theta_\varepsilon = 1 - \varepsilon$ for the ODE (4.6), from which we can deduce that the latter ODE also admits a unique solution subject to the initial condition (2.4). Further, the preceding computations show that for *every* $\rho > 0$ the solution of (4.6) starting at (2.4) is given by (4.1) and (4.2). In particular, at $\rho = \rho_0$ this is also the solution of the ODE (2.5) on $[0, \theta_*(\rho_0) + \delta]$. However, by definition of $\theta_*(\rho)$, necessarily $y_1(\theta, \rho_0)$ of (4.1) is nonpositive for some $\theta \in (\theta_*(\rho_0), \theta_*(\rho_0) + \delta)$, resulting with the desired contradiction.

(d) Simple calculus shows that either $u \mapsto h_\rho(u)$ is monotone increasing and positive on $(0, \infty)$, which happens for all $\rho$ large enough, or $h'_\rho(u) = 0$ has exactly two positive solutions, $u_1 = u_1(\rho)$ corresponding to a local maximum of $h_\rho$ and $u_2 = u_2(\rho) > u_1$ corresponding to a local minimum of $h_\rho$. With $h_\rho(0) = 0$ and $h_\rho(\cdot)$ positive on $[1, \infty)$, while $h_\rho(u_2(\rho)) < 0$ for all $\rho > 0$ small enough, it follows from the definition of $\rho_c$ that $h_{\rho_c}(u_2) = 0$ at $u_2 = u_2(\rho_c) \in (0, 1)$ and $h_{\rho_c}(u)$ is positive at any positive $u \neq u_2$. Hence, by definition $\theta_-(\rho_c) = 1$ while $\theta_*(\rho_c) = 1 - u_2(\rho_c)^l \in (0, 1)$. From part (b) of the proposition we thus have that at $\rho = \rho_c$ the function $y_1(\theta)$ is infinitely continuously differentiable, with $y_1(\theta) = \overline{h}(u(\theta))$ for $\overline{h}(u) = lu^{l-1}h_{\rho_c}(u)$ [cf. (4.1)]. In particular, $y_1(\theta)$ is then zero when $\theta = \theta_c$ or $\theta = 1$ and positive elsewhere (per the preceding analysis of $h_{\rho_c}$). Further, at $\theta = \theta_c$ we have $u(\theta) = u_2(\rho_c)$, an isolated minimizer of $\overline{h}(u)$, and as $u'(\theta_c) > 0$, it follows by elementary calculus that $y_1'(\theta_c) = 0$ and $y_1''(\theta_c) > 0$. Also, $\overline{h}(u) = lu^l(1 + O(u^{l-2}))$ for small $u$, hence $y_1(1-\delta) = l\delta + o(\delta)$ at $\rho = \rho_c$. □

We conclude this section by showing that the discrete recursions corresponding to the mean and covariance of the process $\vec{z}(\cdot)$ of (2.12) are near the solution of the relevant ODEs (at least for $\rho$ near $\rho_c$ and up to time $\tau_n \equiv \lfloor n\theta_c - n^\beta \rfloor$). More precisely, for $\widetilde{\mathbb{A}}_\tau \equiv \mathbb{I}_{\tau < \tau_n} \mathbb{A}(\vec{y}(\tau/n, \rho), \tau/n)$, let

$$(4.7) \quad \vec{y}^*(\tau + 1) = \vec{y}^*(\tau) + n^{-1}\widetilde{\mathbb{A}}_\tau(\vec{y}^*(\tau) - \vec{y}(\tau/n)) + n^{-1}\vec{F}(\vec{y}(\tau/n), \tau/n),$$

starting at $\vec{y}^*(0) \equiv \vec{y}(0, \rho)$ and consider the positive definite matrices

$$(4.8) \quad \mathbb{Q}_\tau = \widetilde{\mathbb{B}}_0^{\tau-1}\mathbb{Q}(0, \rho)(\widetilde{\mathbb{B}}_0^{\tau-1})^\dagger + \sum_{\sigma=0}^{\tau-1} \widetilde{\mathbb{B}}_{\sigma+1}^{\tau-1}\mathbb{G}(\vec{y}(\sigma/n), \sigma/n)(\widetilde{\mathbb{B}}_{\sigma+1}^{\tau-1})^\dagger$$

for $\widetilde{\mathbb{B}}_\sigma^\tau$ of (2.13). Then:

LEMMA 4.3. *Fixing $\beta \in (3/4, 1)$ and $\beta' < 2\beta - 1$, we have for all $n$ large enough and $|\rho - \rho_c| \leq n^{\beta'-1}$*

$$(4.9) \quad \left| n^{1/2}y_1^*(\tau_n) - \frac{\widetilde{F}}{2}n^{2\beta-3/2} - n^{1/2}(\rho - \rho_c)\frac{\partial y_1}{\partial \rho}(\theta_c, \rho_c) \right| \leq Cn^{3\beta-5/2},$$

*the matrices* $\{\widetilde{\mathbb{B}}_\sigma^\tau : \sigma, \tau \leq n\}$ *and their inverses are uniformly bounded with respect to the $L_2$-operator norm (denoted $\|\cdot\|$), and*

$$\|\mathbb{Q}_{\tau_n} - \mathbb{Q}(\theta_c, \rho_c)\| \leq Cn^{\beta-1}, \tag{4.10}$$

*for some finite $C = C(\beta, \beta')$ and all $n$.*

PROOF. Recall part (a) of Proposition 4.2 that $\vec{y}(\theta, \rho) \in \widehat{q}(\varepsilon)$ for $\theta \leq 1 - 2\varepsilon$ and $\rho \in [\varepsilon, 1/\varepsilon]$. Thus, fixing $\beta$, $\beta'$ and $0 < \varepsilon < (1-\theta_c)/2$, it follows that the operator norm of $\widetilde{\mathbb{A}}_\tau$ is uniformly bounded over $\tau \leq \tau_n$, $|\rho - \rho_c| \leq n^{\beta'-1}$ and $n \geq n_0$ (hereafter $n_i$ and $c_i$, $i = 0, 1, \ldots$, are two nondecreasing sequences of finite constants, each depending only on $l$, $\beta$, $\beta'$ and $\varepsilon$). Consequently, the matrices $\widetilde{\mathbb{B}}_\sigma^\tau$ of (2.13) and their inverses are also uniformly bounded with respect to the $L_2$ operator norm for $n \geq n_1$, $\sigma$, $\tau$ and $\rho$ as before.

We proceed to show that $\{\vec{y}^*(\tau), \tau \leq \tau_n\}$ is close to the solution $\vec{y}(\cdot, \rho)$ of the ODE (2.5). To this end, let $D_n^*(\tau) \equiv \vec{y}^*(\tau) - \vec{y}(\tau/n, \rho)$, noting that by definition $D_n^*(\tau+1) = (\mathbb{I} + n^{-1}\widetilde{\mathbb{A}}_\tau)D_n^*(\tau) + \vec{\xi}_n(\tau)$ for $\tau \geq 0$, with $D_n^*(0) = 0$ and

$$\vec{\xi}_n(\tau) = \int_{\tau/n}^{(\tau+1)/n} [\vec{F}(\vec{y}(\tau/n), \tau/n) - \vec{F}(\vec{y}(\theta), \theta)]\, d\theta.$$

By the Lipschitz continuity of $(\theta, \rho) \mapsto \vec{y}$ on $[0, 1-\varepsilon] \times [\varepsilon, 1/\varepsilon]$ (see Proposition 4.2), we know that $\|\vec{y}(\theta) - \vec{y}(\tau/n)\| \leq c_0/n$ for some finite $c_0$, all $\theta \in [\tau/n, \tau/n + 1/n]$ and any $\tau < (1-\varepsilon)n$. Further, since $\|\vec{F}(\vec{x}, \theta)\| \leq 2l$ and $(\vec{x}, \theta) \mapsto \vec{F} = (-1 + (l-1)(\mathfrak{p}_1 - \mathfrak{p}_0), -(l-1)\mathfrak{p}_1)$ is Lipschitz continuous on $\widehat{q}(\varepsilon)$ (see Lemma 4.1), we deduce that for some finite constant $C_* = C_*(l, \varepsilon)$, all $n$, $\tau < (1-\varepsilon)n$ and $\rho \in [\varepsilon, 1/\varepsilon]$,

$$\|\vec{\xi}_n(\tau)\| \leq \frac{1}{n} \sup_{\theta \in [\tau/n, \tau/n+1/n]} \|\vec{F}(\vec{y}(\theta), \theta) - \vec{F}(\vec{y}(\tau/n), \tau/n)\| \leq C_* n^{-2}. \tag{4.11}$$

Since $D_n^*(\tau) = \sum_{\sigma=0}^{\tau-1} \widetilde{\mathbb{B}}_{\sigma+1}^{\tau-1} \vec{\xi}_n(\sigma)$, and $\|\widetilde{\mathbb{B}}_\sigma^\tau\|$ are uniformly bounded, we deduce that

$$\sup_{n \geq n_2} \sup_{|\rho - \rho_c| \leq n^{\beta'-1}} \sup_{\tau \leq \tau_n} n\|\vec{y}^*(\tau) - \vec{y}(\tau/n, \rho)\| \leq c_1 < \infty. \tag{4.12}$$

Let $\theta_n \equiv \tau_n/n$, $\Delta\theta_n \equiv \theta_n - \theta_c = -n^{\beta-1}$ and $\Delta\rho \equiv \rho - \rho_c$. Note that $y_1(\theta, \rho_c) \geq \overline{c}(\Delta\theta_n)^2$ for some $\overline{c} > 0$, all $n$ and $\theta \in [0, \theta_n]$ [see part (d) of Proposition 4.2]. Further, recall that $|\Delta\rho| \leq n^{\beta'-1} = o((\Delta\theta_n)^2)$ by our choice of $\beta' < 2\beta - 1$. The Lipschitz continuity of $\rho \mapsto \vec{y}(\theta, \rho)$ for $\theta \leq \theta_-(\rho)$ thus implies that both $(\theta_c, \rho_c)$ and $(\theta_n, \rho)$ for $n \geq n_4$ and $|\rho - \rho_c| \leq n^{\beta'-1}$ are in the set $\mathcal{A}_\varepsilon \equiv \{(\theta, \rho) : \theta \leq \min(\theta_-(\rho), 1-\varepsilon), \varepsilon \leq \rho \leq 1/\varepsilon\}$ where $(\theta, \rho) \mapsto \vec{y}$ is infinitely continuously differentiable [see part (b) of Proposition 4.2]. Hence, by Taylor expanding $y_1(\cdot)$ around $(\theta_c, \rho_c)$ where $y_1 = \frac{\partial y_1}{\partial \theta} = 0$, we obtain that for some $c_2'$, $c_2$ and all $n$,

(4.13)
$$\left| y_1(\theta_n, \rho) - \Delta\rho \frac{\partial y_1}{\partial \rho} - \frac{1}{2}(\Delta\theta_n)^2 \frac{\partial^2 y_1}{\partial \theta^2} \right|$$
$$\leq c_2'(|\Delta\rho| + |\Delta\theta_n|)(|\Delta\rho| + (\Delta\theta_n)^2) \leq c_2 n^{3(\beta-1)}$$

[with all partial derivatives evaluated at $(\theta_c, \rho_c)$]. Recall that $\widetilde{F} \equiv \frac{\partial^2 y_1}{\partial \theta^2}$, so the left-hand side of (4.9) is bounded above by

$$n^{1/2}\|\vec{y}^*(\tau_n) - \vec{y}(\tau_n/n, \rho)\| + n^{1/2}\left| y_1(\theta_n, \rho) - \Delta\rho \frac{\partial y_1}{\partial \rho} - \frac{1}{2}(\Delta\theta_n)^2 \frac{\partial^2 y_1}{\partial \theta^2} \right|.$$

Thus, controlling the first term via (4.12) and the second term via (4.13) yields the bound of (4.9).

Turning now to the proof of (4.10), recall that the solution $\mathbb{Q}(\cdot)$ of (2.8) is Lipschitz continuous in $(\theta, \rho)$ on the set $\mathcal{A}_\varepsilon$ [see part (b) of Proposition 4.2]. As both $(\theta_c, \rho_c)$ and $(\theta_n, \rho)$, $n \geq n_4$ are in this set, it follows that for some finite $c_3'$, $c_3$ and all $n$,

$$\|\mathbb{Q}(\theta_n, \rho) - \mathbb{Q}(\theta_c, \rho_c)\| \leq c_3'(|\Delta\theta_n| + |\Delta\rho|) \leq c_3 n^{\beta-1}$$

(recall that $\beta' < 2\beta - 1 < \beta$). Further, $\mathbb{Q}(\theta, \rho)$ is given by (4.5), where the matrices $\mathbb{G}(\vec{y}(\zeta), \zeta)$ are bounded and Lipschitz continuous in $\zeta$ (with respect to the $L_2$ operator norm) uniformly in $n \geq n_4$ and $\zeta \leq \theta_n$. The same uniform boundedness applies for $\mathbb{Q}(0, \rho)$ and $\mathbb{B}^\theta_\zeta$, $0 \leq \zeta \leq \theta \leq 1 - \varepsilon$ [see proof of part (a) of Proposition 4.2]. Hence, comparing (4.5) and (4.8) we thus deduce that (4.10) is an immediate consequence of

(4.14) $$\sup_{0 \leq \zeta \leq \theta_n} \|\widetilde{\mathbb{B}}^{\tau_n-1}_{\lceil n\zeta \rceil} - \mathbb{B}^{\theta_n}_\zeta\| \leq c_4 n^{-1},$$

holding for some finite $c_4$ and all $n$. To this end, let $D_n(\sigma, \tau) \equiv \|\widetilde{\mathbb{B}}^{\tau-1}_\sigma - \mathbb{B}^{\tau/n}_{\sigma/n}\|$, noting that by the definition of $\mathbb{B}^\theta_\zeta$ and $\widetilde{\mathbb{B}}^\tau_\sigma$ we have that $D_n(\sigma, \sigma) = 0$ and for all $\tau \geq \sigma$,

(4.15)
$$D_n(\sigma, \tau+1)$$
$$\leq D_n(\sigma, \tau) + n^{-1} \sup_{\theta \in [\tau/n, (\tau+1)/n]} \|\widetilde{\mathbb{A}}_\tau \widetilde{\mathbb{B}}^{\tau-1}_\sigma - \mathbb{A}(\vec{y}(\theta, \rho), \theta)\mathbb{B}^\theta_{\sigma/n}\|.$$

As $(\vec{y}(\theta, \rho), \theta)$, $\theta \leq \theta_n$ and $n \geq n_4$ are in the set $\widehat{q}_+(\varepsilon)$ in which $(\vec{x}, \theta) \mapsto \mathbb{A}(\vec{x}, \theta)$ is bounded and Lipschitz continuous (for the operator norm), it follows that for some $c_5$ finite and all $n$,

(4.16) $$\sup_{\tau < \tau_n} \sup_{\theta \in [\tau/n, (\tau+1)/n]} \|\widetilde{\mathbb{A}}_\tau - \mathbb{A}(\vec{y}(\theta, \rho), \theta)\| \leq c_5 n^{-1}.$$

Further, with $\|\mathbb{A}(\vec{y}(\theta, \rho), \theta)\|$ bounded uniformly in $(\theta, \rho)$, we have from (4.4) the existence of $c_6$ finite, such that

$$\|\mathbb{B}^\theta_\zeta - \mathbb{B}^{\theta'}_{\zeta'}\| \le c_6(|\theta - \theta'| + |\zeta - \zeta'|), \tag{4.17}$$

for any $\rho \in [\varepsilon, 1/\varepsilon]$, $0 \le \zeta \le \theta \le 1 - \varepsilon$, and $0 \le \zeta' \le \theta' \le 1 - \varepsilon$. So, with $\widetilde{\mathbb{A}}_\tau$, $\mathbb{B}^{\tau/n}_{\sigma/n}$ and $\mathbb{A}(\cdot)$ uniformly bounded, by the Lipschitz properties (4.16) and (4.17) we have that

$$\|\widetilde{\mathbb{A}}_\tau \widetilde{\mathbb{B}}^{\tau-1}_\sigma - \mathbb{A}(\vec{y}(\theta), \theta) \mathbb{B}^\theta_{\sigma/n}\| \le c_7 D_n(\sigma, \tau) + c_8 n^{-1},$$

for some $c_7$, $c_8$ finite and all $n$, $\sigma \le \tau \le \tau_n$ and $\theta \in [\tau/n, (\tau+1)/n]$. Plugging this bound in (4.15) we have that $D_n(\sigma, \tau+1) \le (1 + c_7 n^{-1}) D_n(\sigma, \tau) + c_8 n^{-2}$, from which we deduce that for some $c_9$ finite and all $n$,

$$\max_{0 \le \sigma \le \tau \le \tau_n} D_n(\sigma, \tau) \le c_9 n^{-1}.$$

By (4.17), this yields the bound (4.14), hence completing the proof of (4.10) and that of the lemma. $\square$

4.2. *Asymptotic enumeration of the graph ensemble.* Here we show that the initial distribution of the Markov chain $\vec{z}(\cdot)$ of Section 3 is well approximated by a multivariate Gaussian law of mean $n\vec{y}(0)$ and positive definite covariance matrix $n\mathbb{Q}(0)$, with the rescaled mean $\vec{y}(0)$ and covariance $\mathbb{Q}(0)$ given by the initial condition of the corresponding ODE's, namely, (2.4) and (2.10), respectively.

LEMMA 4.4. *For $\vec{x} \in \mathbb{R}^d$ and a positive definite $d$-dimensional matrix $\mathbb{A}$, let $\mathsf{G}_d(\cdot|\vec{x}; \mathbb{A})$ denote the $d$-dimensional normal density of mean $\vec{x}$ and covariance $\mathbb{A}$. Further, let $\vec{z} = (z_1, z_2)$ denote the number of c-nodes of degree 1 and of degree strictly greater than 1 in a random graph from the $\mathcal{G}_l(n, \lfloor n\rho \rfloor)$ ensemble. Then, for any $\varepsilon > 0$ there exist finite, positive constants $\kappa_0$, $\kappa_1$, $\kappa_2$ and $\kappa_3$, such that for all $n$, $r$, and $\rho \in [\varepsilon, 1/\varepsilon]$,*

$$\|\mathbb{E}\vec{z} - n\vec{y}(0)\| \le \kappa_0, \tag{4.18}$$

$$\mathbb{P}\{\|\vec{z} - \mathbb{E}\vec{z}\| \ge r\} \le \kappa_1 e^{-r^2/\kappa_2 n}, \tag{4.19}$$

$$\sup_{\vec{u} \in \mathbb{R}^2} \sup_{x \in \mathbb{R}} \left| \mathbb{P}\{\vec{u} \cdot \vec{z} \le x\} - \int_{\vec{u} \cdot \vec{z} \le x} \mathsf{G}_2(\vec{z}|n\vec{y}(0); n\mathbb{Q}(0))\, d\vec{z} \right| \le \kappa_3 n^{-1/2}. \tag{4.20}$$

PROOF. Set $m = \lfloor n\rho \rfloor$ and $\gamma = l/\rho$. Recall that the description of the ensemble $\mathcal{G}_l(n, m)$ in Section 3.1 provides the following expression for the probability $\mathbb{P}(\vec{z})$ of having exactly $z_1$ c-nodes of degree 1 and $z_2$ c-nodes of degree strictly greater than 1:

$$\mathbb{P}(\vec{z}) = \frac{h(\vec{z}, 0)}{m^{nl}} = \frac{\mathbb{P}_\gamma\{\vec{S}_m = (z_1, z_2, nl)\}}{\mathbb{P}_\gamma\{S^{(3)}_m = nl\}}, \tag{4.21}$$

where $\vec{S}_m = \sum_{i=1}^{m} \vec{X}_i$ for $\vec{X}_i = (\mathbb{I}_{N_i=1}, \mathbb{I}_{N_i \geq 2}, N_i) \in \mathbb{Z}_+^3$ and $N_i$ that are i.i.d. Poisson($\gamma$) random variables. Consequently,

$$\mathbb{E} z_1 = m \mathbb{P}_\gamma \{N_1 = 1\} \frac{\mathbb{P}_\gamma \{S_{m-1}^{(3)} = nl - 1\}}{\mathbb{P}_\gamma \{S_m^{(3)} = nl\}},$$

$$\mathbb{E} z_2 = m - \mathbb{E} z_1 - m \mathbb{P}_\gamma \{N_1 = 0\} \frac{\mathbb{P}_\gamma \{S_{m-1}^{(3)} = nl\}}{\mathbb{P}_\gamma \{S_m^{(3)} = nl\}}.$$

With $|\rho n - m| \leq 1$ and $\vec{y}(0)$ of (2.4) such that $n\vec{y}(0) = n\rho(\mathbb{P}_\gamma\{N_1 = 1\}, \mathbb{P}_\gamma\{N_1 \geq 2\})$, we easily get (4.18) upon using the fact that $S_k^{(3)}$ is a Poisson($k\gamma$) random variable and the sequence $m|e(1-\frac{1}{m})^m - 1|$ is uniformly bounded.

By (4.18), in deriving (4.19) we may and shall replace $\mathbb{E}\vec{z}$ by $\frac{m}{\rho}\vec{y}(0) = (\mathbb{E}S_m^{(1)}, \mathbb{E}S_m^{(2)})$. In view of (4.21), the stated bound (4.19) is then merely

$$\mathbb{P}_\gamma\{|S_m^{(1)} - \mathbb{E}S_m^{(1)}|^2 + |S_m^{(2)} - \mathbb{E}S_m^{(2)}|^2 \geq r^2 | S_m^{(3)} = nl\} \leq \kappa_1 e^{-r^2/\kappa_2 n},$$

which is an immediate consequence of Hoeffding's inequality for the partial sums $(S_m^{(1)}, S_m^{(2)})$ and the uniform lower bound $\mathbb{P}_\gamma\{S_m^{(3)} = nl\} \geq cn^{-1/2}$ with $c > 0$ depending only on $\varepsilon$ and $l$.

Observe next that $\vec{X}_i$ are nondegenerate lattice random variables on $\mathbb{R}^3$, having minimal lattice $\mathbb{Z}^3$, finite moments of all orders and such that

$$\mathsf{cov}(\vec{X}_i) \equiv \mathbb{V} = \begin{pmatrix} p_1(1-p_1) & -p_1 p_{\geq 2} & p_1(1-\gamma) \\ -p_1 p_{\geq 2} & p_{\geq 2}(1-p_{\geq 2}) & \gamma - p_1 - \gamma p_{\geq 2} \\ p_1(1-\gamma) & \gamma - p_1 - \gamma p_{\geq 2} & \gamma \end{pmatrix},$$

with $p_1 = \mathbb{P}_\gamma(N_i = 1) = \gamma e^{-\gamma}$ and $p_{\geq 2} = \mathbb{P}_\gamma(N_i \geq 2) = 1 - e^{-\gamma} - \gamma e^{-\gamma}$. Thus, upon bounding $(1+\|\vec{u}\|^3)P_1(-\mathsf{G}_3(\cdot|\vec{0}, \mathbb{V}):\{\xi_\nu\})(\vec{u})$ for the correction term $P_1$ of [7], (7.19) (with $\{\xi_\nu\}$ denoting the cumulants of the law of $\vec{X}_1$), uniformly in $\gamma \in [\varepsilon', 1/\varepsilon']$ and $\vec{u} \in \mathbb{R}^3$, it follows from Corollary 22.3 of [7] (with $s = 3$ there), that for some finite $c = c(\varepsilon')$, any such $\gamma$, all $m$ and $\vec{z} \in \mathbb{Z}^2$,

$$(4.22) \quad |\mathbb{P}_\gamma\{\vec{S}_m = \vec{z}_e\} - \mathsf{G}_3(\vec{z}_e | m\vec{x}_e; m\mathbb{V})| \leq \frac{cm^{-2}}{1 + m^{-3/2}\|\vec{z}_e - m\vec{x}_e\|^3},$$

where $\vec{z}_e = (\vec{z}, nl)$ and $\vec{x}_e \equiv \rho^{-1}(\vec{y}(0), l) = \mathbb{E}_\gamma \vec{X}_1$. Applying the same argument for $S_m^{(3)} \in \mathbb{R}^1$, and possibly enlarging $c(\varepsilon')$ as needed we further have that

$$(4.23) \quad |\mathbb{P}_\gamma\{S_m^{(3)} = nl\} - \mathsf{G}_1(nl | m\gamma; mV_{33})| \leq cm^{-1}.$$

Next, summing the bound of (4.22) over $\vec{z} \in \mathbb{Z}^2$, we deduce that for some finite $c' = c'(\varepsilon)$ any $\gamma$ and $m$,

$$(4.24) \quad \sum_{\vec{z} \in \mathbb{Z}^2} |\mathbb{P}_\gamma\{\vec{S}_m = \vec{z}_e\} - \mathsf{G}_3(\vec{z}_e | m\vec{x}_e; m\mathbb{V})| \leq c'm^{-1}.$$

Further, $\mathbb{P}_\gamma\{S_m^{(3)} = nl\} = \sum_{\vec{z}} \mathbb{P}_\gamma\{\vec{S}_m = \vec{z}_e\}$, hence we get from (4.21) and the bounds of (4.23) and (4.24) that for some finite $\kappa = \kappa(\varepsilon')$ and any $\gamma$ and $m$,

$$(4.25) \quad \sum_{\vec{z} \in \mathbb{Z}^2} \left| \mathbb{P}(\vec{z}) - \frac{\mathsf{G}_3(\vec{z}_e|m\vec{x}_e m\mathbb{V})}{\mathsf{G}_1(nl|m\gamma; mV_{33})} \right| \leq \frac{cm^{-1} + c'm^{-1}}{\mathsf{G}_1(nl|m\gamma; mV_{33})} \leq \kappa n^{-1/2}$$

[with the rightmost inequality due to the uniform lower bound on $m^{1/2}\mathsf{G}_1(nl|m\gamma; mV_{33})$ for $|nl - m\gamma| \leq l/\varepsilon'$]. The ratio $\mathsf{G}_3(\cdots)/\mathsf{G}_1(\cdots)$ appearing in (4.25) is the conditional distribution of $(z_1, z_2)$, given $z_3 = nl$, under the (joint) law $\mathsf{G}_3(\cdots)$, which is thus a Gaussian distribution of mean $n'\vec{y}(0)$ and the positive definite covariance matrix $n'\widetilde{\mathbb{V}}$, with $n' \equiv m/\rho$ and the entries of the two-dimensional matrix $\widetilde{\mathbb{V}}$ given by $\widetilde{V}_{ij} = \rho[V_{ij} - V_{i3}V_{j3}/V_{33}]$. Upon substituting the expressions for $p_1$ and $p_{\geq 2}$, we see that $\widetilde{\mathbb{V}}$ coincides with $\mathbb{Q}(0)$ of (2.10).

So, it follows from (4.25) that

$$\sup_{\vec{u} \in \mathbb{R}^2} \sup_{x \in \mathbb{R}} \left| \mathbb{P}\{\vec{u} \cdot \vec{z} \leq x\} - \sum_{\vec{u} \cdot \vec{z} \leq x} \mathsf{G}_2(\vec{z}|n'\vec{y}(0); n'\mathbb{Q}(0)) \right| \leq \kappa n^{-1/2}.$$

We thus arrive at (4.20) upon observing first that

$$\sup_{h \leq 1} \sup_{\vec{u} \in \mathbb{R}^2} \sup_{x \in \mathbb{R}} \frac{1}{h} \left| \sum_{\vec{u} \cdot \vec{z} \leq x} \mathsf{G}_2(\vec{z}|h^{-2}\vec{y}(0); h^{-2}\mathbb{Q}(0)) - \int_{\vec{u} \cdot \vec{z} \leq x} \mathsf{G}_2(\vec{z}|h^{-2}\vec{y}(0); h^{-2}\mathbb{Q}(0)) \, d\vec{z} \right|$$

is uniformly bounded in $\gamma$ by the Euler–MacLaurin sum formula (cf. Theorem A.4.3 of [7] for the Schwartz function $\mathsf{G}_2(\cdot|\vec{0}; \mathbb{Q}(0))$, where the correction in $\Lambda_1(\vec{x})$ of [7], (A.4.20), to the Gaussian distribution is then at most $\kappa'h$ for a finite $\kappa'(\varepsilon)$, all $\vec{x} \in \mathbb{R}^2$ and $\gamma$), then noting that $\sqrt{n} \sup_{\vec{u}} \sup_x |G(\vec{u}, x; n) - G(\vec{u}, x; n')|$ is bounded in $\gamma$, $n$ and $|n' - n| \leq 1/\varepsilon$ for the Gaussian distribution function $G(\vec{u}, x; r) \equiv \int_{\vec{u} \cdot \vec{z} \leq x} \mathsf{G}_2(\vec{z}|r\vec{y}(0); r\mathbb{Q}(0)) \, d\vec{z}$. □

4.3. *Asymptotic transition probabilities.* We next prove an approximated formula for the transition probabilities $W^+_\tau(\Delta \vec{z}|\vec{z})$, that we often use in the sequel. This formula is valid throughout $\mathcal{Q}_+(\varepsilon) \equiv \mathcal{Q}(\varepsilon) \cap \{z_1 \geq 1\} \subseteq \mathbb{Z}^3$, where for each $\varepsilon > 0$,

$$\mathcal{Q}(\varepsilon) \equiv \{(\vec{z}, \tau) : -nl + n\varepsilon \leq z_1; n\varepsilon \leq z_2;$$
$$0 \leq \tau \leq n(1 - \varepsilon); n\varepsilon \leq (n - \tau)l - \max(z_1, 0) - 2z_2\}$$

is a finite subset of $\mathbb{Z}^3$. As many of our approximations involve the rescaled variables $\vec{x} \equiv n^{-1}\vec{z}$ and $\theta \equiv \tau/n$, we note in passing that if $(\vec{z}, \tau) \in \mathcal{Q}(\varepsilon)$, then necessarily $(\vec{x}, \theta)$ is in the set $\widehat{q}(\varepsilon)$ of Lemma 4.1 and if further $(\vec{z}, \tau) \in \mathcal{Q}_+(\varepsilon)$, then also $(\vec{x}, \theta) \in \widehat{q}_+(\varepsilon)$.

LEMMA 4.5. *For each $\theta \in [0,1)$ let $K_\theta : \mathbb{R}^2 \to \mathcal{K}_\theta$ denote the projection onto the convex set $\mathcal{K}_\theta \equiv \{\vec{x} \in \mathbb{R}^2_+ : x_1 + 2x_2 \leq l(1-\theta)\}$. Recall that each $\theta \in [0,1)$ and $\vec{x} \in \mathcal{K}_\theta$ specifies by (2.2) a well-defined probability vector $(\mathfrak{p}_0, \mathfrak{p}_1, \mathfrak{p}_2)$. For such $\theta, \vec{x}$ define the transition kernel*

$$(4.26) \qquad \widehat{W}_\theta(\Delta \vec{z} | \vec{x}) \equiv \binom{l-1}{q_0 - 1, q_1, q_2} \mathfrak{p}_0^{q_0-1} \mathfrak{p}_1^{q_1} \mathfrak{p}_2^{q_2},$$

*where $q_0 = -\Delta z_1 - \Delta z_2 \geq 1$, $q_1 = -\Delta z_2 \geq 0$, $q_2 = l + \Delta z_1 + 2\Delta z_2 \geq 0$. For any $\vec{x} \in \mathbb{R}^2$, set $\widehat{W}_\theta(\cdot | \vec{x}) \equiv \widehat{W}_\theta(\cdot | K_\theta(\vec{x}))$. That is, $\Delta z_1 = -1 - \tilde{q}_0 + q_1$ and $\Delta z_2 = -q_1$, with $(\tilde{q}_0, q_1, q_2)$ having the multinomial law of parameters $l-1$, $\mathfrak{p}_0$, $\mathfrak{p}_1$, $\mathfrak{p}_2$ that correspond to the projection of $\vec{x}$ onto $\mathcal{K}_\theta$.*

*Then, there exists a positive constant $C = C(l, \varepsilon)$, such that, for any $\rho \in [\varepsilon, 1/\varepsilon]$, $(\vec{z}, \tau) \in \mathcal{Q}_+(\varepsilon)$, $\Delta z_1 \in \{-l, \ldots, l-2\}$, $\Delta z_2 \in \{-(l-1), \ldots, 0\}$, and all $n$,*

$$|W^+_\tau(\Delta \vec{z} | \vec{z}) - \widehat{W}_{\tau/n}(\Delta \vec{z} | n^{-1} \vec{z})| \leq \frac{C}{n}.$$

PROOF. Following the notation of Lemma 3.1, for each $\vec{q} = (p_0, q_0, q_1, q_2) \in \mathcal{D}$, let

$$(4.27) \quad c_l(\vec{q}) = \binom{p_0 + q_1 + q_2}{p_0, q_1, q_2} \mathsf{coeff}[(e^{\mathbf{x}} - 1 - \mathbf{x})^{p_0} (e^{\mathbf{x}} - 1)^{q_1 + q_2}, \mathbf{x}^{l - q_0}],$$

and for $\vec{z} = (z_1, z_2)$ let

$$g_l(\vec{z}) = \sum_{\vec{q} \in \mathcal{D}} \binom{z_1 - 1}{q_0 - 1} \binom{z_2}{p_0 + q_1 + q_2} c_l(\vec{q}).$$

Using the identities $z_0 = z'_0 - q_0 - p_0$, $z'_1 - q_1 = z_1 - q_0$ and $z'_2 - q_2 = z_2 - p_0 - q_1 - q_2$ of (3.3), it follows after elementary algebra that $g_l(\vec{z})$ equals the sum over $\mathcal{D}$ in (3.2) times the term $\binom{m}{z_0, z_1, z_2} / \binom{m}{z'_0, z'_1, z'_2}$.

Next note that for any $\lambda > 0$, and integers $t, s \geq 1$,

$$(4.28) \qquad p_\lambda(t, s) = \mathsf{coeff}[(e^{\mathbf{x}} - 1 - \mathbf{x})^t, \mathbf{x}^s] \lambda^s (e^\lambda - 1 - \lambda)^{-t}$$

is precisely

$$p_\lambda(t, s) = \mathbb{P}_\lambda\left\{\sum_{i=1}^t N_i = s\right\},$$

where $\{N_i\}$ are i.i.d. random variables, with $\mathbb{P}_\lambda(N_1 = k) = \mathbb{P}(\overline{N}_\lambda = k | \overline{N}_\lambda \geq 2)$ and $\overline{N}_\lambda$ a Poisson random variable of parameter $\lambda > 0$. It is not hard to explicitly compute

$$f_1(\lambda) = \mathbb{E}_\lambda(N_1) = \frac{\lambda(e^\lambda - 1)}{e^\lambda - 1 - \lambda},$$

$$f_2(\lambda)^2 = \mathsf{Var}_\lambda(N_1) = \frac{\lambda}{(e^\lambda - 1 - \lambda)^2}[(e^\lambda - 1)^2 - \lambda^2 e^\lambda],$$

and the normalized $k$th moment $f_k(\lambda) = \mathbb{E}_\lambda (N_1 - f_1(\lambda))^k / f_2(\lambda)^k$, $k \geq 3$.

The behavior of $f_1$ was already considered in the proof of Lemma 4.1. Moreover $f_2 : \mathbb{R}_+ \to \mathbb{R}_+$ is bounded away from zero and infinity when $\lambda$ is bounded away from zero and infinity, respectively, resulting with $f_k(\lambda)$ that are also bounded away from infinity for each $k$.

Using (4.28) and writing explicitly the remaining terms in the expression (3.2), it is not hard to verify that

$$\text{(4.29)} \quad W^+_\tau(\Delta \vec{z} | \vec{z}) = n^{l-1} \binom{l(n-\tau)-1}{l-1}^{-1} \\ \times \lambda^{l+\Delta z_1}(e^\lambda - 1 - \lambda)^{\Delta z_2} \frac{p_\lambda(z'_2, (n-\tau-1)l - z'_1)}{p_\lambda(z_2, (n-\tau)l - z_1)} \hat{g}_l(\vec{z}),$$

where $z'_1 = z_1 + \Delta z_1$, $z'_2 = z_2 + \Delta z_2$ and $\hat{g}_l(\vec{z}) \equiv g_l(\vec{z})/n^{(l-1)}$.

Let $\xi \equiv ((n-\tau)l - z_1)/z_2$. Since $(n-\tau)l \geq z_1 + 2z_2$ we have that $\xi \geq 2$, and there exists a unique nonnegative solution of $f_1(\lambda) = \xi$. Further, as long as $(\vec{z}, \tau) \in \mathcal{Q}_+(\varepsilon)$ we get that $2 + (\varepsilon/\rho) \leq \xi \leq l/\varepsilon$ and hence $\varepsilon^2 \leq \lambda \leq l/\varepsilon$ (for $\rho \leq 1/\varepsilon$). We show in the sequel that this implies that there exists a positive constant $\tilde{C} = \tilde{C}(l, \varepsilon)$, such that, for any $\Delta z_1 \in \{-l, \ldots, l-2\}$ and $\Delta z_2 \in \{-(l-1), \ldots, 0\}$,

$$\text{(4.30)} \quad \left| \frac{p_\lambda(z_2 + \Delta z_2, (n-\tau)l - z_1 - l - \Delta z_1)}{p_\lambda(z_2, (n-\tau)l - z_1)} - 1 \right| \leq \frac{\tilde{C}}{n}.$$

Further, the positive term $\lambda^{l+\Delta z_1}(e^\lambda - 1 - \lambda)^{\Delta z_2}$ does not depend on $n$, whereas elementary calculus implies that

$$\text{(4.31)} \quad n^{l-1} \binom{(n-\tau)l-1}{l-1}^{-1} = \frac{(l-1)!}{[l(1-\theta)]^{l-1}} (1 + R_n),$$

where $|R_n| \leq \bar{C}(l)/(n\varepsilon)$ in $\mathcal{Q}_+(\varepsilon)$.

We turn to the asymptotic of $\hat{g}_l(\vec{z})$ for $(\vec{z}, \tau) \in \mathcal{Q}_+(\varepsilon)$. To this end, note that the condition $2p_0 + q_0 + q_1 + q_2 \leq l$ implies that the set $\mathcal{D}$ is at most of size $l^4$ and that the nonnegative coefficients $c_l(\vec{q})$ of (4.27) are bounded, uniformly in $\vec{q}$ by some $K = K(l) < \infty$ that is independent of $z_1$ and $z_2$ (hence independent of $n$). On $\mathcal{Q}_+(\varepsilon)$ the contribution to $\hat{g}_l(\vec{z})$ of the term indexed by $\vec{q}$ is at most $Kn^{-(l-1)}(nl)^{p_0+q_0+q_1+q_2-1}$. As $2p_0 + q_0 + q_1 + q_2 \leq l$, the sum over terms with either $p_0 > 0$ or $q_2 < l - q_0 - q_1$ is at most $Kl^{l-p_0+3}n^{-1}$.

Consider now $\vec{q}$ with $p_0 = 0$ and $q_2 = l - q_0 - q_1$, in which case $q_1 = -\Delta z_2$ and $q_0 = -\Delta z_1 - \Delta z_2 \geq 1$ are uniquely determined by $\Delta \vec{z}$. Note that $c_l(\vec{q}) = \binom{l-q_0}{q_1}$ for these choices of $p_0$ and $q_2$, resulting with

$$\text{(4.32)} \quad \hat{g}_l(\vec{z}) = n^{-(q_0-1)} \frac{(z_1-1)!}{(z_1-q_0)!} \frac{1}{(l-1)!} x_2^{l-q_0} \binom{l-1}{q_0-1, q_1, q_2} + \tilde{R}_n,$$

for some $|\tilde{R}_n| \leq \tilde{K}(l,\varepsilon)/n$. Since

$$x_1^{q_0-1}\left(1 - \frac{l^2}{n}\right) \leq n^{-(q_0-1)}\frac{(z_1-1)!}{(z_1-q_0)!} \leq x_1^{q_0-1},$$

replacing $n^{-(q_0-1)}(z_1-1)!/(z_1-q_0)!$ in (4.32) by $x_1^{q_0-1}$ and collecting together (4.29), (4.30), (4.31) and (4.32), results with the statement (4.26) of the lemma (note that $2q_1 + q_2 = l + \Delta z_1$).

We complete the proof of the lemma by showing that (4.30) is a consequence of a local CLT for the sum $S_k$ of i.i.d. lattice random variables $X_i = (N_i - \xi)/f_2(\lambda)$. Indeed, $X_i$ have zero mean, unit variance and for some finite $C_k$ we have that $|\mathbb{E}(X_1^k)| = |f_k(\lambda)| \leq C_k$ for all $(\vec{z},\tau) \in \mathcal{Q}_+(\varepsilon)$. Further, $p_\lambda(z_2,(n-\tau)l-z_1) = \mathbb{P}(S_k = 0)$ and $p_\lambda(z_2 + \Delta z_2, (n-\tau)l - z_1 - l - \Delta z_1) = \mathbb{P}(S_{k'} = \eta)$ for $k = z_2$, $k' - k = \Delta z_2 \in \{-(l-1),\ldots,0\}$ and $\eta = -(l + \Delta z_1 + \xi \Delta z_2)/f_2(\lambda)$. Note that $\eta$ is uniformly bounded by some $c_1 = c_1(l,\varepsilon)$ on $\mathcal{Q}_+(\varepsilon)$ and in the lattice of span $b = f_2(\lambda)^{-1}$ of possible values of $S_{k'}$. Thus, for some finite $c_2 = c_2(l,\varepsilon)$, all $\eta$ and $k'$ as above, we have by Theorem 5.4 and (5.27) of [21] that

$$\left| f_2(\lambda)\sqrt{k'}\mathbb{P}(S_{k'} = \eta) - \phi\left(\frac{\eta}{\sqrt{k'}}\right) + \frac{f_3(\lambda)}{6\sqrt{k'}}\phi^{(3)}\left(\frac{\eta}{\sqrt{k'}}\right) \right| \leq \frac{c_2}{k'},$$

where $\phi(u) = e^{-u^2/2}/\sqrt{2\pi}$ and $\phi^{(3)}(u)$ denotes its third derivative. The same applies for $k$ and $\eta = 0$, yielding that

$$|f_2(\lambda)\sqrt{k}\mathbb{P}(S_k = 0) - \phi(0)| \leq \frac{c_2}{k}.$$

In particular, with $k \geq n\varepsilon$, we see that $\mathbb{P}(S_k = 0) \geq c_3/\sqrt{n}$ for some $c_3 > 0$ and all $n \geq n_0$, both depending only upon $l$ and $\varepsilon$. As $\phi(u)$ is an even function with uniformly bounded derivatives of any order, $k, k' \geq \varepsilon n$, $|\eta| \leq c_1$ and $|k - k'| \leq l$, it follows that for some finite $c_4 = c_4(l,\varepsilon)$,

$$\left| \frac{\sqrt{k}}{\sqrt{k'}}\phi\left(\frac{\eta}{\sqrt{k'}}\right) - \phi(0) - \frac{\sqrt{k}f_3(\lambda)}{6k'}\phi^{(3)}\left(\frac{\eta}{\sqrt{k'}}\right) \right| \leq \frac{c_4}{n},$$

from which (4.30) now directly follows. □

We often rely on the following regularity property of $(\vec{x},\theta) \mapsto \widehat{W}_\theta(\cdot|\vec{x})$ for the transition kernels of (4.26).

LEMMA 4.6. *With $\|\cdot\|_{\mathrm{TV}}$ denoting the total variation norm and $\|\cdot\|$ the Euclidean norm in $\mathbb{R}^2$, there exist positive constants $L = L(l,\varepsilon)$ such that for any $\theta,\theta' \in [0, 1-\varepsilon]$ and $\vec{x},\vec{x}' \in \mathbb{R}^2$,*

(4.33) $$\|\widehat{W}_{\theta'}(\cdot|\vec{x}') - \widehat{W}_\theta(\cdot|\vec{x})\|_{\mathrm{TV}} \leq L(\|\vec{x}' - \vec{x}\| + |\theta' - \theta|).$$

PROOF. With $(\vec{x}, \theta) \mapsto K_\theta(\vec{x})$ Lipschitz continuous, given that one finite set supports the kernels $\widehat{W}_\theta(\cdot|\vec{x})$ for all $(\vec{x}, \theta)$ and that $\widehat{W}_\theta(\Delta \vec{z}|\vec{x})$ of (4.26) is a smooth function of $(\mathfrak{p}_0, \mathfrak{p}_1, \mathfrak{p}_2)$ for $\vec{x} \in \mathcal{K}_\theta$, we get (4.33) out of the Lipschitz continuity of $(\mathfrak{p}_0, \mathfrak{p}_1, \mathfrak{p}_2)$ on $\widehat{q}(\varepsilon)$, proved in Lemma 4.1. □

4.4. *Absence of small cores.* A considerable simplification of our analysis comes from the observation that a typical large random hypergraph does not have a nonempty core of size below a certain threshold. For the convenience of the reader, we next adapt a result of [30] (and its proof) to the context of our graph ensemble.

LEMMA 4.7. *A subset of v-nodes of a hypergraph is called a stopping set if the restriction of the hypergraph to this subset has no c-node of degree 1. For $l \geq 3$ and any $\varepsilon > 0$ there exist $\kappa(l, \varepsilon) > 0$ and $C(l, \varepsilon)$ finite such that for any $m \geq \varepsilon n$ the probability that a random hypergraph from the ensemble $\mathcal{G}_l(n, m)$ has a stopping set of less than $m\kappa(l, \varepsilon)$ v-nodes is at most $C(l, \varepsilon) m^{1 - l/2}$.*

REMARK 4.8. Since the core is the stopping set including the maximal number of v-nodes, the lemma implies that for $m \geq \varepsilon n$ the probability that a random hypergraph from the ensemble $\mathcal{G}_l(n, m)$ has a nonempty core of size less than $m\kappa(l, \varepsilon)$ is at most $C(l, \varepsilon) m^{1 - l/2}$. With $n \leq m/\varepsilon$, upon changing $\kappa$ to $\kappa/\varepsilon$ and increasing $C$ as needed, it further follows that the probability of having a nonempty core with less than $n\kappa$ v-nodes is at most $C n^{1 - l/2}$.

PROOF. Let $N(s, r)$ denote the number of stopping sets in our random hypergraph which involve exactly $s$ v-nodes and $r$ c-nodes. Then, necessarily $r \leq \lfloor ls/2 \rfloor$ and

$$\mathbb{E} N(s, r) = \binom{n}{s} \binom{m}{r} \frac{1}{m^{sl}} \mathsf{coeff}[(e^{\mathbf{x}} - 1 - \mathbf{x})^r, \mathbf{x}^{sl}] (sl)!$$

(multiply the number of sets of $s$ v-nodes and $r$ c-nodes by the probability that such a set forms a stopping set, with $\mathsf{coeff}[(e^{\mathbf{x}} - 1 - \mathbf{x})^r, \mathbf{x}^{sl}](sl)!$ counting the number of ways of connecting the $s$ v-nodes to these $r$ c-nodes so as to form a stopping set, while $m^{sl}$ is the total number of ways of connecting the $s$ v-nodes in our graph ensemble). It is easy to see that for any integers $r, t \geq 1$,

$$\mathsf{coeff}[(e^{\mathbf{x}} - 1 - \mathbf{x})^r, \mathbf{x}^t] \leq (e^{\mathbf{x}} - 1 - \mathbf{x})^r |_{\mathbf{x}=1} \leq 1.$$

Hence, for some $\zeta = \zeta(l, \varepsilon)$ finite, any $m \geq \varepsilon n$, $sl \leq m$ and $r \leq \lfloor ls/2 \rfloor$,

$$\mathbb{E} N(s, r) \leq \binom{n}{s} \binom{m}{r} \frac{(sl)!}{m^{sl}} \leq \frac{n^s}{s!} \frac{m^{\lfloor sl/2 \rfloor}}{\lfloor sl/2 \rfloor!} \frac{(sl)!}{m^{sl}}$$

$$\leq \frac{n^s}{s!} \left( \frac{sl}{m} \right)^{\lceil sl/2 \rceil} \leq \left[ \zeta \left( \frac{s}{m} \right)^{l/2 - 1} \right]^s.$$

Thus, fixing $0 < \kappa < 1/l$ (so $sl \leq m$ whenever $s \leq \kappa m$), for $l \geq 3$, the probability that a random hypergraph from the ensemble $\mathcal{G}_l(n,m)$ has a stopping set of size at most $m\kappa$ is bounded above by

$$\mathbb{E}\left[\sum_{s=1}^{m\kappa}\sum_{r=1}^{\lfloor ls/2 \rfloor} N(s,r)\right] \leq \zeta m^{1-l/2}\sum_{s=1}^{\infty} sl(\zeta\kappa^{l/2-1})^{s-1} \leq 4\zeta l m^{1-l/2},$$

provided $\zeta\kappa^{l/2-1} \leq 1/2$. □

**5. Auxiliary processes and proof of Proposition 2.1.** In this section we provide relations between two auxiliary inhomogeneous $\mathbb{Z}^2$-valued Markov processes whose distributions are denoted, respectively, as $\mathbb{P}_{n,\rho}(\cdot)$ and $\widehat{\mathbb{P}}_{n,\rho}(\cdot)$. In both cases, we denote the process as $\{\vec{z}(\tau) = (z_1(\tau), z_2(\tau)), 0 \leq \tau \leq n\}$, and use for both the same initial condition

$$\mathbb{P}_{n,\rho}(\vec{z}(0) = \vec{z}) = \widehat{\mathbb{P}}_{n,\rho}(\vec{z}(0) = \vec{z}) = \mathbb{P}_{\mathcal{G}_l(n,m)}(\vec{z}(G) = \vec{z}) = \frac{h(\vec{z},0)}{m^{nl}},$$

if $\vec{z} \in \mathbb{Z}_+^2$ is such that $z_1 + 2z_2 \leq nl$, and $\mathbb{P}_{n,\rho}(\vec{z}(0) = \vec{z}) = \widehat{\mathbb{P}}_{n,\rho}(\vec{z}(0) = \vec{z}) = 0$ otherwise. Here $\mathbb{P}_{\mathcal{G}_l(n,m)}(\cdot)$ is the uniform distribution on the graph ensemble $\mathcal{G}_l(n,m)$ and $m \equiv \lfloor n\rho \rfloor$.

Turning to specify the transition kernels, recall the triangles $\mathcal{K}_\theta \equiv \{\vec{x} \in \mathbb{R}_+^2 : x_1 + 2x_2 \leq l(1-\theta)\}$, $\theta \in [0,1)$, and set

$$W_\tau(\Delta\vec{z}|\vec{z}) = \begin{cases} W_\tau^+(\Delta\vec{z}|\vec{z}), & \text{if } z_1 \geq 1, n^{-1}\vec{z} \in \mathcal{K}_{\tau/n}, \\ \widehat{W}_{\tau/n}(\Delta\vec{z}|n^{-1}\vec{z}), & \text{otherwise,} \end{cases}$$

for $W_\tau^+(\cdot|\cdot)$ of (3.2) and the simpler kernel $\widehat{W}_\theta(\cdot|\cdot)$ of (4.26). The transition probabilities are then

(5.1) $\quad \mathbb{P}_{n,\rho}(\vec{z}(\tau+1) = \vec{z} + \Delta\vec{z}|\vec{z}(\tau) = \vec{z}) = W_\tau(\Delta\vec{z}|\vec{z}),$

(5.2) $\quad \widehat{\mathbb{P}}_{n,\rho}(\vec{z}(\tau+1) = \vec{z} + \Delta\vec{z}|\vec{z}(\tau) = \vec{z}) = \widehat{W}_{\tau/n}(\Delta\vec{z}|n^{-1}\vec{z}),$

for $\tau = 0, 1, \ldots, n-1$. While the Markov process of Lemma 3.1 describing the evolution under the decimation algorithm has $n^{-1}\vec{z}(\tau) \in \mathcal{K}_{\tau/n}$, this is not necessarily the case for the two auxiliary processes we consider here. Nevertheless, the Markov process of Lemma 3.1 coincides with the one associated with $\mathbb{P}_{n,\rho}(\cdot)$ up to the first time $\tau$ at which $z_1(\tau) = 0$, that is, when the decimation algorithm terminates at the core of the hypergraph.

We next provide a coupling that keeps the process of distribution $\mathbb{P}_{n,\rho}(\cdot)$ "very close" to its "approximation" by the process of distribution $\widehat{\mathbb{P}}_{n,\rho}(\cdot)$ as long as the former belongs to $\mathcal{Q}(\eta)$ for some $\eta > 0$. We shall see in Corollary 5.4 that up to an exponentially small probability (as $n \to \infty$), this is indeed the case for $\tau \leq (1-\varepsilon)n$, allowing us to focus on the properties of the simpler distribution $\widehat{\mathbb{P}}_{n,\rho}(\cdot)$.

LEMMA 5.1. *There exist finite $C_* = C_*(l, \varepsilon)$ and positive $\lambda_* = \lambda_*(l, \varepsilon)$, and a coupling between $\{\vec{z}(\tau)\} \stackrel{d}{=} \mathbb{P}_{n,\rho}(\cdot)$ and $\{\vec{z}'(\tau)\} \stackrel{d}{=} \widehat{\mathbb{P}}_{n,\rho}(\cdot)$, such that for any $n$, $\rho \in [\varepsilon, 1/\varepsilon]$ and $r > 0$,*

$$\text{(5.3)} \qquad \mathbb{P}\left\{\sup_{\tau \leq \tau_*} \|\vec{z}(\tau) - \vec{z}'(\tau)\| > r\right\} \leq C_* e^{-\lambda_* r},$$

*where $\tau_* \leq n$ denotes the first time such that $(\vec{z}(\tau_*), \tau_*) \notin \mathcal{Q}(\varepsilon)$.*

PROOF. To construct the coupling between the two processes, start with $\vec{z}'(0) = \vec{z}(0)$, which is possible since $\vec{z}(0)$ and $\vec{z}'(0)$ are identically distributed. Then, for $\tau = 0, 1, \ldots, n-1$, with $\vec{z}(\tau) = \vec{z}$ and $\vec{z}'(\tau) = \vec{z}'$, set $\vec{z}(\tau+1) = \vec{z} + \Delta\vec{z}$ and $\vec{z}'(\tau+1) = \vec{z}' + \Delta\vec{z}'$, where the joint distribution (coupling) of $(\Delta\vec{z}, \Delta\vec{z}')$ is chosen such that

$$\text{(5.4)} \qquad \mathbb{P}(\Delta\vec{z} \neq \Delta\vec{z}' | \vec{z}, \vec{z}') = \|W_\tau(\cdot | \vec{z}) - \widehat{W}_{\tau/n}(\cdot | n^{-1}\vec{z}')\|_{\text{TV}}.$$

Clearly, it suffices to show that $\Delta_n(\lambda_*) \leq C_*$ for some $\lambda_* > 0$ and $C_* < \infty$ that depend only on $l$ and $\varepsilon$, where

$$Z(\tau) \equiv \sup_{\sigma \leq \tau \wedge \tau_*} \|\vec{z}(\sigma) - \vec{z}'(\sigma)\|, \qquad \Delta_\tau(\lambda) \equiv \mathbb{E}[e^{\lambda Z(\tau)}],$$

for $\tau = \{0, \ldots, n\}$ and $\lambda \geq 0$. To this end, note first that by our definition of $W_\tau(\cdot | \vec{z})$, we have from Lemma 4.5 that for some finite $\tilde{c} = \tilde{c}(l, \varepsilon)$, any $(\vec{z}, \tau) \in \mathcal{Q}(\varepsilon)$, and all $n$,

$$\text{(5.5)} \qquad \|W_\tau(\cdot | \vec{z}) - \widehat{W}_{\tau/n}(\cdot | n^{-1}\vec{z})\|_{\text{TV}} \leq \frac{\tilde{c}}{n}$$

[since the kernels $W_\tau(\cdot | \vec{z})$ and $\widehat{W}_{\tau/n}(\cdot | n^{-1}\vec{z})$ are nonzero for at most $2l^2$ points]. Further, with $\|\Delta\vec{z}\| \leq 2l$ and $\|\Delta\vec{z}'\| \leq 2l$, we have that for any $0 \leq \lambda \leq 1/(4l)$ (so $e^{4l\lambda} \leq 1 + 8l\lambda$), $\sigma = 0, 1, \ldots, n-1$ and realizations of the two processes,

$$e^{\lambda Z(\sigma+1)} \leq \{1 + 8l\lambda \mathbb{I}_{\{\Delta\vec{z}(\sigma) \neq \Delta\vec{z}'(\sigma), \sigma < \tau_*\}}\} e^{\lambda Z(\sigma)}.$$

As $\tau_*$ is a stopping time and our coupling satisfies (5.4), upon considering the expectation of the preceding inequality we get that

$$\text{(5.6)} \quad \begin{aligned} \Delta_{\sigma+1}(\lambda) &\leq \Delta_\sigma(\lambda) \\ &\quad + 8l\lambda \mathbb{E}\{\|W_\sigma(\cdot | \vec{z}(\sigma)) - \widehat{W}_{\sigma/n}(\cdot | n^{-1}\vec{z}'(\sigma))\|_{\text{TV}} \mathbb{I}_{\sigma < \tau_*} e^{\lambda Z(\sigma)}\}. \end{aligned}$$

Recall that as long as $(\vec{z}(\sigma), \sigma) \in \mathcal{Q}(\varepsilon)$, by (5.5) and Lemma 4.6

$$\text{(5.7)} \quad \|W_\sigma(\cdot | \vec{z}(\sigma)) - \widehat{W}_{\sigma/n}(\cdot | n^{-1}\vec{z}'(\sigma))\|_{\text{TV}} \leq \frac{\tilde{c}}{n} + \frac{L}{n}\|\vec{z}(\sigma) - \vec{z}'(\sigma)\|.$$

Since $\|\vec{z}(\sigma) - \vec{z}'(\sigma)\|\mathbb{I}_{\sigma<\tau_*} \leq Z(\sigma)$, combining the bounds of (5.6) and (5.7), we deduce that

$$\Delta_{\sigma+1}(\lambda) \leq [1 + 8l\widetilde{c}n^{-1}\lambda]\mathbb{E}\{(1 + n^{-1}8lL\lambda Z(\sigma))e^{\lambda Z(\sigma)}\}$$
(5.8)
$$\leq [1 + 8l\widetilde{c}n^{-1}\lambda]\Delta_\sigma(\lambda(1 + 8lLn^{-1})).$$

Since $\Delta_0(\lambda) = 1$, taking $\lambda = \lambda_* = \exp(-8lL)/(4l) \leq 1/(4l)$, and applying the inequality (5.8) for the monotone increasing sequence $\{\lambda_\sigma, \sigma \geq 0\}$ with $\lambda_0 = \lambda_*$ and $\lambda_{\sigma+1} = \lambda_\sigma(1 + 8lLn^{-1})$, such that $\lambda_n = \lambda_*(1 + 8lL/n)^n \leq 1/(4l)$, we get that

$$\Delta_n(\lambda_*) \leq \prod_{\sigma=0}^{n-1}(1 + 8l\widetilde{c}n^{-1}\lambda_\sigma) \leq \exp\left\{8l\widetilde{c}n^{-1}\sum_{\sigma=0}^{n-1}\lambda_\sigma\right\} \leq \exp\{8l\widetilde{c}\lambda_n\} \leq \exp\{2\widetilde{c}\},$$

completing the proof of the lemma. □

We turn to establish some of the asymptotic (in $n \to \infty$) properties of our approximating processes [of distribution $\widehat{\mathbb{P}}_{n,\rho}(\cdot)$].

LEMMA 5.2. *For any $l \geq 3$ and $\varepsilon > 0$ there exist positive, finite constants $\eta \leq \varepsilon$, and $C_0, C_1, C_2, C_3$, such that, for any $n$, $\rho \in [\varepsilon, 1/\varepsilon]$ and $\tau \in \{0, \ldots, \lfloor n(1-\varepsilon) \rfloor\}$,*

(a) *$\vec{z}(\tau)$ is exponentially concentrated around its mean*

(5.9) $$\widehat{\mathbb{P}}_{n,\rho}\{\|\vec{z}(\tau) - \mathbb{E}\vec{z}(\tau)\| \geq r\} \leq 4e^{-r^2/C_0 n}.$$

(b) *$\vec{z}(\tau)$ is close to the solution of the ODE (2.5),*

(5.10) $$\mathbb{E}\|\vec{z}(\tau) - n\vec{y}(\tau/n)\| \leq C_1\sqrt{n \log n}.$$

(c) *$(\vec{z}(\tau), \tau) \in \mathcal{Q}(\eta)$ with high probability; more precisely,*

(5.11) $$\widehat{\mathbb{P}}_{n,\rho}\{(\vec{z}(\tau), \tau) \notin \mathcal{Q}(\eta)\} \leq C_2 e^{-C_3 n}.$$

PROOF. (a) For $\tau = 0$, upon taking $C_0$ large enough, this is an immediate consequence of (4.19). Turning to the general case, applying the Azuma–Hoeffding inequality for Doob's martingale

$$Z(\sigma) = \mathbb{E}[\vec{z}(\tau)|\vec{z}(0), \ldots, \vec{z}(\sigma)], \qquad \sigma \in \{0, \ldots, \tau\},$$

we see that for some $c_0 = c_0(\varepsilon)$ finite, any $n, r > 0$, $\tau = 1, \ldots, n(1-\varepsilon)$ and $\rho \in [\varepsilon, 1/\varepsilon]$,

(5.12) $$\widehat{\mathbb{P}}_{n,\rho}\{\|\vec{z}(\tau) - \mathbb{E}[\vec{z}(\tau)|\vec{z}(0)]\| \geq r\} \leq 4\exp(-r^2/(2c_0^2\tau)),$$

provided $\operatorname{ess\,sup} \|Z(\sigma) - Z(\sigma - 1)\| \leq c_0$ for all $1 \leq \sigma \leq \tau$. To this end, with $\vec{z}(\cdot)$ a Markov process, we have the bound

$$
\begin{aligned}
&\operatorname{ess\,sup} \|Z(\sigma) - Z(\sigma - 1)\| \\
&\quad \leq \sup_{\vec{z}^{(1)}, \vec{z}^{(2)}} \|\mathbb{E}[\vec{z}(\tau)|\vec{z}(\sigma) = \vec{z}^{(1)}] - \mathbb{E}[\vec{z}(\tau)|\vec{z}(\sigma) = \vec{z}^{(2)}]\|,
\end{aligned}
\tag{5.13}
$$

where the preceding supremum is over all $\vec{z}^{(1)}$, $\vec{z}^{(2)}$ such that some trajectories $\{\vec{z}(0) \ldots \vec{z}(\sigma - 1), \vec{z}(\sigma) = \vec{z}^{(1)}\}$ and $\{\vec{z}(0) \ldots \vec{z}(\sigma - 1), \vec{z}(\sigma) = \vec{z}^{(2)}\}$ are both of positive probability. In particular, $\|\vec{z}^{(1)} - \vec{z}^{(2)}\| \leq 4l$. Fixing such $\vec{z}^{(1)}$ and $\vec{z}^{(2)}$, let $\vec{z}^{(1)}(\nu)$ and $\vec{z}^{(2)}(\nu)$ denote the realizations of two Markov processes of same transition kernels $\widehat{W}_\theta(\cdot|\cdot)$, starting at $\vec{z}^{(1)}(\sigma) = \vec{z}^{(1)}$ and $\vec{z}^{(2)}(\sigma) = \vec{z}^{(2)}$, respectively, where for $\nu = \sigma, \ldots, \tau - 1$ the joint distribution (coupling) of $\Delta \vec{z}^{(1)}(\nu) \equiv \vec{z}^{(1)}(\nu + 1) - \vec{z}^{(1)}(\nu)$ and $\Delta \vec{z}^{(2)}(\nu) \equiv \vec{z}^{(2)}(\nu + 1) - \vec{z}^{(2)}(\nu)$ is chosen such that

$$
\begin{aligned}
&\mathbb{P}(\Delta \vec{z}^{(1)}(\nu) \neq \Delta \vec{z}^{(2)}(\nu)|\vec{z}^{(1)}(\nu), \vec{z}^{(2)}(\nu)) \\
&\quad = \|\widehat{W}_{\nu/n}(\cdot|n^{-1}\vec{z}^{(1)}(\nu)) - \widehat{W}_{\nu/n}(\cdot|n^{-1}\vec{z}^{(2)}(\nu))\|_{\mathrm{TV}}.
\end{aligned}
$$

With $\Delta(\nu) \equiv \mathbb{E}\|\vec{z}^{(1)}(\nu) - \vec{z}^{(1)}(\nu)\|$, the right-hand side of (5.13) is upper-bounded by the supremum of $\Delta(\tau)$ over all possible pairs of initial conditions such that $\Delta(\sigma) = \|\vec{z}^{(1)} - \vec{z}^{(2)}\| \leq 4l$. Further, due to the Markov property of $\vec{z}$ and the preceding coupling, for $\sigma \leq \nu < \tau$ we have by (4.33) that

$$
\begin{aligned}
\Delta(\nu + 1) &\leq \mathbb{E}\|\vec{z}^{(1)}(\nu) - \vec{z}^{(2)}(\nu)\| + \mathbb{E}\{\mathbb{E}[\|\Delta \vec{z}^{(1)}(\nu) - \Delta \vec{z}^{(2)}(\nu)\| | \vec{z}^{(1)}(\nu), \vec{z}^{(2)}(\nu)]\} \\
&\leq \Delta(\nu) + 4l\mathbb{E}\{\|\widehat{W}_{\nu/n}(\cdot|n^{-1}\vec{z}^{(1)}(\nu)) - \widehat{W}_{\nu/n}(\cdot|n^{-1}\vec{z}^{(2)}(\nu))\|_{\mathrm{TV}}\} \\
&\leq \left(1 + \frac{4lL}{n}\right)\Delta(\nu).
\end{aligned}
$$

With $\tau \leq n$, it thus follows that $\Delta(\tau) \leq \exp(4lL)\Delta(\sigma) \leq 4l\exp(4lL) =: c_0$, as claimed.

Further, the preceding argument shows that $\psi(\vec{z}) \equiv \mathbb{E}[\vec{z}(\tau)|\vec{z}(0) = \vec{z}]$ is a uniformly Lipschitz continuous function of $\vec{z}$, of Lipschitz constant $\|\psi\|_{\mathrm{L}} = \exp(4lL)$ that is independent of $\tau$, $n$ and $\rho$. Hence, from (4.19) we have that

$$\mathbb{P}\{\|\psi(\vec{z}(0)) - \psi(\mathbb{E}\vec{z}(0))\| \geq r\|\psi\|_{\mathrm{L}}\} \leq \mathbb{P}\{\|\vec{z}(0) - \mathbb{E}\vec{z}(0)\| \geq r\} \leq \kappa_1 e^{-r^2/\kappa_2 n}.$$

Integrating this over $r \geq 0$, we have that $\|\mathbb{E}\psi(\vec{z}(0)) - \psi(\mathbb{E}\vec{z}(0))\| \leq c\sqrt{n}$ for some finite constant $c$ depending only on $\varepsilon$ and $l$, yielding that

$$
\begin{aligned}
&\widehat{\mathbb{P}}_{n,\rho}\{\|\mathbb{E}[\vec{z}(\tau)|\vec{z}(0)] - \mathbb{E}[\vec{z}(\tau)]\| \geq r\} \\
&\quad = \mathbb{P}\{\|\psi(\vec{z}(0)) - \mathbb{E}\psi(\vec{z}(0))\| \geq r\} \leq C_1' e^{-r^2/c_2' n},
\end{aligned}
$$

for some $C_1'$ and $c_2'$ which depend only on $\varepsilon$ and $l$, which, together with (5.12), concludes the proof of (5.9).

(b) Since $\|\vec{z}(\tau)\| \leq 2nl$, choosing $r = \sqrt{C_0 n \log n}$ in (5.9) we find that

(5.14) $$\mathbb{E}\|\vec{z}(\tau) - \mathbb{E}\vec{z}(\tau)\| \leq c_1 \sqrt{n \log n},$$

for some finite $c_1(\varepsilon)$. Denote by $\Delta_m(\tau) \equiv \|\mathbb{E}\vec{z}(\tau) - n\vec{y}(\tau/n)\|$ the error made in replacing the expectation of the process $\vec{z}(\tau)$ of distribution $\widehat{\mathbb{P}}_{n,\rho}(\cdot)$ with the (rescaled) solution of the ODE. Then, fixing $\tau \leq n(1-\varepsilon)$, we have by the Markov property of $\vec{z}(\cdot)$ that

$$\Delta_m(\tau+1)$$
$$= \|\mathbb{E}\vec{z}(\tau) - n\vec{y}(\tau/n) + \mathbb{E}\{\mathbb{E}[\Delta\vec{z}(\tau)|\vec{z}(\tau)]\} - n[\vec{y}(\tau/n + 1/n) - \vec{y}(\tau/n)]\|.$$

Recall that for $\theta \leq 1 - \varepsilon$,

$$\frac{d\vec{y}}{d\theta} = \vec{F}(\vec{y}, \theta) = \sum_{\Delta\vec{z}} \widehat{W}_\theta(\Delta\vec{z}|\vec{y}) \Delta\vec{z},$$

so by the triangle inequality we get that

$$\Delta_m(\tau+1) \leq \Delta_m(\tau) + \left\|\mathbb{E}\left\{\sum_{\Delta\vec{z}}[\widehat{W}_{\tau/n}(\Delta\vec{z}|n^{-1}\vec{z}(\tau)) - \widehat{W}_{\tau/n}(\Delta\vec{z}|n^{-1}\mathbb{E}\vec{z}(\tau))]\Delta\vec{z}\right\}\right\|$$
$$+ \left\|\mathbb{E}\left\{\sum_{\Delta\vec{z}}[\widehat{W}_{\tau/n}(\Delta\vec{z}|n^{-1}\mathbb{E}\vec{z}(\tau)) - \widehat{W}_{\tau/n}(\Delta\vec{z}|\vec{y}(\tau/n))]\Delta\vec{z}\right\}\right\|$$
$$+ n\left\|\int_{\tau/n}^{\tau/n+1/n}[\vec{F}(\vec{y}(\theta), \theta) - \vec{F}(\vec{y}(\tau/n), \tau/n)]\,d\theta\right\|$$
$$\equiv \Delta_m(\tau) + \delta_m^{(0)}(\tau) + \delta_m^{(1)}(\tau) + \delta_m^{(2)}(\tau).$$

Recall as in (4.11) that $\delta_m^{(2)}(\tau) \leq C_* n^{-1}$. Since $\|\Delta\vec{z}\| \leq 4l$, we have by (4.33) that

$$\delta_m^{(1)}(\tau) \leq 4l\mathbb{E}\|\widehat{W}_{\tau/n}(\Delta\vec{z}|n^{-1}\mathbb{E}\vec{z}(\tau)) - \widehat{W}_{\tau/n}(\Delta\vec{z}|\vec{y}(\tau/n))\|_{\mathrm{TV}}$$
$$\leq \frac{4lL}{n}\|\mathbb{E}\vec{z}(\tau) - n\vec{y}(\tau/n)\| = \frac{4lL}{n}\Delta_m(\tau).$$

Similarly, by (4.33) and (5.14), for some $c_2 = c_2(\varepsilon)$ finite,

$$\delta_m^{(0)}(\tau) \leq \frac{4lL}{n}\mathbb{E}\|\vec{z}(\tau) - \mathbb{E}\vec{z}(\tau)\| \leq c_2\sqrt{\frac{\log n}{n}},$$

so putting these estimates together, we obtain the inequality

$$\Delta_m(\tau+1) \leq \left(1 + \frac{4lL}{n}\right)\Delta_m(\tau) + c_3\sqrt{\frac{\log n}{n}}.$$

Further, recall (4.18) of Lemma 4.4 that $\Delta_m(0)$ is bounded in $n$ and $m = \lfloor \rho n \rfloor$, provided $\rho \in [\varepsilon, 1/\varepsilon]$. Thus, we easily get (5.10) upon applying the preceding recursion for $\tau = 0, \ldots, n-1$.

(c) In the course of proving part (a) of Proposition 4.2 we have seen that there exists $\eta = \eta(\varepsilon, l) > 0$ such that if $\rho \in [\varepsilon, 1/\varepsilon]$ and $\theta \leq (1 - \varepsilon)$, then $y_1(\theta) \geq -l + 2\eta$, $y_2(\theta) \geq 2\eta$ and $(1-\theta)l - \max(y_1(\theta), 0) - 2y_2(\theta) \geq 2\eta$. Consequently, taking $\eta \leq \varepsilon$, for such $\rho$ and $\tau \in \{0, \ldots, \lfloor n(1-\varepsilon) \rfloor\}$, if $\|\vec{z}(\tau) - n\vec{y}(\tau/n)\| \leq n\eta/3$, then clearly $(\vec{z}(\tau), \tau) \in \mathcal{Q}(\eta)$. We thus get (5.11) upon considering (5.9) and (5.10) for $r = n\eta/6$ and $n$ such that $C_1\sqrt{n \log n} \leq n\eta/6$. □

The first consequence of Lemma 5.2 is the existence of "critical time window." That is, for $\rho$ near $\rho_c$ a typical trajectory $\{\vec{z}(\tau); 0 \leq \tau \leq (1-\varepsilon)n\}$ does not traverse the $z_1 = 0$ plane if $\tau$ is not near $n\theta_c$.

COROLLARY 5.3. *Fixing $\beta \in (3/4, 1)$, $\beta' < 2\beta - 1$ and $\varepsilon > 0$, let $I_n \equiv [0, n\theta_c - n^\beta] \cup [n\theta_c + n^\beta, n(1-\varepsilon)]$. Then, for some $C_4$ finite, $\eta$ positive, all $n$ and $|\rho - \rho_c| \leq n^{\beta'-1}$,*

$$\widehat{\mathbb{P}}_{n,\rho}\left\{\min_{\tau \in I_n} z_1(\tau) \leq n^{\beta'}\right\} \leq C_4 e^{-n^\eta}.$$

PROOF. From part (d) of Proposition 4.2, we have that $ny_1(\tau/n, \rho_c) \geq cn^{2\beta-1}$ for some $c > 0$, all $n$ and $\tau \in I_n$. Since $\rho \mapsto \vec{y}(\theta, \rho)$ is Lipschitz continuous [by Proposition 4.2, part (a)], there exists a finite constant $c'$ such that $\|\vec{y}(\theta, \rho) - \vec{y}(\theta, \rho_c)\| \leq c' n^{\beta'-1}$ for any $\theta \in [0, 1-\varepsilon]$ and $|\rho - \rho_c| \leq n^{\beta'-1}$.

By part (b) of Lemma 5.2, we thus get that for $\beta' < 2\beta - 1$, $\beta > 3/4$, some positive $C = C(\beta, \beta')$ and all $n$ large enough, if $\tau \in I_n$ and $|\rho - \rho_c| \leq n^{\beta'-1}$, then

$$\mathbb{E}z_1(\tau) \geq ny_1(\tau/n, \rho) - C_1\sqrt{n \log n}$$
$$\geq ny_1(\tau/n, \rho_c) - c'n^{\beta'} - C_1\sqrt{n \log n} \geq 2Cn^{2\beta-1}.$$

Applying now Lemma 5.2, part (a), we see that for any $\eta < (4\beta - 3)/2$, some $C' = C'(\beta, \beta', \eta)$ finite and all $n$ large enough

$$\widehat{\mathbb{P}}_{n,\rho}\{z_1(\tau) \leq n^{\beta'}\} \leq \widehat{\mathbb{P}}_{n,\rho}\{\|\vec{z}(\tau) - \mathbb{E}\vec{z}(\tau)\| \geq Cn^{2\beta-1}\} \leq C'e^{-n^{2\eta}},$$

whenever $\tau \in I_n$ and $|\rho - \rho_c| \leq n^{\beta'-1}$. To conclude, recall that there are at most $n$ integers $\tau \in I_n$. □

The second consequence of Lemma 5.2 is that with high probability also the process $\{\vec{z}(\tau)\}$ of distribution $\mathbb{P}_{n,\rho}(\cdot)$ belongs to the set $\mathcal{Q}(\eta)$ as long as $\tau/n$ is bounded away from 1.

COROLLARY 5.4. *For any $\varepsilon > 0$, there exists $\eta > 0$ and positive, finite constants $C_5, C_6$ such that if $\rho \in [\varepsilon, 1/\varepsilon]$, then*

$$(5.15) \quad \mathbb{P}_{n,\rho}\{(\vec{z}(\tau), \tau) \in \mathcal{Q}(\eta) \,\forall 0 \leq \tau \leq n(1-\varepsilon)\} \geq 1 - C_5 e^{-C_6 n}.$$

PROOF. From part (c) of Lemma 5.2 we have that for some $\eta' \in (0, \varepsilon)$, positive and finite $c_5$ and $c_6$,

$$(5.16) \quad \widehat{\mathbb{P}}_{n,\rho}\{(\vec{z}(\tau), \tau) \in \mathcal{Q}(\eta') \,\forall 0 \leq \tau \leq n(1-\varepsilon)\} \geq 1 - c_5 e^{-c_6 n}.$$

Applying the coupling of Lemma 5.1 with $\eta = \eta'/4$ for the value of $\varepsilon$ in the statement of this lemma, we also have that

$$(5.17) \quad \mathbb{P}\left\{\sup_{\tau \leq \tau_f} \|\vec{z}(\tau) - \vec{z}'(\tau)\| > \eta n\right\} \leq c_7 e^{-c_8 n},$$

where $\tau_f \leq n$ denotes the first time such that $(\vec{z}(\tau), \tau) \notin \mathcal{Q}(\eta)$. Further, if $\tau_f \leq n(1-\varepsilon)$ and $\sup_{\tau \leq \tau_f} \|\vec{z}(\tau) - \vec{z}'(\tau)\| \leq \eta n$, then necessarily $(\vec{z}'(\tau), \tau) \notin \mathcal{Q}(4\eta) = \mathcal{Q}(\eta')$ for $\tau = \tau_f \leq n(1-\varepsilon)$, an event whose probability is at most $c_5 e^{-c_6 n}$ [by (5.16)]. Combining the latter bound with (5.17) we find that

$$\mathbb{P}_{n,\rho}\{\tau_f \leq n(1-\varepsilon)\} \leq c_5 e^{-c_6 n} + c_7 e^{-c_8 n},$$

yielding (5.15) for $C_5 = c_5 + c_7$ and $C_6 = \min(c_6, c_8)$, both finite and positive. □

PROOF OF PROPOSITION 2.1. For $\{\vec{z}(\tau); \tau \geq 0\}$ distributed according to $\mathbb{P}_{n,\rho}(\cdot)$ let $\tau_{**}$ denote the first time at which $z_1(\tau) \leq 0$. Since by construction (and using Lemma 3.1), the sequence $\{z_1(\tau); 0 \leq \tau \leq \tau_{**}\}$ is distributed as the number of c-nodes of the graphs $G(\tau)$ having degree 1 under the decimation algorithm, we see that $P_l(n, \rho) = \mathbb{P}_{n,\rho}\{\tau_{**} \leq n - 1\}$.

Further, the core of the initial graph $G(0)$ includes at most $n - \tau_{**}$ vertices. Consequently, by Lemma 4.7 (cf. Remark 4.8), we can choose $D < \infty$ and $0 < \kappa < 1 - \theta_c$ such that

$$\mathbb{P}_{n,\rho}\{\tau_{**} \leq n(1-\kappa)\} \leq P_l(n, \rho) \leq \mathbb{P}_{n,\rho}\{\tau_{**} \leq n(1-\kappa)\} + \tfrac{1}{4}\delta_n,$$

for $|\rho - \rho_c| \leq n^{\beta'-1}$ and $\delta_n \equiv D n^{-1/2} (\log n)^2$.

By Corollary 5.4, there exist $0 < \eta \leq \varepsilon < \kappa$ and finite, positive $C_5, C_6$, such that $\{\vec{z}(\tau), 0 \leq \tau \leq n(1-\varepsilon)\} \subseteq \mathcal{Q}(\eta)$ with probability at least $1 - C_5 e^{-C_6 n}$, for all $n$. Hence, we have that

$$\mathbb{P}_{n,\rho}\left\{\min_{0 \leq \tau \leq \tau_*} z_1(\tau) \leq 0\right\} \leq P_l(n, \rho) \leq \mathbb{P}_{n,\rho}\left\{\min_{0 \leq \tau \leq \tau_*} z_1(\tau) \leq 0\right\} + C_5 e^{-C_6 n} + \tfrac{1}{4}\delta_n,$$

for

$$\tau_* = n(1-\kappa) \wedge \min\{\tau : (\vec{z}(\tau), \tau) \notin \mathcal{Q}(\eta)\}.$$

By Lemma 5.1, there exist $A > 0$, and a coupling of the process $\{\vec{z}(\tau)\}$ with a process $\{\vec{z}'(\tau)\}$ of distribution $\widehat{\mathbb{P}}_{n,\rho}(\cdot)$, such that, with probability larger than $1 - 1/2n$, up to time $\tau_*$ the distance between these two processes is at most $\varepsilon_n \equiv A \log n$. Therefore, enlarging $D$ if necessary, we have that

$$\mathbb{P}\left\{\min_{0 \leq \tau \leq \tau_*} z_1'(\tau) \leq -\varepsilon_n\right\} - \frac{1}{n} \leq P_l(n, \rho) \leq \mathbb{P}\left\{\min_{0 \leq \tau \leq \tau_*} z_1'(\tau) \leq \varepsilon_n\right\} + \frac{1}{n} + \frac{1}{4}\delta_n.$$

We have seen that $\tau_* < n(1-\kappa)$ with probability of at most $C_5 e^{-C_6 n}$. Hence, enlarging $D$ once more, we find that

$$\widehat{\mathbb{P}}_{n,\rho}\left\{\min_{0 \leq \tau \leq n(1-\kappa)} z_1(\tau) \leq -\varepsilon_n\right\} - \tfrac{1}{2}\delta_n \leq P_l(n, \rho)$$

$$\leq \widehat{\mathbb{P}}_{n,\rho}\left\{\min_{0 \leq \tau \leq n(1-\kappa)} z_1(\tau) \leq \varepsilon_n\right\} + \tfrac{1}{2}\delta_n.$$

With $\theta_c < (1 - \kappa)$, the set $[0, n(1 - \kappa)]$ is the disjoint union of $J_n$ as in the statement of the proposition and a set of the form $I_n$ of Corollary 5.3. Thus, bounding the probability of the event $\min_{\tau \in I_n} z_1(\tau) \leq \varepsilon_n$ via the latter corollary yields the thesis of the proposition [enlarging $D$ as needed for absorbing the term $C_4 \exp(-n^\eta)$ into $\delta_n$]. □

Let $\overline{\mathbb{P}}_{n,\rho}(\cdot)$ denote the law of the $\mathbb{R}^2$-valued Markov chain $\{\vec{z}'(\tau)\}$ of (2.12), where $\vec{z}'(0)$ has the uniform distribution $\mathbb{P}_{\mathcal{G}_l(n,m)}(\cdot)$ on the graph ensemble $\mathcal{G}_l(n,m)$ for $m \equiv \lfloor n\rho \rfloor$, and

$$\overline{\mathbb{P}}_{n,\rho}(\vec{z}'(\tau+1) = \vec{z}'(\tau) + \Delta_\tau + \widetilde{\mathbb{A}}_\tau(n^{-1}\vec{z}'(\tau) - \vec{y}(\tau/n))|\vec{z}'(\tau) = \vec{z}')$$
$$= \widehat{W}_{\tau/n}(\Delta_\tau|\vec{y}(\tau/n)).$$

We conclude this section by providing a coupling that keeps the process $\{\vec{z}'(\cdot)\}$ "sufficiently close" to $\{\vec{z}(\cdot)\}$ of distribution $\widehat{\mathbb{P}}_{n,\rho}(\cdot)$ throughout the time interval $J_n$ of interest to us.

PROPOSITION 5.5. *Fixing $\beta \in (3/4, 1)$ and $\beta' < 2\beta - 1$, for any $\delta > \beta - 1/2$ there exist finite constants $\alpha$, $c$ and a coupling of the processes $\{\vec{z}(\cdot)\}$ of distribution $\widehat{\mathbb{P}}_{n,\rho}(\cdot)$ and $\{\vec{z}'(\cdot)\}$ of distribution $\overline{\mathbb{P}}_{n,\rho}(\cdot)$ such that for all $n$ and $|\rho - \rho_c| \leq n^{\beta'-1}$,*

$$(5.18) \qquad \mathbb{P}\left\{\sup_{\tau \in J_n} \|\vec{z}(\tau) - \vec{z}'(\tau)\| \geq cn^\delta\right\} \leq \frac{\alpha}{4n}.$$

The key to Proposition 5.5 is the following elementary martingale concentration property.

LEMMA 5.6. *Consider an $\mathbb{R}^d$-valued discrete-time martingale $(Z_s, \mathcal{F}_s)$ with $Z_0 = 0$ and $U_s = Z_{s+1} - Z_s$ such that for some finite $\Gamma$ and a stopping time $\tau_*$ for $\mathcal{F}_s$*

(5.19)
$$\mathbb{E}[\|U_s\|^2 e^{\lambda \cdot U_s} \mid \mathcal{F}_s] \leq \Gamma \mathbb{E}[e^{\lambda \cdot U_s} \mid \mathcal{F}_s] < \infty$$

*whenever $s < \tau_*, \|\lambda\| < 1$.*

*Then, for any $0 \leq a < t\Gamma\sqrt{d}$,*
$$\mathbb{P}\{\|Z_{\min(t,\tau_*)}\| \geq a\} \leq 2d \exp\left\{-\frac{a^2}{2d\Gamma t}\right\}.$$

PROOF. Recall that for real-valued variable $V$, if $\mathbb{E}[V] = 0$ and $\mathbb{E}[V^2 \exp(uV)] \leq \kappa \mathbb{E}[\exp(uV)] < \infty$ for all $u \in [0,1]$, then $\mathbb{E}[\exp(V)] \leq \exp(\kappa/2)$ (bound the value of $\phi(1)$ for $\phi(u) \equiv \log \mathbb{E}[\exp(uV)]$ using $\phi(0) = \phi'(0) = 0$ and $\phi''(u) \leq \kappa$). In the special case of $d = 1$ and $\tau_* = \infty$, we have from (5.19) that the preceding assumptions hold for $\kappa = \Gamma\lambda^2$, $\|\lambda\| < 1$ and $V$ having the law of $\lambda U_s$ conditional on $\mathcal{F}_s$. Consequently, then $\mathbb{E}[\exp(\lambda U_s)|\mathcal{F}_s] \leq \exp(\Gamma\lambda^2/2)$, implying that $\mathbb{E}[M_t] \leq \mathbb{E}[M_0] = 1$ for the supermartingale $M_s = \exp(\lambda Z_s - \Gamma\lambda^2 s/2)$. Considering $a \in [0, t\Gamma)$ and $\lambda = a/(\Gamma t)$, we thus deduce that $\mathbb{P}\{Z_t \geq a\} \leq \exp\{-a^2/(2\Gamma t)\}$ in case $Z_s$ is a real-valued martingale for $\mathcal{F}_s$ and (5.19) holds for all $s < \infty$. The stated bound for $\mathbb{R}^d$-valued martingale $Z_t$ of coordinates $Z_{t,i}$ follows upon noting that the event $\{\|Z_t\| \geq a\}$ is contained in the union of the events $\{uZ_{t,i} \geq a/\sqrt{d}\}$ for $u = -1, 1$ and $i = 1, \ldots, d$, with $uZ_{s,i}$ real-valued martingales. Finally, we get the thesis in the general case, where $\mathbb{P}(\tau_* < \infty) > 0$, upon considering the (stopped) martingale $Z_{\min(s,\tau_*)}$. □

PROOF OF PROPOSITION 5.5. We couple the processes $\{\vec{z}'(\cdot)\} \stackrel{d}{=} \overline{\mathbb{P}}_{n,\rho}(\cdot)$, and $\{\vec{z}(\tau)\} \stackrel{d}{=} \widehat{\mathbb{P}}_{n,\rho}(\cdot)$ in a joint Markov process, by letting $\vec{z}'(0) = \vec{z}(0)$ and for $\tau = 0, 1, 2, \ldots, n-1$,
$$\mathbb{P}(\Delta\vec{z}(\tau) \neq \Delta_\tau | \mathcal{F}_\tau) = \|\widehat{W}_{\tau/n}(\cdot | n^{-1}\vec{z}(\tau)) - \widehat{W}_{\tau/n}(\cdot | \vec{y}(\tau/n))\|_{\mathrm{TV}},$$
where $\Delta\vec{z}(\tau) \equiv \vec{z}(\tau + 1) - \vec{z}(\tau)$ and $\mathcal{F}_\tau$ denotes the $\sigma$-algebra generated by $\{\vec{z}(\sigma), \vec{z}'(\sigma), \sigma \leq \tau\}$.

Fixing $\varepsilon < (1 - \theta_c)/2$, let $\tau_* \leq n$ denote the first value of $\tau$ such that $\|\vec{z}(\tau) - n\vec{y}(\tau/n)\| > K\sqrt{n \log n}$, with a finite $K = K(\varepsilon)$ such that by parts (a) and (b) of Lemma 5.2, for any $n$ and $\rho \in [\varepsilon, 1/\varepsilon]$,
$$\widehat{\mathbb{P}}_{n,\rho}\{\tau_* \leq n(1-\varepsilon)\} \leq n^{-1}.$$
Fix $\beta \in (3/4, 1)$, $\beta'$ and $\delta > \beta - 1/2$. With at most $n$ values for $\tau$ in $J_n$ we thus obtain (5.18) once we show that some $c < \infty$, all $n$ large enough

(5.20)
$$\sup_{\tau \in J_n} \mathbb{P}\{\tau < \tau_*, \|\vec{z}'(\tau) - \vec{z}(\tau)\| \geq cn^\delta\} \leq n^{-2}.$$

To this end, consider Doob's decomposition of the adapted process $N_s \equiv (\widetilde{\mathbb{B}}_0^{s-1})^{-1}(\vec{z}'(s) - \vec{z}(s))$ as the sum of an $\mathcal{F}_s$-martingale $\{Z_s\}$, null at zero, and the predictable sequence

$$V_{\tau+1} = \sum_{s=0}^{\tau} \Delta V_s \equiv \sum_{s=0}^{\tau} \mathbb{E}[N_{s+1} - N_s | \mathcal{F}_s].$$

It follows from our coupling that $\Delta V_s = (\widetilde{\mathbb{B}}_0^{s-1})^{-1} \vec{R}(n^{-1}\vec{z}(s), \vec{y}(s/n), s/n)$, where

$$\vec{R}(\vec{x}', \vec{x}, \theta) \equiv \vec{F}(\vec{x}, \theta) + \mathbb{I}_{\theta < \theta_n} \mathbb{A}(\vec{x}, \theta)[\vec{x}' - \vec{x}] - \vec{F}(\vec{x}', \theta)$$

(with $\theta_n \equiv \tau_n/n = \lfloor n\theta_c - n^\beta \rfloor/n$), and that for $\Delta_t^* \equiv \Delta_t - \Delta \vec{z}(t)$,

(5.21) $$U_t \equiv Z_{t+1} - Z_t = (\widetilde{\mathbb{B}}_0^{t-1})^{-1}\{\Delta_t^* - \mathbb{E}[\Delta_t^* | \mathcal{F}_t]\}.$$

Since $A_{ab}(\vec{x}, \theta) = \partial_{x_b} F_a(\vec{x}, \theta)$ with $\vec{F}(\vec{x}, \theta)$ having Lipschitz continuous derivatives on $\widehat{q}_+(\varepsilon)$, it follows that $\|\vec{R}(\vec{x}', \vec{x}, \theta)\| \leq c_0 \|\vec{x}' - \vec{x}\|^2$ for some $c_0 = c_0(\varepsilon)$ finite, provided $\theta < \theta_n$ and both $(\vec{x}, \theta)$ and $(\vec{x}', \theta)$ are in $\widehat{q}_+(\varepsilon)$. By the Lipschitz continuity of $\vec{F}(\vec{x}, \theta)$ we also have that $\|\vec{R}(\vec{x}', \vec{x}, \theta)\| \leq c_0 \|\vec{x}' - \vec{x}\|$ in case $\theta \geq \theta_n$, as soon as $(\vec{x}, \theta)$ and $(\vec{x}', \theta)$ are in $\widehat{q}(\varepsilon)$.

Recall Lemma 4.3 that for some finite $n_0$ and $\kappa$ we have that $\|(\widetilde{\mathbb{B}}_0^{\tau-1})^{-1}\| \leq \kappa$ for all $\tau$, $\rho \in [\varepsilon, 1/\varepsilon]$ and $n \geq n_0$. In the course of proving part (a) of Proposition 4.2 we have seen that the distance of $(\vec{y}(\theta, \rho), \theta)$ from the complement of $\widehat{q}(\varepsilon)$ is bounded away from zero, uniformly in $\theta \leq 1 - 2\varepsilon$ and $\rho \in [\varepsilon, 1/\varepsilon]$. Further, $y_1(\theta, \rho) \geq \kappa' n^{2(\beta-1)}$ for some $\kappa' > 0$, all $n$, $|\rho - \rho_c| \leq n^{\beta'-1}$, and $\theta \leq \theta_n$ (cf. proof of Corollary 5.3). Consequently, for some finite $n_1 = n_1(K, \varepsilon)$ and all $n \geq n_1$, the event $\{s < \tau_*\}$ implies that both $(n^{-1}\vec{z}(s), s/n)$ and $(\vec{y}(s/n), s/n)$ are in $\widehat{q}(\varepsilon)$ when $s \leq n\theta_c + n^\beta$, and in case $s < \tau_n$ they are also in $\widehat{q}_+(\varepsilon)$. We deduce that if $n \geq n_1$ and $\{s < \tau_*\}$, then

$$\|\Delta V_s\| \leq c_0 \|(\widetilde{\mathbb{B}}_0^{s-1})^{-1}\| \|n^{-1}\vec{z}(s) - \vec{y}(s/n)\|^2 \leq c_0 \kappa K^2 n^{-1} \log n,$$

when $s < \tau_n$, whereas $\|\Delta V_s\| \leq c_0 \kappa K n^{-1/2}(\log n)^{1/2}$ for $s \in J_n$. Hence, for some finite $c_1$ and all $n \geq n_1$, the event $\{\tau < \tau_*\}$ implies for $\tau \in J_n$ that $\|V_\tau\| \leq c_1 n^{\beta - 1/2}(\log n)^{1/2}$. Fixing $\eta \in (1/4, \beta - 1/2)$, since $\|\vec{z}'(\tau) - \vec{z}(\tau)\| \leq \|\widetilde{\mathbb{B}}_0^{\tau-1}\|[\|V_\tau\| + \|Z_\tau\|]$ and $\|\widetilde{\mathbb{B}}_0^{\tau-1}\|$ are bounded uniformly in $n$, $\tau$ and $\rho$, we thus get (5.20) by considering Lemma 5.6 at $\tau \in J_n$ and $a = n^\eta$, provided we show that for some $c_2$ finite, the martingale differences $U_t$ of (5.21) satisfy the inequality (5.19) with $\Gamma = c_2 n^{-1/2}(\log n)^{1/2}$ (as indeed $n^\eta \leq \tau_n \Gamma \sqrt{d}$ for all $n$ large enough and $n^{2\eta}/2d\Gamma n \to \infty$). To this end, note first that by the total variation bound of Lemma 4.6 and the definition of $\tau_*$, for $t < \tau_*$ our coupling of $(\vec{z}, \vec{z}')$ results with

$$\mathbb{P}(\Delta_t^* \neq 0 | \mathcal{F}_t) \leq L \|n^{-1}\vec{z}(t) - \vec{y}(t/n)\| \leq LKn^{-1/2}(\log n)^{1/2} \equiv u_n.$$

Further, the bounded support of $\widehat{W}_\theta(\cdot|\vec{x})$ implies that $\|\Delta_t^*\| \leq 4l$, so for $t < \tau_*$ also

$$\|\mathbb{E}[\Delta_t^*|\mathcal{F}_t]\| \leq 4l\mathbb{P}(\Delta_t^* \neq 0|\mathcal{F}_t) \leq 4l\min(u_n, 1).$$

From the preceding estimates we deduce that $U_t$ of (5.21) is such that $\|U_t\| \leq 8l\kappa$ and when $t < \tau_*$, also

$$\mathbb{P}(\|U_t\| > 4l\kappa \min(u_n, 1)|\mathcal{F}_t) \leq \mathbb{P}(\Delta_t^* \neq 0|\mathcal{F}_t) \leq u_n.$$

These two facts easily imply that if $\|\lambda\| \leq 1$ and $t < \tau_*$, then the inequality (5.19) holds for $\Gamma = 2(8l\kappa)^2 e^{16l\kappa} u_n$ which as we have already seen, completes the proof of the proposition. $\square$

**6. Gaussian approximation and proof of Proposition 2.2.** This section is devoted to the proof of Proposition 2.2. Specifically, building on Proposition 5.5, in Section 6.1 we approximate the Markov process $\{\vec{z}(\tau)\}$ of distribution $\widehat{\mathbb{P}}_{n,\rho}(\cdot)$ by a Brownian motion with a quadratic shift, when $\tau$ is within the window $J_n$ around the critical time. Then, in Section 6.2 we show how the one-sided Brownian motion with quadratic shift can be replaced by a two-sided motion, once its initial condition is appropriately mapped to the distribution of the two-sided motion at the critical time. Finally, in Section 6.3 we show that this distribution [which a priori depends on the law of $\vec{z}(0)$] is also well approximated by a Gaussian law and complete the proof of Proposition 2.2.

6.1. *Local approximation by a Brownian motion with quadratic shift.* Our goal here is to approximate the probabilities of interest to us in terms of the minimal value of the Brownian motion with quadratic shift

$$(6.1) \quad X_n(\tau) \equiv n^{1/3}[X(n^{-2/3}(\tau - 0.5 - n\theta_c)) - X(n^{-2/3}(\tau_n - 0.5 - n\theta_c))],$$

within $J_n \equiv [n\theta_c - n^\beta, n\theta_c + n^\beta]$, for the process $\{X(t)\}$ of Proposition 2.2. As stated in the following lemma, while doing this we also approximate the law of $z_1(\tau_n)$ by that of the sum of $\vec{u}_n \cdot [\vec{z}(0) - n\vec{y}(0)]$, where $\vec{u}_n$ denotes the first row of $\widetilde{\mathbb{B}}_0^{\tau_n - 1}$, and an independent normal random variable of mean $ny_1^*(\tau_n)$ and variance $n[(\mathbb{Q}_{\tau_n})_{11} - \vec{u}_n^\dagger \mathbb{Q}(0, \rho)\vec{u}_n]$.

LEMMA 6.1. *Fixing $\beta \in (3/4, 1)$ and $A > 0$, set $\beta' < 2\beta - 1$, $\varepsilon_n = A\log n$ as in Proposition 2.1 and $Y_n \equiv \inf_{t \in J_n} X_n(t)$ for $X_n(\cdot)$ of (6.1). Then, for any $\delta > 3\beta - 2$, there exist positive, finite constants $\alpha$ and $C$ such that for any $n$ and $|\rho - \rho_c| \leq n^{\beta' - 1}$,*

$$\mathbb{P}\{\xi_n^* + \xi_n + Y_n \leq -Cn^\delta\} - \frac{\alpha}{n} \leq \widehat{\mathbb{P}}_{n,\rho}\left\{\min_{\tau \in J_n} z_1(\tau) \leq \pm\varepsilon_n\right\}$$
$$\leq \mathbb{P}\{\xi_n^* + \xi_n + Y_n \leq Cn^\delta\} + \frac{\alpha}{n},$$

where $\xi_n \equiv \vec{u}_n \cdot [\vec{z}(0) - n\vec{y}(0)]$, *the normal random variable* $\xi_n^*$ *of mean* $ny_1^*(\tau_n)$ *and variance* $n[(\mathbb{Q}_{\tau_n})_{11} - \vec{u}_n^\dagger \mathbb{Q}(0, \rho) \vec{u}_n]$, *and* $X_n(\cdot)$, *are mutually independent*.

PROOF. The strategy we follow is to progressively simplify the process $\{\vec{z}(\tau)\}$ of distribution $\widehat{\mathbb{P}}_{n,\rho}(\cdot)$ until obtain the stated bounds of the lemma, where each simplification is justified by a coupling argument. The first and most important step of this program has already been done in Proposition 5.5. Since the chain $\{\vec{z}'(\cdot)\}$ of law $\overline{\mathbb{P}}_{n,\rho}(\cdot)$ has independent increments for $\tau \geq \tau_n$, we can apply Sakhanenko's refinement of the Hungarian construction, to [22, 32] for the uniformly bounded (by $4l$) independent increments $\xi_i = z_1'(\tau_n + i) - z_1'(\tau_n + i - 1)$. We then deduce the existence of a real-valued Gaussian process $b_n(\tau)$, independent of $z_1'(\tau_n)$, such that $b_n(\tau_n) = 0$; its independent increments $\Delta b_n(\tau) \equiv b_n(\tau + 1) - b_n(\tau)$ have mean and variance

$$\mathbb{E}\Delta b_n(\tau) = F_1(\vec{y}(\tau/n, \rho), \tau/n), \qquad \operatorname{Var} \Delta b_n(\tau) = G_{11}(\vec{y}(\tau/n, \rho), \tau/n)$$

[matching the corresponding moments of $z_1'(\tau + 1) - z_1'(\tau)$], such that for some finite $c_0$, $\alpha$ and all $n$, $\rho$,

$$(6.2) \qquad \mathbb{P}\left\{ \sup_{\tau \in J_n} |z_1'(\tau) - z_1'(\tau_n) - b_n(\tau)| \geq c_0 \log n \right\} \leq \frac{\alpha}{4n}$$

(the latter follows by Chebyshev's inequality from [32]; see, e.g., [34], Theorem A).

Considering the representation (2.14) for $\vec{z}'(\tau_n)$ we see that $z_1'(\tau_n) - \xi_n$ is the sum of the uniformly bounded real-valued independent variables $(\widetilde{\mathbb{B}}_{\sigma+1}^{\tau_n-1} \Delta_\sigma)_1$, $\sigma = 0, \ldots, \tau_n - 1$ plus a nonrandom constant. Hence, similarly to the derivation of (6.2) we obtain that for some finite $c_0$, $\alpha$ and all $n$, $\rho$,

$$(6.3) \qquad \mathbb{P}\{|z_1'(\tau_n) - \xi_n - \xi_n^*| \geq c_0 \log n\} \leq \frac{\alpha}{4n},$$

where $\xi_n^*$ is a normal random variable, independent of $\xi_n$ and $b_n(\cdot)$, whose mean and variance match those of $z_1'(\tau_n) - \xi_n$. It is not hard to verify that the latter mean and variance are indeed $ny_1^*(\tau_n)$ and $n[(\mathbb{Q}_{\tau_n})_{11} - \vec{u}_n^\dagger \mathbb{Q}(0, \rho) \vec{u}_n]$, as stated.

We clearly have the representation

$$b_n(\tau) = \sum_{\sigma=\tau_n}^{\tau-1} \mathbb{E}\Delta b_n(\sigma) + B\left( \sum_{\sigma=\tau_n}^{\tau-1} \operatorname{Var} \Delta b_n(\sigma) \right),$$

for a standard Brownian motion $B(\cdot)$. Further, the real-valued Gaussian process $\{X_n(t), t \geq \tau_n\}$ of (6.1) admits the representation

$$X_n(t) = \widetilde{F} \int_{\tau_n - 0.5}^{t - 0.5} (\sigma/n - \theta_c) \, d\sigma + B(\widetilde{G}(t - \tau_n)),$$

for the same standard Brownian motion $B(\cdot)$, where $\widetilde{G} = G_{11}$ and $\widetilde{F} = \frac{dF_1}{d\theta}$, both evaluated at $\theta = \theta_c$ and $\vec{y} = \vec{y}(\theta_c, \rho_c)$ [so $\widetilde{F}$ is as defined in (2.7)]. Combining (5.18), (6.2) and (6.3) we establish the thesis of the lemma upon showing that the preceding coupling of $b_n(\cdot)$ and $X_n(\cdot)$ is such that for some $\alpha$, $c_1$ finite and all $n$,

$$(6.4) \qquad \mathbb{P}\left\{\sup_{t \in J_n} |b_n([t]) - X_n(t)| \geq 3c_1 n^\delta\right\} \leq \frac{\alpha}{4n}.$$

The sup in (6.4) is taken over all *real-valued* $t \in J_n = [n\theta_c - n^\beta, n\theta_c + n^\beta]$, while in the sequel we use $\tau \in J_n$ to denote an *integer* in the same interval.

With $\{X_n(\tau + t) - X_n(\tau) : t \in [0,1]\}$ having the same law as $\{B(\widetilde{G}t) + a_{n,\tau}(t) : t \in [0,1]\}$ for nonrandom $a_{n,\tau}(t)$ which are bounded uniformly in $t \in [0,1]$, $n$ and $\tau \in J_n$, we obviously get (6.4) upon showing that

$$(6.5) \qquad \sup_{\tau \in J_n} \mathbb{P}\{|b_n(\tau) - X_n(\tau)| \geq 2c_1 n^\delta\} \leq n^{-2}.$$

Inequality (6.5) is a direct consequence of having a finite $\kappa$ such that, for $\Delta X_n(\tau) = X_n(\tau + 1) - X_n(\tau)$,

$$e_1(\tau) \equiv |\operatorname{Var}\Delta b_n(\tau) - \operatorname{Var}\Delta X_n(\tau)| = |G_{11}(\vec{y}(\tau/n, \rho), \tau/n) - \widetilde{G}| \leq \kappa n^{\beta-1},$$

$$e_2(\tau) \equiv |\mathbb{E}\Delta b_n(\tau) - \mathbb{E}\Delta X_n(\tau)|$$
$$= |F_1(\vec{y}(\tau/n, \rho), \tau/n) - (\tau/n - \theta_c)\widetilde{F}| \leq \kappa n^{2(\beta-1)},$$

for all $\tau \in J_n$ and $|\rho - \rho_c| \leq n^{\beta'-1}$. Indeed, since $\delta > \beta + 2(\beta - 1) > 1 - \beta$ and the interval $J_n$ is of length $2\lceil n^\beta \rceil$, taking $c_1$ large enough so $c_1 n^\delta \geq \kappa n^{2(\beta-1)}|J_n|$ for all $n$, the stated bound on $e_2(\cdot)$ guarantees that $|\mathbb{E}b_n(\tau) - \mathbb{E}X_n(\tau)| \leq c_1 n^\delta$ for all $\tau \in J_n$, whereas the corresponding bound on $e_1(\cdot)$ guarantees that $\operatorname{Var}(b_n(\tau) - X_n(\tau)) \leq c_1 n^\delta n^{1-\beta}$, leading (by standard Gaussian tail estimates) to (6.5).

Turning to bound $e_1(\tau)$ and $e_2(\tau)$, recall that $\tau \in J_n$ and $|\rho - \rho_c| \leq n^{\beta'-1}$ imply that $(\tau/n, \rho) \in [0, 1-\varepsilon] \times [\varepsilon, 1/\varepsilon]$, so $|\vec{y}(\tau/n, \rho) - \vec{y}(\tau/n, \rho_c)| \leq C_1|\rho - \rho_c|$ for some constant $C_1 = C_1(\varepsilon)$ by the Lipschitz continuity of $\rho \mapsto \vec{y}(\theta, \rho)$ [see part (a) of Proposition 4.2]. Further, then $(\vec{y}(\tau/n, \rho), \tau/n) \in \widehat{q}(\varepsilon)$, so by Lemma 4.1 we have the Lipschitz continuity of $\rho \mapsto F_1(\vec{y}(\tau/n, \rho), \tau/n)$ and $\rho \mapsto G_{11}(\vec{y}(\tau/n, \rho), \tau/n)$. That is, for some constant $C_2 = C_2(\varepsilon)$ and all such $\tau$, $\rho$, $n$,

(6.6) $\quad e_1(\tau) \leq |G_{11}(\vec{y}(\tau/n, \rho_c), \tau/n) - G_{11}(\vec{y}(\theta_c, \rho_c), \theta_c)| + C_2|\rho - \rho_c|,$

(6.7) $\quad e_2(\tau) \leq |\widehat{F}_1(\tau/n) - (\tau/n - \theta_c)\widetilde{F}| + C_2|\rho - \rho_c|,$

where $\widehat{F}_1(\theta) \equiv F_1(\vec{y}(\theta, \rho_c), \theta)$. Similarly, the Lipschitz continuity of $\theta \mapsto \vec{y}(\theta, \rho_c)$ on $[0, 1-\varepsilon]$ [from part (a) of Proposition 4.2] together with that of $(\vec{x}, \theta) \mapsto G_{11}(\vec{x}, \theta)$ on $\widehat{q}(\varepsilon)$ (by Lemma 4.1) result in

$$|G_{11}(\vec{y}(\tau/n, \rho_c), \tau/n) - G_{11}(\vec{y}(\theta_c, \rho_c), \theta_c)| \leq C_3|\tau/n - \theta_c| \leq C_3 n^{\beta-1},$$

for some $C_3 = C_3(\varepsilon)$ and all $\tau \in J_n$. Thus, with $\beta' - 1 < 2(\beta - 1)$, we get from (6.6) that $e_1(\tau) \leq \kappa n^{\beta-1}$ for all $\tau \in J_n$ and $|\rho - \rho_c| \leq n^{\beta'-1}$, as stated. As for bounding $e_2(\tau)$, recall that $\theta \mapsto \vec{y}(\theta, \rho_c)$ is infinitely continuously differentiable on $[0, 1-\varepsilon]$ [cf. parts (b) and (d) of Proposition 4.2]. Further, as $(\vec{y}(\theta, \rho_c), \theta) \in \widehat{q}_+(\varepsilon)$ for all $\theta \in [0, 1-\varepsilon]$ [by (a) and (d) of Proposition 4.2], from Lemma 4.1 we have that $\widehat{F}_1(\cdot)$ is differentiable on $[0, 1-\varepsilon]$ with a Lipschitz continuous derivative. Recall that $dy_1/d\theta = 0$ at $\theta = \theta_c$ and $\rho = \rho_c$ [see part (d) of Proposition 4.2]. Hence, $\widehat{F}_1(\theta_c) = 0$ [in view of the ODE (2.5)], and with $\widetilde{F} = \widehat{F}_1'(\theta_c)$ we deduce that for some $C_4 = C_4(\varepsilon)$ and all $\theta \in [0, 1-\varepsilon]$,

$$|\widehat{F}_1(\theta) - (\theta - \theta_c)\widetilde{F}| \leq C_4|\theta - \theta_c|^2.$$

Combining (6.7) with the latter bound (for $\theta = \tau/n$ and $\tau \in J_n$, so $|\theta - \theta_c| \leq n^{\beta-1}$), we conclude that $e_2(\tau) \leq \kappa n^{2(\beta-1)}$ for all $\tau \in J_n$ and $|\rho - \rho_c| \leq n^{\beta'-1}$, as stated. □

6.2. *Brownian computations.* We show in the sequel that for large $s$ and $u$ the distribution of

$$V_{s,u} = \inf_{t \in [-s, u]} X(t) - X(-s)$$

is well approximated by that of $V_* - \widetilde{X}(-s)$ for $V_* \equiv \inf_{t \in \mathbb{R}} X(t)$ and $\{\widetilde{X}(t)\}$ an independent copy of $\{X(t)\}$. More precisely, we prove:

LEMMA 6.2. *With the preceding definitions, for $0 < \varphi < 4(1-\psi)/3 - 1$, all $s, u$ large enough and any nonrandom $v$,*

$$\mathbb{P}\{V_* - \widetilde{X}(-s) \geq v + 2s^{-\psi}\} - 5e^{-s^\varphi} \leq \mathbb{P}\{V_{s,u} \geq v\}$$
$$\leq \mathbb{P}\{V_* - \widetilde{X}(-s) \geq v - 2s^{-\psi}\} + 5e^{-(s \wedge u)^\varphi}.$$

PROOF. Conditioning upon the value of $X(-s)$ we have on account of the independence of $\{X(-t): t \geq 0\}$ and $\{X(t): t \geq 0\}$ that

$$\mathbb{P}\{V_{s,u} \geq v\} = \mathbb{E}[p_s(v + \widetilde{X}(-s), \widetilde{X}(-s))q_{-u,0}(v + \widetilde{X}(-s))],$$

where for $s > \theta \geq 0$ and any $a, b$,

$$p_s(a, b) \equiv \mathbb{P}\left\{\inf_{-s \leq t \leq 0} X(t) \geq a \Big| X(-s) = b\right\},$$

$$q_{-s,-\theta}(a) \equiv \mathbb{P}\left\{\inf_{-s \leq t \leq -\theta} X(t) \geq a\right\}.$$

Recall that the law of $\{X(t): -s \leq t \leq 0\}$ conditional upon $\{X(-s) = b\}$ is merely the law of $\{X_{b,s}(t) \equiv X(t) - \frac{t}{s}(b - X(-s)): -s \leq t \leq 0\}$. Thus, in particular,

(6.8) $\quad \mathbb{P}\{V_{s,u} \geq v\} = \mathbb{E}[p^{(s)}_{-s,0}(v + \widetilde{X}(-s), \widetilde{X}(-s)) q_{-u,0}(v + \widetilde{X}(-s))],$

where

$$p^{(c)}_{-s,-\theta}(a,b) \equiv \mathbb{P}\Big\{\inf_{-s \leq t \leq -\theta} X_{b,c}(t) \geq a\Big\}.$$

Similarly,

(6.9) $\quad \mathbb{P}\{V_* - \widetilde{X}(-s) \geq v\} = \mathbb{E}[q_{-\infty,0}(v + \widetilde{X}(-s))^2].$

Fixing $0 < \varphi < 4(1-\psi)/3 - 1$, choose $(\varphi+1)/2 < \kappa < 2(1-\psi)/3$. Then, setting $\rho = 1 - \psi - \kappa$ and $\theta = s^\rho$ (so $s^{-\psi} = \frac{\theta}{s} s^\kappa$), it follows that if $|b - \frac{1}{2}\widetilde{F}s^2| \leq s^\kappa$, then for all $s$ large enough

(6.10)
$$\mathbb{P}\Big\{\sup_{-\theta \leq t \leq 0} |X_{b,s}(t) - X(t)| \geq 2s^{-\psi}\Big\}$$
$$\leq \mathbb{P}\{|X(-s) - \tfrac{1}{2}\widetilde{F}s^2| \geq s^\kappa\} \leq e^{-s^\varphi}.$$

Consequently, for any value of $a$,

(6.11) $\quad q_{-\theta,0}(a + 2s^{-\psi}) - e^{-s^\varphi} \leq p^{(s)}_{-\theta,0}(a,b) \leq q_{-\theta,0}(a - 2s^{-\psi}) + e^{-s^\varphi}.$

Relying upon these bounds we next show that if $|b - \frac{1}{2}\widetilde{F}s^2| \leq s^\kappa$, then for $s$ large enough and all $a$,

(6.12) $\quad q_{-\infty,0}(a + 2s^{-\psi}) - 3e^{-s^\varphi} \leq p^{(s)}_{-s,0}(a,b) \leq q_{-\infty,0}(a - 2s^{-\psi}) + 3e^{-s^\varphi}.$

Indeed, with $X_{b,s}(0) = X(0) = 0$, clearly (6.12) holds for $a > 0$ [as then $q_{-\infty,0}(a) = p^{(s)}_{-s,0}(a,b) = 0$]. Next recall that for $c \geq s \geq \theta \geq 0$ and any $a, b$,

(6.13)
$$p^{(c)}_{-\theta,0}(a,b) - [1 - p^{(c)}_{-s,-\theta}(a,b)]$$
$$\leq p^{(c)}_{-s,0}(a,b) \leq p^{(c)}_{-\theta,0}(a,b),$$

(6.14) $\quad q_{-\infty,0}(a) \leq q_{-\theta,0}(a) \leq q_{-\infty,0}(a) + [1 - q_{-\infty,-\theta}(a)].$

Combining these with the monotonicity in $a$ of the functions $p^{(c)}$ and $q$, we thus get (6.12) also for $a \leq 0$ out of (6.11) as soon as we show that for $\theta = s^\rho$ and all large $s$

(6.15) $\quad 1 - q_{-\infty,-\theta}(2s^\kappa) = \mathbb{P}\Big\{\inf_{t \leq -\theta} X(t) < 2s^\kappa\Big\} \leq e^{-s^\varphi},$

(6.16) $\quad 1 - p^{(s)}_{-s,-\theta}(0,b) = \mathbb{P}\Big\{\inf_{-s \leq t \leq -\theta} X_{b,s}(t) < 0\Big\} \leq 2e^{-s^\varphi}.$

Now, since $2\rho > \kappa > \varphi$, it follows by standard Gaussian tail estimates that for $\theta = s^\rho$ and $s$ large

$$\mathbb{P}\Big\{\inf_{t \le -\theta} X(t) < 2s^\kappa\Big\}$$
$$\le \sum_{\tau=\lfloor\theta\rfloor}^\infty \Big[\mathbb{P}\Big\{\sqrt{\widetilde{G}}W(\tau) \le -\frac{\widetilde{F}}{6}\tau^2\Big\}$$
$$+ \mathbb{P}\Big\{\sqrt{\widetilde{G}}\inf_{0 \le t \le 1}[W(\tau+t) - W(\tau)] \le -\frac{\widetilde{F}}{6}\tau^2\Big\}\Big]$$
$$\le 3\sum_{\tau=\lfloor\theta\rfloor}^\infty e^{-\widetilde{F}^2\tau^2/(72\widetilde{G})} \le e^{-s^\varphi},$$

thus establishing (6.15). Further, as $|X_{b,s}(t) - X(t)| \le |X(-s) - \frac{\widetilde{F}}{2}s^2| + |b - \frac{\widetilde{F}}{2}s^2|$, we deduce from (6.10) that if $|b - \frac{1}{2}\widetilde{F}s^2| \le s^\kappa$, then

$$\mathbb{P}\Big\{\sup_{-s \le t \le 0}|X_{b,s}(t) - X(t)| \ge 2s^\kappa\Big\} \le e^{-s^\varphi},$$

which together with (6.15) implies the bound (6.16).

We now apply in (6.8) standard Gaussian tail estimates for $|\widetilde{X}(-s) - \frac{1}{2}\widetilde{F}s^2| > s^\kappa$, and the bounds of (6.12) otherwise. With the $[0,1]$-valued $p_{-s,0}^{(c)}(a,b)$ and $q_{-s,0}(a)$ monotone in $s$ and $a$, this results in

$$\mathbb{E}[q_{-\infty,0}(v + \widetilde{X}(-s) + 2s^{-\psi})^2] - 4e^{-s^\varphi}$$
(6.17)
$$\le \mathbb{P}\{V_{s,u} \ge v\}$$
$$\le \mathbb{E}[q_{-\infty,0}(v + \widetilde{X}(-s) - 2s^{-\psi})q_{-u,0}(v + \widetilde{X}(-s) - 2s^{-\psi})] + 4e^{-s^\varphi}.$$

Finally, if $a > 0$, then $q_{-u,0}(a) = q_{-\infty,0}(a) = 0$, whereas for $a \le 0$, taking $\theta = u$ in (6.14) we find by (6.15) and the monotonicity of $q_{-\infty,-u}(a)$ that $q_{-u,0}(a) \le q_{-\infty,0}(a) + \exp(-u^\varphi)$ for $u$ large. Combining this upper bound on $q_{-u,0}$ with (6.9) and (6.17) provides the thesis of the lemma. $\square$

Since $Y_n$ of Lemma 6.1 has the same law as $n^{1/3}V_{s,u}$ for $s = s(n) = n^{-2/3}[n\theta_c + 0.5 - \tau_n]$ and $u = u(n) = n^{-2/3}[n^\beta - 0.5]$, we have the following immediate corollary of Lemmas 6.1 and 6.2.

COROLLARY 6.3. *Fixing $\beta \in (3/4, 1)$ and $A > 0$, set $J_n = [n\theta_c - n^\beta, n\theta_c + n^\beta]$ and $\beta' < 2\beta - 1$, $\varepsilon_n = A\log n$ as in Proposition 2.1. Let $\{\widetilde{X}(t)\}$ denote an i.i.d. copy of the process $\{X(t)\}$ of Proposition 2.2. Then, for any $\nu <$* $\min(7/3-3\beta, \beta/4-1/6)$, there exist $c$ finite such that for all $n$ and $|\rho-\rho_c| \leq n^{\beta'-1}$,

$$\mathbb{P}\left\{\widetilde{\xi}_n + \inf_t X(t) \leq -n^{-\nu}\right\} - \frac{c}{n} \leq \widehat{\mathbb{P}}_{n,\rho}\left\{\min_{\tau \in J_n} z_1(\tau) \leq \pm\varepsilon_n\right\}$$
$$\leq \mathbb{P}\left\{\widetilde{\xi}_n + \inf_t X(t) \leq n^{-\nu}\right\} + \frac{c}{n},$$

where $\widetilde{\xi}_n \equiv n^{-1/3}(\xi_n + \xi_n^*) - \widetilde{X}(-n^{\beta-2/3})$ [and $\xi_n$ and $\xi_n^*$ of Lemma 6.1 are independent of both $\{X(\cdot)\}$ and $\{\widetilde{X}(\cdot)\}$].

PROOF. Fixing $\nu < \min(7/3-3\beta, \beta/4-1/6)$, set $\delta > 3\beta-2$ of Lemma 6.1 and $\psi \in (0, 1/4)$ of Lemma 6.2 such that $\nu < 1/3 - \delta$ and $\nu < \psi(\beta - 2/3)$. Conditioning on the values of $\xi_n$ and $\xi_n^*$ we apply Lemma 6.2 for the values of $s = s(n)$ and $u = u(n)$ indicated above, taking there $v(n) = n^{-1/3}[\pm C n^\delta - \xi_n - \xi_n^*]$ for the finite constant $C$ of Lemma 6.1. With $s(n) \wedge u(n) \geq n^{\beta-2/3} - 2$ and $\beta > 2/3$, the error terms $5\exp(-(s(n) \wedge u(n))^\psi)$ are accommodated within $c/(2n)$ for some finite $c$ and all $n$. Further, enlarging $c$ if needed, with $|s(n) - n^{\beta-2/3}| \leq 2n^{-2/3}$ and $\nu < \frac{2}{3}(\beta - 2/3) < 1/3$, it is easy to see that for all $n$,

$$\mathbb{P}\left\{|\widetilde{X}(-s(n)) - \widetilde{X}(-n^{\beta-2/3})| \geq \frac{1}{2}n^{-\nu}\right\} \leq \frac{c}{2n}.$$

Our choice of $\psi$ and $\delta$ is such that $Cn^{\delta-1/3} + 2s(n)^{-\psi} \leq \frac{1}{2}n^{-\nu}$ for all $n \geq n_0$, so adding to $c$ the constant $\alpha$ of Lemma 6.1 and making sure that $c \geq 2n_0$, upon taking the expectation over $\xi_n$ and $\xi_n^*$ our thesis follows from the latter lemma. □

6.3. *Proof of Proposition 2.2.* Fixing $\beta \in (3/4, 1)$, $r \in \mathbb{R}$ and $\beta' < 2\beta - 1$, we have that $|\rho_n - \rho_c| \leq n^{\beta'-1}$ for $\rho_n = \rho_c + rn^{-1/2}$ and all $n$ large enough. Further, taking $\beta = 10/13 \in (3/4, 1)$ which maximizes the bound $\nu_0 \equiv \min(5/2-3\beta, \beta/4) - 1/6$ on $\nu$ in Corollary 6.3 leads to $\nu_0 = 5/26 - 1/6 > 0$. Thus, fixing $A > 0$, the statement (2.11) of the proposition is a consequence of Corollary 6.3, once we show that for any $1/6 < \eta < \nu + 1/6 < 5/26$ and all $n$ large enough,

$$\left|\mathbb{P}\left\{\widetilde{\xi}_n + \inf_t X(t) \leq \pm n^{-\nu}\right\} - \mathbb{P}\left\{n^{1/6}\xi(r) + \inf_t X(t) \leq 0\right\}\right| \leq n^{-\eta},$$

where $\xi(r)$ denotes a normal random variable of mean $(\frac{\partial y_1}{\partial \rho})r$ and variance $Q_{11}$ (both evaluated at $\theta = \theta_c$ and $\rho = \rho_c$), independent of $X(\cdot)$. Conditioning on $\widetilde{\xi}_n$ and $\xi(r)$, by the independence of $\{X(t) : t \geq 0\}$ and $\{X(t) : t \leq 0\}$, this is equivalent to

(6.18) $\qquad |\mathbb{E}[q(\pm n^{-\nu} - \widetilde{\xi}_n)^2] - \mathbb{E}[q(-n^{1/6}\xi(r))^2]| \leq n^{-\eta},$

where $q(a) \equiv \mathbb{P}(\inf_{t \leq 0} X(t) \geq a)$. With $q^2(a)$ a $[0,1]$-valued monotone non-increasing function that approaches zero as $a \to \infty$, we have that for any random variables $Y, Z$ and nonrandom $v$,

$$|\mathbb{E}[q(vZ)^2] - \mathbb{E}[q(vY)^2]| \leq \sup_x |\mathbb{P}(Z \leq x) - \mathbb{P}(Y \leq x)|.$$

Applying this for $v = -n^{1/6}$, $Z = n^{-1/6}\widetilde{\xi}_n \pm n^{-(\nu+1/6)}$ and $Y = \xi(r)$ of bounded density, we deduce that (6.18) holds for all $n$ large enough, thus completing the proof of Proposition 2.2 as soon as we show that for $\eta < 5/2 - 3\beta$ and all $n$ large enough,

$$(6.19) \quad \sup_{x \in \mathbb{R}} |\mathbb{P}(n^{-1/6}\widetilde{\xi}_n \leq x) - \mathbb{P}(\xi(r) \leq x)| \leq n^{-\eta}.$$

To this end, recall that $n^{-1/6}\widetilde{\xi}_n = n^{-1/2}\xi_n + n^{-1/2}\xi_n^* - n^{-1/6}\widetilde{X}(-n^{\beta-2/3})$, where the latter three random variables are independent of each other. Hence, in view of (4.20) we have that

$$\sup_{x \in \mathbb{R}} |\mathbb{P}(n^{-1/6}\widetilde{\xi}_n \leq x) - \mathbb{P}(\zeta_n \leq x)| \leq \kappa_3 n^{-1/2}$$

where $\zeta_n$ is obtained upon replacing $\xi_n$ with a normal random variable of zero mean and variance $n\vec{u}_n^\dagger \mathbb{Q}(0, \rho_n)\vec{u}_n$ [for the positive definite initial condition $\mathbb{Q}(0,\rho)$ of the ODE (2.8) at $\rho = \rho_n$]. With $\mathbb{E}\zeta_n = n^{1/2}y_1^*(\tau_n) - \frac{\widetilde{F}}{2}n^{2\beta-3/2}$, it follows from (4.9) that $|\mathbb{E}\zeta_n - \mathbb{E}\xi(r)| \leq Cn^{3\beta-5/2}$. Similarly, it follows from (4.10) that for some $C$ finite and all $n$,

$$|\text{Var}(\zeta_n) - \text{Var}(\xi(r))| = |(\mathbb{Q}_{\tau_n})_{11} + \widetilde{G}n^{\beta-1} - Q_{11}(\theta_c, \rho_c)|$$
$$\leq Cn^{\beta-1} \leq Cn^{3\beta-5/2}.$$

With $\text{Var}(\xi(r)) > 0$ independent of $n$, our thesis (6.19) easily follows from these bounds on the difference in the mean and variance of the normal random variables $\zeta_n$ and $\xi(r)$.

DEPARTMENTS OF STATISTICS AND
MATHEMATICS
STANFORD UNIVERSITY
STANFORD, CALIFORNIA 94305
USA
E-MAIL: amir@math.stanford.edu

DEPARTMENTS OF ELECTRICAL ENGINEERING
AND STATISTICS
STANFORD UNIVERSITY
STANFORD, CALIFORNIA 94305
USA
E-MAIL: montanari@stanford.edu